\documentstyle[11pt,leqno,amscd,amssymb,pstricks,latexsym,amsbsy,xypic,mathrsfs,verbatim]{amsart}

\input xy
\xyoption{all}

\newcommand{\scr}[1]{\mathscr #1}
\newtheorem{Thm}{Theorem}[section]
\newtheorem{Lem}[Thm]{Lemma}
\newtheorem{Cor}[Thm]{Corollary}
\newtheorem{Prop}[Thm]{Proposition}

 \theoremstyle{definition}

\newtheorem{Rem}[Thm]{Remark}
\newtheorem{Rems}[Thm]{Remarks}

\numberwithin{equation}{section}

\def\mod{{\rm mod\,}}
\def\ind{{\rm ind}}

\def\Hom{{\rm Hom}}
\def\End{{\rm End}}
\def\Ext{{\rm Ext}}
\def\Aut{{\rm Aut}}

\def\Ker{{\rm Ker}}
\def\op{{\rm op}}

\def\co{{\mathcal O}}
\def\vx{{\vec{x}}}
\def\vy{{\vec{y}}}
\def\vc{{\vec{c}}}

\def\coh{{\rm coh\,}}
\def\vect{{\rm vect}}

\def\bbZ{{\mathbb Z}}
\def\bbN{{\mathbb N}}

\def\bbC{{\mathbb C}}
\def\bbX{{\mathbb X}}
\def\bbL{{\mathbb L}}
\def\bbP{{\mathbb P}}
\def\bbH{{\mathbb H}}

\def\cR{{\mathcal R}}

\def\cL{{\mathcal L}}
\def\cD{{\mathcal D}}
\def\cE{{\mathcal E}}

\def\fL{{\mathfrak L}}
\def\fg{{\mathfrak g}}

\def\fI{{\mathfrak I}}
\def\fkF{{\mathfrak F}}

\def\ad{{\rm ad}}
\def\az{\alpha}
\def\bz{\beta}

\def\dz{\delta}

\def\lz{\lambda}
\def\llz{\Lambda}

\def\oz{\omega}
\def\vez{\varepsilon}

\def\dz{\delta}

\def\lra{\longrightarrow}

\def\lan{\langle}
\def\ran{\rangle}

\def\bbF{{\mathbb F}}
\def\refl{{\rm r}}
\def\Ref{{\rm R}}
\def\scrS{{\scr S}}
\def\vez{{\varepsilon}}
\def\lmto{{\longmapsto}}
\def\rmx{{\rm x}}
\def\bfx{{\bf x}}
\def\bfy{{\bf y}}
\def\bfv{{\boldsymbol v}}
\def\re{{\rm re}}
\def\im{{\rm im}}
\def\b1{{\bf 1}}
\def\whI{{\widehat I}}
\def\wtilde{\widetilde}
\def\bfU{{\bf U}}
\def\what{\widehat}
\def\bfPsi{{\boldsymbol \Psi}}
\def\scrA{{\scr A}}
\def\calU{{\cal U}}

\begin{document}

\title[Mutations and weighted projective lines]
{Applications of mutations in the derived \\ categories of weighted projective lines \\ to Lie and quantum algebras}
\author{Bangming Deng, Shiquan Ruan and Jie Xiao}
\address{Yau Mathematical Sciences Center, Tsinghua University,
Beijing 100084, China.} \email{bmdeng@@math.tsinghua.edu.cn}
\address{School of Mathematical Sciences, Xiamen University,
Xiamen 361005, China.} \email{sqruan@@xmu.edu.cn}
\address{Department of Mathematical Sciences, Tsinghua University,
Beijing 100084, China.} \email{jxiao@@math.tsinghua.edu.cn}

\begin{abstract} Let $\coh\bbX$ be the category of coherent sheaves over a weighted
projective line $\bbX$ and let $D^b(\coh\bbX)$ be its bounded derived category. The
present paper focuses on the study of the right and left mutation functors arising
in $D^b(\coh\bbX)$ attached to certain line bundles. As applications, we first show
that these mutation functors give rise to simple reflections for the Weyl group of
the star shaped quiver $Q$ associated with $\bbX$. By further dealing with the Ringel--Hall
algebra of $\bbX$, we show that these functors provide a realization for
Tits' automorphisms of the Kac--Moody algebra $\frak{g}_Q$ associated with $Q$, as well as for Lusztig's
symmetries of the quantum enveloping algebra of ${\frak g}_Q$.
 \end{abstract}

 \thanks{Supported by the National Natural Science Foundation of China.}

%\date{\today}



\maketitle

\section{Introduction}

The weighted projective lines and their coherent sheaf categories
were introduced by Geigel and Lenzing \cite{GL87} whose aim was to provide
geometric counterparts to canonical algebras and their module categories
studied by Ringel \cite{R84}. More precisely, a weighted projective line
over a field $\bbF$ is characterized by a weight sequence ${\bf p}=(p_1,\ldots,p_t)$
of positive integers and a sequence ${\boldsymbol\lambda}=(\lz_1,\ldots,\lz_t)$ of
pairwise distinct points in the projective line $\bbP_\bbF^1$. Geigle and Lenzing
\cite{GL87} proved that the category $\coh\bbX$ of coherent sheaves over the weighted projective
line $\bbX=\bbX({\bf p},{\boldsymbol\lambda})$ associated with $({\bf p},{\boldsymbol\lambda})$
is derived equivalent to the category of finite dimensional modules over the canonical
algebra $\llz=\llz({\bf p},{\boldsymbol\lambda})$ of the same type. An analogue of Kac's Theorem
which describes the dimension types of indecomposable coherent sheaves over $\bbX$ was obtained by
Crawley-Boevey \cite{CB10}. The Weyl groups and Artin groups associated to $\bbX$ have been
studied by Shiraishi, Takahashi and Wada in \cite{STW}.

Motivated by the work of Kapranov \cite{Kap}, Schiffimann \cite{Sch04} studied
the Hall algebra $H(\bbX)$ of a weighted projective line $\bbX$ over a finite field $\bbF$
and showed that it provides a realization of the quantum enveloping
algebra of a certain nilpotent subalgebra of the affinization of the corresponding
Kac--Moody algebra. In particular, he obtained a geometric construction of the quantum enveloping
algebra of the elliptic Lie algebra of type $D^{(1,1)}_4$, $E^{(1,1)}_6$, $E^{(1,1)}_7$ or
$E^{(1,1)}_8$ in the sense of Saito--Yoshii (\cite{Saito,SY}). Independently, by considering
the root categories of tubular algebras
(a special class of canonical algebras), Lin and Peng \cite{LinP} obtained a realization of these
elliptic Lie algebras. By further dealing with properties of the reduced Drinfeld double
$\cD C(\coh\bbX)$ of the composition subalgebra of $H(\bbX)$, Burban and Schiffmann \cite{BurSch} obtained
a new realization of the quantum enveloping algebras of affine Lie algebras of simply-laced
types as well as some new embeddings between them, and derived new results on the structure of
the quantized enveloping algebras of the elliptic Lie algebras. Based on Schiffmann's work \cite{Sch04},
Dou, Jiang and Xiao \cite{DJX} explicitly constructed a collection of generators of $\cD C(\coh\bbX)$
and verified that they satisfy all the Drinfeld relations of the corresponding quantum loop algebra.

The notion of reflection functors for the categories of representations of quivers
was introduced by Bernstein, Gelfand and Ponomarev \cite{BGP} whose aim was to obtain
a simple and elegant proof of Gabriel's theorem. Indeed, in the Grothedieck group level, BGP-reflection
functors give rise to simple reflections in the associated Weyl groups. Further,
it was shown in \cite{SV,XY} that BGP-reflection functors induce isomorphisms of double Ringel--Hall
algebras of quivers with distinct orientations, and their restrictions to the associated quantum enveloping algebras coincide
with Lusztig's symmetries. In some sense, BGP-reflection functors provide a categorification
of Lusztig's symmetries. However, this construction does not work in general since BGP-reflection functors can only be
defined for sinks or sources of quivers. One of the motivations in the present paper is to
generalize this construction to arbitrary vertices for the star shaped quivers associated with
weighted projective lines.

In order to study interrelations of various categories with the module categories
over finite dimensional algebras, Bondal \cite{Bondal} introduced mutations in a
triangulated category in terms of exceptional collections, which generalize reflection functors
defined in \cite{BGP}, as well as tilting modules in \cite{BB76}. Let $D^b(\coh\bbX)$ be the bounded
derived category of $\coh\bbX$ over a weighted projective line $\bbX$. Each line bundle
$\co(\vec{x})\in\coh\bbX$ becomes an exceptional object in $D^b(\coh\bbX)$.
By applying Bondal's construction to $D^b(\coh\bbX)$ in terms of $\co(\vec{x})$,
we obtain mutually inverse triangulated functors between the bounded derived categories
$D^b(\co(\vx)\!^{\bot})$ and $D^b({}^{\bot}\co(\vx))$, which are called the right and left
mutation functors, respectively. Here $\co(\vx)\!^{\bot}$ and ${}^{\bot}\co(\vx)$
are the right and left perpendicular subcategories of $\co(\vx)$ in $\coh\bbX$, respectively.

Now let $Q=(I,Q_1)$ be the star shaped quiver attached to $\bbX$, where $I$ and
$Q_1$ are the vertex set and arrow set, respectively. The main purpose of this
paper is to apply these mutation functors to obtain
a realization of simple reflections, Tits' automorphisms and Lusztig's symmetries associated
with $Q$. We emphasize that this approach works for an arbitrary vertex in $Q$ (which is neither
a sink nor a source) and, moreover, our construction provides new simple root bases in the root
system of $Q$ which do not arise from simple representations of $Q$.
More precisely, since both $D^b(\co(\vx)\!^{\bot})$ and $D^b({}^{\bot}\co(\vx))$ are equivalent
to the bounded derived category of the module category of the path algebra of $Q$, their Grothendieck
groups are isomorphic to the free abelian group $\bbZ I$ with basis $I$. By identifying $\bbZ I$
with a quotient group of the Grothendieck group of $\coh\bbX$, the mutation
functors induced by $\co(\vx)$ give rise to invertible linear transformations on $\bbZ I$.
We first show that such linear transformations attached to a special family
of line bundles $\co(\vx)$ realize all simple reflections in the Weyl group $W(Q)$ of
the quiver $Q$ (or the Kac--Moody algebra $\fg_Q$ of $Q$). Secondly, following the idea in \cite{PX00},
there are (generic) Hall Lie algebras associated with the root categories of $D^b(\co(\vx)\!^{\bot})$
and $D^b({}^{\bot}\co(\vx))$ both of which are isomorphic to the Hall Lie algebra of $Q$. In particular, their
composition Lie subalgebras are all isomorphic to the Kac--Moody algebra $\fg_Q$. We then
show that the mutation functors mentioned above induce automorphisms of $\fg_Q$ which
coincide with the Tits' automorphisms of $\fg_Q$, up to sign. Finally, be working with the
double Ringel--Hall algebra of $\bbX$, we prove that these mutation functors also
provide a realization for Lusztig's symmetries of the quantum enveloping algebra ${\bf U}_v(\fg_Q)$
of ${\frak g}_Q$.

The paper is organized as follows. Section 2 mainly recalls the notions of a weighted projective
line $\bbX$ and its category $\coh\bbX$ of coherent sheaves. In section 3 we show that
the mutation functors arising in the derived category $D^b(\coh\bbX)$ of $\coh\bbX$
give rise to simple reflections in the Weyl group of the star shaped quiver $Q$ associated with $\bbX$.
 Section 4 treats the Ringel--Hall Lie algebra $\fL(\cR_{\bbX})$ of the
root category $\cR_{\bbX}$ of $\coh\bbX$ and introduces its generic composition Lie algebra
$\fL_c(\cR_{\bbX})$. It is then shown that the mutation functors
induce isomorphisms between certain subalgebras of $\fL_c(\cR_{\bbX})$ which
 give rise to Tits' automorphisms of Kac--Moody algebra $\fg_Q$ associated with $Q$.
In section 5 we deal with the double Ringel--Hall algebra of $\coh\bbX$ and some
subalgebras of it and show that the mutation functors provide a realization of
Lusztig's symmetries of the quantum enveloping algebra $\bfU_v(\fg_Q)$ of $\fg_Q$.
Some basic calculations and results concerning the linear quiver, which are needed
in sections 3 and 5, are collected in the appendix A.

\bigskip

\section{Preliminaries}

In this section we recall some basic notions needed in the paper, including the
category of coherent sheaves over a weighted projective line, the root
system, Kac--Moody algebra and quantum enveloping algebra associated with a quiver.
We refer to \cite{GL87,Kac90, L93,Jan,CK} for details.

\subsection{}
Following \cite{GL87}, a \emph{weighted projective line} $\bbX=\bbX_\bbF$ over a field $\bbF$ is given by
a weight sequence ${\bf p}=(p_{1},\ldots, p_{t})$ of positive integers, and a collection
${\boldsymbol\lambda}=(\lambda_{1},,\ldots, \lambda_{t}) $ of
distinct closed points (of degree $1$) in the projective line $\bbP^{1}_\bbF$ which can be
normalized as $\lambda_{1}=\infty, \lambda_{2}=0, \lambda_{3}=1$.
More precisely, let $\bbL=\bbL(\bf p)$ be the rank one abelian group with
generators $\vec{x}_{1}, \ldots, \vec{x}_{t}$ and the relations
\[ p_{1}\vx_1=\cdots=p_{t}\vx_t=:\vec{c},\]
where $\vec{c}$ is called the \emph{canonical element} of
$\mathbb{L}$. Each element $\vx\in\bbL$ has the \emph{normal form} $\vx=\sum\limits_{i=1}^t l_i\vx_i+l\vc$ with $0\leq l_i\leq p_i-1$ and $l\in\bbZ$. Denote by $S$ the commutative algebra
$$S=S({\bf p},{\boldsymbol\lambda})=\bbF[X_{1},\cdots, X_{t}]/{\frak a}:
= \bbF[{\rm x}_{1}, \ldots, {\rm x}_{t}],$$
where ${\frak a}=(f_{3},\ldots,f_{t})$ is the ideal generated by
$f_{i}=X_{i}^{p_{i}}-X_{2}^{p_{2}}+\lambda_{i}X_{1}^{p_{1}}$ for
$3\leq i\leq t$. Then $S$ is $\mathbb{L}$-graded by setting
$$\mbox{deg}(\rmx_{i})=\vx_i\; \text{ for $1\leq i\leq t$.}$$
Finally, the weighted projective line associated with $\bf p$ and $\boldsymbol\lz$
is defined to be
$$\bbX=\bbX_\bbF=\rm{Spec}^{\bbL}S,$$
 the spectrum of $\bbL$-graded homogeneous prime ideals of $S$.

The category of coherent sheaves on $\bbX$ is defined to be the quotient category
$$\coh\bbX={\rm mod}^{\mathbb{L}}S/\mbox{mod}_{0}^{\mathbb{L}}S,$$
where ${\rm mod}^{\mathbb{L}}S$ is the category of finitely generated $\mathbb{L}$-graded
$S$-modules, while $\mbox{mod}_{0}^{\mathbb{L}}S$ is the Serre subcategory of $\mathbb{L}$-graded
$S$-modules of finite length. The grading shift gives twists $E(\vec{x})$
for every sheaf $E$ and $\vec{x}\in\bbL$.

Moreover, $\coh\bbX$ is a hereditary abelian category with Serre
duality of the form
$$D\Ext^1(X, Y)\cong\Hom(Y, X(\vec{\omega})),$$
where $D=\Hom_\bbF(-,\bbF)$, and $\vec{\omega}:=(t-2)\vec{c}-\sum_{i=1}^{t}\vec{x}_{i}\in\bbL$ is
called the \emph{dualizing element}. This implies
the existence of almost split sequences in $\coh\bbX$ with the
Auslander--Reiten translation $\tau$ given by the grading shift with
$\vec{\omega}$.

It is known that $\coh\bbX$ admits a splitting torsion pair $({\rm coh}_0\bbX,\vect{\bbX})$,
where ${\rm coh}_0\bbX$ and $\vect{\bbX}$ are full subcategories
of torsion sheaves and vector bundles, respectively.
The free module $S$ yields a structure sheaf $\co\in \vect{\bbX}$, and each
object in $\vect{\bbX}$ has a finite filtration by line bundles, that is, sheaves of the form
$\co(\vec{x})$. The subcategory ${\rm coh}_0\bbX$ admits ordinary simple
sheaves $S_x$ for each $x\in\bbH_\bbF:=\bbP^{1}_\bbF\backslash\{\lz_1,\ldots,\lz_t\}$
and exceptional simple sheaves $S_{ij}$ for $1\leq i\leq t$ and $0\leq j\leq p_i-1$ satisfying that for each $r\in\bbZ$
$$\Hom(\co(r\vec{c}), S_{ij})\not=0 \Longleftrightarrow j=0,$$
and that the nonzero extensions between these simple sheaves are given by
$$\Ext^1(S_x, S_x)\cong\bbF(x),\; \Ext^1(S_{ij}, S_{ik})\cong\bbF\;\text{ for $k\equiv j-1$ (mod\, $p_i$),}$$
 where $\bbF(x)$ denotes the finite extension of $\bbF$ with $[\bbF(x):\bbF]$ the degree of $x$.
For each simple sheaf $S$ and $n\geq 1$, there is a unique sheaf $S^{(n)}$
with length $n$ and top $S$, which is uniserial. Indeed, the sheaves $S^{(n)}$ form a complete
set of indecomposable objects in ${\rm coh}_0\mbox{-}\bbX$. For convenience,
we also use the notation $S_{ij}$ for $1\leq i\leq t$ and $j\in\bbZ$ to denote
the simple sheaf $S_{il}$ with $j\equiv l$ (mod $p_i$) and $0\leq l\leq p_i-1$.

The Grothendieck group $K_0(\coh\bbX)$ of $\coh\bbX$ is a free abelian group
with a basis $[\co(\vec{x})]$ with $0\leq\vec{x}\leq \vec{c}$.
The Euler form associated with $\bbX$ is the
bilinear form on $K_0(\coh\bbX)$ defined by
$$\langle [X],[Y]\rangle=\dim_\bbF\Hom(X,Y)-\dim_\bbF\Ext^{1}(X,Y)$$
for any $X,Y\in\coh\bbX$. Its symmetrization is defined by
$$([X],[Y])=\langle [X],[Y]\rangle +\langle [Y],[X]\rangle.$$
Recall that $\delta:=[\co(\vc)]-[\co]$ satisfies that $(\delta, [X])=0$ for any $X\in\coh\bbX$.

\subsection{} Let $Q=(I=Q_0,Q_1)$ be a quiver with vertex set $I=Q_0$ and arrow set $Q_1$.
For an arrow $\rho$ in $Q$, we denote by $t\rho$ and $h\rho$ the tail and head of
$\rho$, respectively. We always assume that $Q$ is acyclic, i.e., $Q$ contains no oriented cycles.
Let $\bbZ I$ be the free abelian group with basis $I$ whose elements will
be often written as $\bfx=(x_i)$. For each $i\in I$, set $\az_i=(\delta_{ij})_{j\in I}$,
called a {\it simple root}. By definition, the Euler form
$\langle-,-\rangle: \bbZ I\times\bbZ I\rightarrow\bbZ$ associated with $Q$ is defined by
$$\langle \bfx,\bfy\rangle=\sum_{i\in I}x_iy_i-\sum_{\rho\in Q_1}x_{t\rho}y_{h\rho},$$
 where $\bfx=(x_i),\bfy=(y_i)\in\bbZ I$. Its symmetrization
$$(\bfx,\bfy):=\langle \bfx,\bfy\rangle+\langle \bfy,\bfx\rangle$$
is called the {\it symmetric} Euler form of $Q$. It is easily checked that
$$(\bfx,\bfy)={\bfx}^t C_Q{\bfy},$$
where elements in $\bbZ I$ are considered as column vectors,
${\bfx}^t$ denotes the transpose of $\bfx$, and $C_Q=(c_{ij})_{i,j\in I}$ is
the generalized Cartan matrix attached to $Q$ with entries given by
$$c_{ij}=\left\{ \begin{array}{ll}  2, \;\;  &\mbox{if}\;\; i=j; \\
    -|\{\mbox{arrows between $i$ and $j$}\}|, \;\; &\mbox{if}\;\; i\ne j.
        \end{array}\right.$$
For each vertex $i\in I$, we define an element $r_i\in\Aut(\bbZ I)$ by
\begin{equation}\label{reflection}
\refl_i(\mu)=\mu-(\mu,\alpha_i)\alpha_i,\;\forall\,\mu\in\bbZ I,
\end{equation}
or, simply by
$\refl_i(\alpha_j)=\alpha_j-c_{ij}\alpha_i,\;\forall\, j\in I.$
 Since there are no loops at the vertex $i$, then $\refl_i(\alpha_i)=-\alpha_i$
and $\refl_i$ fixes the set $\{\mu\in\bbZ I\;|\;(\mu,\alpha_i)=0\}$. Moreover,
$\refl_i$ has order 2. We call $\refl_i$ a {\it simple reflection}. It is easy to see that,
for each $i\in I$, we have
$$(\refl_i(\mu),\refl_i(\nu))=(\mu,\nu),\;\forall\,\mu,\nu\in\bbZ I.$$
 The subgroup of $\Aut(\bbZ I)$ generated by the simple reflections $\refl_i$
is called the {\it Weyl group} of $Q$, denoted by $W=W(Q)$.

With $Q$ we can associate a Kac--Moody Lie algebra $\fg=\fg_Q$ over the complex
 field $\bbC$ which is generated by $e_i, f_i, h_i$ ($i\in I$) with defining relations:
$$\left\{\aligned
{} &[h_i,h_j]=0,\; [e_i, f_j]=\delta_{ij} h_i,\\
&[h_i, e_j]=c_{ij} e_j, \;  [h_i, f_j]=-c_{ij}f_j,\\
 &(\ad e_i)^{1-c_{ij}}(e_j)=0,\; (\ad f_i)^{1-c_{ij}}(f_j)=0\; (\forall\, i\neq j).
\endaligned\right.$$
 We remark that $\fg$ is in fact the derived algebra of the Kac--Moody algebra
 associated with $C_Q$ defined in \cite{Kac90}. It is well known that
 $\fg=\oplus_{\az\in\bbZ I}\fg_\az$ is graded by $\bbZ I$ with $\deg e_i=\alpha_i,
\deg f_i=-\alpha_i$ and $\deg h_i=0$ for all $i\in I$, and the associated root
system is defined to be the set
$$\Delta=\Delta(Q)=\{0\neq \alpha\in\bbZ I\mid \fg_{\alpha}\neq 0\}$$
 which is divided into the disjoint union of two subsets $\Delta^{\re}$ (of real roots)
 and $\Delta^{\im}$ (of imaginary roots). More generally, we can
define a reflection $\refl_\az\in W(Q)$ with respect to each real root $\az\in \Delta(Q)$ by
$$\refl_\az(\mu)=\mu-(\mu,\alpha)\alpha,\;\forall\,\mu\in\bbZ I.$$

Furthermore, we have the loop algebra $\cL\fg$ of $\fg$ which is a complex
Lie algebra generated by $e_{ir},f_{ir}, h_{ir}$ ($i\in I,r\in\bbZ$) and $c$
subject to some relations; see, for example, \cite[\S1.3]{Sch04}. Then $\fg$ is identified
with the Lie subalgebra of $\cL\fg$ generated by $e_{i0}$, $f_{i0}$ and $h_{i0}$ for $i\in I$.

Now let $\bbC(\bfv)$ be the field of rational functions in an indeterminate $\bfv$.
By definition, the quantum enveloping algebra ${\bf U}_\bfv(\fg_Q)$ of $\fg_Q$ is
a $\bbC(\bfv)$-algebra generated by $E_i,F_i,K_i^{\pm 1}$ ($i\in I$) with the following relations:
$$\left\{\aligned
{}& K_iK_j=K_jK_i,\; K_iK_i^{-1}=1,\;\;\forall\, i,j\in I;\\
& K_i E_j=\bfv^{c_{ij}}E_jK_i,\; K_iF_j=\bfv^{-c_{ij}}F_jK_i,\;\;\forall\,i,j\in I;\\
&E_iF_j-F_jE_i=\dz_{ij}\frac{K_i-K_{-i}}{\bfv-\bfv^{-1}},\;\;\forall\,i,j\in I;\\
&\sum_{s=0}^{1-c_{ij}}(-1)^s E_i^{(s)} E_j E_i^{(1-c_{ij}-s)}=0,\;\;\forall\,i\not=j\in I;\\
&\sum_{s=0}^{1-c_{ij}}(-1)^s F_i^{(s)} F_j F_i^{(1-c_{ij}-s)}=0,\;\;\forall\,i\not=j\in I,
\endaligned\right.$$
 where $E_i^{(s)}=E_i^s/[s]!$ and $F_i^{(s)}=F_i^s/[s]!$ are divided powers. Here
$$[s]!=\prod_{t=1}^s [t]\;\text{ with }\; [t]=\frac{\bfv^t-\bfv^{-t}}{\bfv-\bfv^{-1}}.$$
Similarly, we can also define the quantum enveloping algebra ${\bf U}_v(\fg_Q)$ of $\fg_Q$
over $\bbC$ by specializing $\bfv$ to a non-root of unity $v\in\bbC$. We are mainly concerned
with $v=\sqrt{q}$ for a prime power $q$.

\subsection{} Given the weighted projective line $\bbX$ associated with the pair
$({\bf p},{\boldsymbol\lambda})$, there is an associated star shaped quiver $Q=Q_{\bbX}$
whose vertex set $I$ consists of vertices $\star$ and $\oz_{ij}$ for
$1\leq i\leq t$ and $1\leq j\leq p_i-1$ and whose arrows are given as follows:
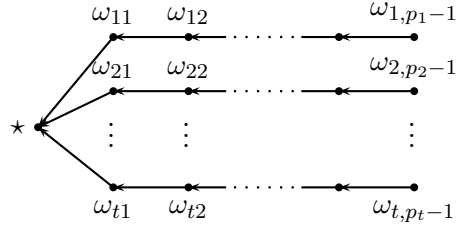
\begin{center}
\begin{pspicture}(-1,0.3)(8,3.6)
\psset{xunit=1cm,yunit=.8cm} \psdot*[dotsize=3pt](0,2)
\psdot*[dotsize=3pt](1,3.5) \psdot*[dotsize=3pt](2,3.5)
\psdot*[dotsize=3pt](4,3.5) \psdot*[dotsize=3pt](5,3.5)
\psdot*[dotsize=3pt](1,1) \psdot*[dotsize=3pt](2,1)
\psdot*[dotsize=3pt](4,1) \psdot*[dotsize=3pt](5,1)
\psdot*[dotsize=3pt](1,2.6) \psdot*[dotsize=3pt](2,2.6)
\psdot*[dotsize=3pt](4,2.6) \psdot*[dotsize=3pt](5,2.6)
\psline{<-}(0,2)(1,3.5) \psline{<-}(1,3.5)(2,3.5) \psline{<-}(2,3.5)(2.5,3.5)
\psline[linestyle=dotted,linewidth=1pt](2.6,3.5)(3.4,3.5)
\psline{<-}(4,3.5)(5,3.5)  \psline{-}(3.5,3.5)(4,3.5)
\psline{<-}(0,2)(1,2.6) \psline{<-}(1,2.6)(2,2.6) \psline{<-}(2,2.6)(2.5,2.6)
\psline[linestyle=dotted,linewidth=1pt](2.6,2.6)(3.4,2.6)
\psline{<-}(4,2.6)(5,2.6)  \psline{-}(3.5,2.6)(4,2.6)
\psline{<-}(0,2)(1,1) \psline{<-}(1,1)(2,1) \psline{<-}(2,1)(2.5,1)
\psline[linestyle=dotted,linewidth=1pt](2.6,1)(3.4,1)
\psline{-}(3.5,1)(4,1)\psline{<-}(4,1)(5,1)
\uput[l](1.2,2){$\vdots$} \uput[l](2.2,2){$\vdots$} \uput[l](5.2,2){$\vdots$}
\uput[l](0,2){$\star$} \uput[u](1,3.5){$\oz_{11}$}
\uput[u](2,3.5){$\oz_{12}$} \uput[u](5,3.5){$\oz_{1,p_1-1}$}
\uput[u](1,2.6){$\oz_{21}$} \uput[u](2,2.6){$\oz_{22}$} \uput[u](5,2.6){$\oz_{2,p_2-1}$}
\uput[d](1,1){$\oz_{t1}$} \uput[d](2,1){$\oz_{t2}$} \uput[d](5,1){$\oz_{t,p_t-1}$}
\end{pspicture}
\end{center}
 For convenience and simplicity, we sometimes use the notation $\star=\oz_{i0}$ for
 each $1\leq i\leq t$ and put
$$\widehat I=\{(i,j)\mid 1\leq i\leq t, 1\leq j\leq p_i-1\}.$$
Thus, for the star shaped quiver $Q$, we have the associated Weyl group $W(Q)$, Kac--Moody
Lie algebra $\fg_Q$ and quantum enveloping algebra ${\bf U}_\bfv(\fg_Q)$ of $Q$
which are the central objects studied in the present paper.
Note that by \cite[\S1]{CB10}, the loop algebra $\cL\fg_Q$ of $\fg_Q$ is $\what\Gamma$-graded with
$$\deg e_{\oz r}=\az_\oz+r\dz,\; \deg f_{\oz r}=-\az_\oz+r\dz,\;
\deg h_{\oz r}=r\dz,\; \deg c=0,$$
where $\what\Gamma=\bbZ I\oplus\bbZ\dz$, $\oz\in I$, and $r\in\bbZ$. Moreover, $\cL\fg_Q$
admits a root space decomposition
$$\cL\fg_Q=(\cL\fg_Q)_0\oplus \bigoplus_{\zeta\in\widehat\Delta} (\cL\fg_Q)_\zeta\;\,
\text{with $(\cL\fg_Q)_0=(\fg_Q)_0\oplus\bbC c$,}$$
  where
 $$\widehat\Delta=\{\az+r\dz\mid \az\in\Delta(Q),r\in\bbZ\}\cup\{r\dz\mid 0\not=r\in\bbZ\}.$$
\medskip

{\it Throughout the paper, let $\bbX$ be a weighted projective line
of weight type ${\bf p}=(p_1, \ldots, p_t)$ over a field $\bbF$ and $Q=Q_\bbX=(I,Q_1)$ be the associated
star shaped quiver given above. Given an abelian category $\scr A$, we denote by $D^b(\scr A)$
the bounded derived category of $\scr A$ with the shift functor $[1]$ and by $\cR_{\scr A}=D^b({\scr A})/[2]$ its root
category; see \cite{Ha} or \cite{PX00} for the definitions.
By $\ind \cR_{\scr A}$ we denote the set of representatives of the
isoclasses of indecomposable objects in $\cR_{\scr A}$. The Grothendieck group of $\scr A$ is denoted by $K_0(\scr A)$.
Each object $X\in\scr A$ gives a
stalk complex $C^\bullet=(C^i)$ with $C^0=X$ and $C^i=0$ for all $i\not=0$, which is still
denoted by $X$.}

\bigskip

\section{Mutation functors and simple reflections}

In this section we apply the construction of right and left mutation functors defined
in \cite{Bondal} to the bounded derived category $D^b(\coh\bbX)$ of a weighted projective line
$\bbX$. The main aim is to show that the mutation functors attached to certain line bundles
over $\bbX$ give rise to all simple
reflections in the Weyl group $W(Q)$ of the star shaped quiver $Q=Q_\bbX$ associated with $\bbX$.
In other words, these mutation functors provide a catgorification of simple reflections in $W(Q)$.

\subsection{}  For each fixed $\vx\in \bbL$, set
$$\co(\vx)\!^{\bot}=\{X\in \coh\bbX\mid \Hom(\co(\vx), X)=0=\Ext^1(\co(\vx), X)\}$$
and
$${}^{\bot}\!\co(\vx)=\{X\in \coh\bbX\mid \Hom(X, \co(\vx))=0=\Ext^1(X, \co(\vx))\}.$$
 By \cite{GL91}, both $\co(\vx)\!^{\bot}$ and ${}^{\bot}\!\co(\vx)$ are again abelian
hereditary categories. Moreover, we have
$$D^b(\co(\vx)\!^{\bot})=\{X^\bullet\in D^b(\coh\bbX)\mid \Hom(\co(\vx), X^\bullet[n])=0,\;\forall\,n\in\bbZ\},$$
 that is, $D^b(\co(\vx)\!^{\bot})$ coincides with the right perpendicular subcategory
 of $\co(\vx)$ in $D^b(\coh\bbX)$. Similarly, we obtain
$$D^b({}^{\bot}\!\co(\vx))=\{X^\bullet\in D^b(\coh\bbX)\mid \Hom(X^\bullet[n], \co(\vx))=0,\;\forall\,n\in\bbZ\}$$
 which is the left perpendicular subcategory of $\co(\vx)$ in $D^b(\coh\bbX)$.

For an $\bbF$-algebra $A$, let ${\rm mod}\,A$ be the category of finite dimensional left
$A$-modules. Also, by $\bbF Q$ we denote the path algebra of $Q$ over the field $\bbF$.
Then we have the following result.

\begin{Lem}\label{derived equivalent to star shaped quiver}
For each $\vx\in \bbL$, there are triangle equivalences
$${\frak F}_\vx: D^b(\co(\vx)\!^{\bot})\lra D^b({\rm mod}\, \bbF Q) \;\text{ and }\;
{}_\vx{\frak F}: D^b({}^{\bot}\!\co(\vx))\lra D^b({\rm mod}\, \bbF Q^{\op}),$$
 where $Q^{\op}$ denotes the opposite quiver of $Q$.
\end{Lem}

\begin{pf} By \cite[Prop.~14]{GL87}, $T_\vx=\bigoplus_{0\leq \vy\leq \vc}\co(\vx-\vy)$ is a
tilting bundle in $\coh\bbX$. Thus,
$$T_\vx^r:=\bigoplus_{0< \vy\leq \vc}\co(\vx-\vy)$$
is a tilting object in $\co(\vx)\!^{\bot}$. Then $\Hom(T_\vx^r,-)$ induces an equivalence
$$\fkF_\vx:={\bf R}\Hom(T_\vx^r,-): D^b(\co(\vx)\!^{\bot})\lra D^b(\mod\End(T_\vx^r)^\op).$$
By \cite[\S4]{GL87} again, $\End(T_\vx)^\op$ is the canonical algebra of the same type as
$\bbX$. This implies that $\End(T_\vx^r)^\op\cong \bbF Q$, which gives the first equivalence.

Similarly, the tilting object $T_\vx^l:=\bigoplus_{0<\vy\leq \vc}\co(\vx+\vy)$ in ${}^{\bot}\!\co(\vx)$
gives rise to the equivalence
$${}_\vx{\frak F}: D^b({}^{\bot}\!\co(\vx))\lra D^b({\rm mod}\, \bbF Q^{\op}).$$
\end{pf}

We now consider the case $\vx=\vc$. Then the functor
$\Hom(T_\vc^r,-): \co(\vc)\!^{\bot}\to \mod\bbF Q$
takes
$$\co\lmto P_\star,\;\co(j\vx_i)\lmto P_{\oz_{ij}},\;\forall\,(i,j)\in\whI,$$
  where $\whI=\{(i,j)\mid 1\leq i\leq t,1\leq j\leq p_i-1\}$, and $P_\oz$
denotes the indecomposable projective $\bbF Q$-module corresponding to $\oz\in I$.
For each $(i,j)\in\whI$, there is an exact sequence
$$0\lra\co((j-1)\vx_i)\lra \co(j\vx_i)\lra S_{ij}\lra 0$$
 in $\coh(\bbX)$. Applying $\Hom(T_\vc^r,-)$ gives an exact sequence
$$0\lra P_{\oz_{i,j-1}}\lra P_{\oz_{ij}}\lra S_{\oz_{ij}}\lra 0$$
 in $\mod\bbF Q$, where $S_{\oz_{ij}}$ denotes the simple $\bbF Q$-module corresponding
to the vertex $\oz_{ij}$. Consequently, the functor
${\frak F}_\vc: D^b(\co(\vc)\!^{\bot})\to D^b({\rm mod}\, \bbF Q)$ takes
\begin{equation}\label{value-F_c}
\co\lmto P_\star=S_\star,\;\co(j\vx_i)\lmto P_{\oz_{ij}},\;S_{ij}\lmto S_{\oz_{ij}},\;\forall\;(i,j)\in\whI.
\end{equation}

\subsection{} Since each line bundle $\co(\vx)$ is exceptional in $D^b(\coh\bbX)$, it follows
from \cite[Thm.~3.2]{Bondal} (see also \cite{M})
that the full subcategory $\lan\co(\vx)\ran$ of $D^b(\coh\bbX)$ generated by $\co(\vx)$ is admissible,
that is, the inclusion functor $\lan\co(\vx)\ran\hookrightarrow D^b(\coh\bbX)$ admits right and left
adjoint functors. Then the inclusion functors
$$j_{\ast,\vx}: D^b(\co(\vx)\!^{\bot})\lra D^b(\coh\bbX)\;\text{ and }\;
r_{\ast,\vx}: D^b({}^{\bot}\!\co(\vx))\lra D^b(\coh\bbX)$$
have respectively left and right adjoint functors, which we denote
by $j_{\vx}^{\ast}$ and $r^{!}_{\vx}$, respectively. Moreover, the restriction functor
$${\bf{R}}_{\vx}:=r^{!}_{\vx}|_{D^b(\co(\vx)\!^{\bot})}: D^b(\co(\vx)\!^{\bot})\lra D^b({}^{\bot}\!\co(\vx))$$
 is an equivalence whose quasi-inverse is given by
${\bf{L}}_{\vx}:=j^{\ast}_{\vx}|_{D^b({}^{\bot}\!\co(\vx))}$. The functors ${\bf{R}}_{\vx}$ and
${\bf{L}}_{\vx}$ are called the right mutation and left mutation functors, respectively,
which are determined by the following triangles
\begin{equation} \label{approaximation-1}
 {\bf{R}}_{\vx}(X)\lra X\lra \Hom^\bullet(X, \co(\vx))\otimes \co(\vx) \;\;\text{and}
\end{equation}
\begin{equation} \label{approaximation-2}
D\Hom^\bullet(\co(\vx), X)\otimes \co(\vx)\lra X \lra  {\bf{L}}_{\vx}(X).
\end{equation}

By using basic facts in $\coh\bbX$ given in \cite{GL87}, we obtain the following result.

\begin{Lem} \label{mutation formula for special objects} For each $\vx\in \bbL$ with normal form
$\vx=\sum\limits_{i=1}^tl_i\vx_i+l\vc$, we have
  \begin{itemize}
    \item[(1)]
         ${\bf{R}}_{\vx}(\co(\vx-k\vx_i))=
         \left\{\begin{array}{lll} S_{i, l_i}^{(k)}[-1], && 1\leq k\leq p_i-1;\\
              \co(\vx+\vc)[-1], && k=p_i,
             \end{array} \right.$
     %${\bf{R}}_{\vx}(\co(\vx-k\vx_i))=S_{i, l_i}^{(k)}[1]$ for $1\leq k\leq p_i-1$ and $\co(\vx+\vc)[1]$ for $k=p_i$;
    \item[(2)] for each $S_{ij}^{(k)}\in \co(\vx)\!^{\bot}$,
    ${\bf{R}}_{\vx}(S_{ij}^{(k)})=
    \left \{ \begin{array}{lll}
               \co(\vx+k\vx_i), && j=l_i+k;\\
                S_{ij}^{(k)}, && \text{else}.
             \end{array} \right.$
    %$\{\bf{R}}_{\vx}(S_{i, l_i}^{(k)})=\co(\vx+k\vx_i)$ for $j=l_i+k$ and $S_{i, l_i}^{(k)}$ else.
  \end{itemize}
\end{Lem}

\begin{pf}
  (1) In $D^b(\coh\bbX)$, there are triangles
  $$ \xymatrix{
  \co(\vx+\vc)[-1]\ar[r] & \co(\vx-\vc) \ar[r]^{(x_1^2, x_2^2)\quad} &  \co(\vx)\oplus\co(\vx)
  \ar[r]^{\quad(x_2^2, -x_1^2)^t}  & \co(\vx+\vc)} \;\;\text{ and}$$
  $$ \xymatrix{
  S_{i, l_i}^{(k)}[-1]\ar[r] & \co(\vx-k\vx_i) \ar[r]^{\;\quad x_i^k} &  \co(\vx)  \ar[r] & S_{i, l_i}^{(k)}},$$
where $1\leq k\leq p_i-1$. The assertion then follows from the isomorphisms
$$\Hom^\bullet(\co(\vx-\vc), \co(\vx))\otimes \co(\vx)\cong \co(\vx)\oplus\co(\vx)\;\;\text{ and}$$
$$\Hom^\bullet(\co(\vx-k\vx_i), \co(\vx))\otimes \co(\vx)\cong \co(\vx),\;\,\forall\, 1\leq k\leq p_i-1.$$

  (2) This follows from the (mutation) triangle:
  $$\xymatrix{\co(\vx+k\vx_i) \ar[r] & S_{i, l_i+k}^{(k)}\ar[r] & \co(\vx)[1] \ar[r]^{x_i^k\;\quad} &  \co(\vx+k\vx_i)[1].}$$
\end{pf}
%Hence the right mutation functor (see \cite{Bondal} for detailed definition)
%$${\bf{R}}_{\vx}: D^b(\coh\bbX)\to D^b(\coh\bbX); \quad X\mapsto R_{\vx};$$
%defined via the universal triangle
%\begin{equation}\label{approaximation}
% {\bf{R}}_{\vx}(X)\to X\to \Hom^.(X, \co(\vx))\otimes \co(\vx)
%\end{equation}
%induces a triangulated equivalence
%${\bf{R}}_{\vx}:\co(\vx)\!^{\bot}\cong {}^{\bot}\!\co(\vx),$ and hence a triangulated equivalence between their root categories ${\bf{R}}_{\vx}: \cR(\co(\vx)\!^{\bot})\cong\cR({}^{\bot}\!\co(\vx))$.
%which induces equivalences between their root categories $$\cR(\co(\vx)\!^{\bot})\cong %\cR(kQ)\cong \cR({}^{\bot}\!\co(\vx)).$$
%\begin{Lem}
%Let $\vx, \vy\in\bbL$. Then for any indecomposable object $X$ in $\cR( {}^{\bot}\!\co(\vx)) )$, there exists a unique indecomposable object $Y$ in $\cR( \co(\vy)\!^{\bot} )$, such that $\overline{[X]}=\overline{[Y]}$ in $K_0(\coh\bbX)/\bbZ\delta$.
%\end{Lem}
%
%\begin{Def} Keep the notation above. The corresponding
%$$\psi: \cR( {}^{\bot}\!\co(\vx)) ) \to \cR( \co(\vy)\!^{\bot} ); X\mapsto Y;$$ is called a Fourier transformation.
%\end{Def}
%
%Up to Fourier transformations, we can compose any two equivalence functors ${\bf{R}}_{\vx}$ and ${\bf{R}}_{\vy}$ on the root categories in the following way: $${\bf{R}}_{\vx}{\bf{R}}_{\vy}: \cR( \co(\vy)\!^{\bot}) \xrightarrow{{\bf{R}}_{\vy}} \cR( {}^{\bot}\!\co(\vy)) ) \xrightarrow{\psi} \cR( \co(\vx)\!^{\bot} )\xrightarrow{{\bf{R}}_{\vx}} \cR( {}^{\bot}\!\co(\vx)) ).$$

\subsection{} By \cite{GL87,Sch04}, the Grothendieck group $K_0(\coh\bbX)$ can be identified
with the free abelian group $\what\Gamma=\bbZ I\oplus\bbZ\dz$ via
$$[\co(r\vc)]=\az_\star+r\dz,\;[S_x]=\dz,\;[S_{ij}]=
\left \{ \begin{array}{lll}
               \az_{\oz_{ij}}, && \text{if $j\not=0$};\\
                \dz-\sum_{l=1}^{p_i-1}\az_{\oz_{il}}, && \text{if $j=0$},
             \end{array} \right.$$
 where $S_x$ is the simple sheaf associated with $x\in{\mathbb P}^1_\bbF\backslash\{\lz_1,\ldots,\lz_t\}$.

For each $\vx\in \bbL$, we denote respectively by $\phi_{\vx}^{r}$ and $\phi_{\vx}^{l}$ the
compositions of natural embeddings and projections
$$\aligned
&\phi_{\vx}^{r}: K_0(\co(\vx)\!^{\bot})\lra K_0(\coh\bbX)\lra K_0(\coh\bbX)/\bbZ\delta\;\;\text{and}\\
&\phi_{\vx}^{l}: K_0({}^{\bot}\!\co(\vx))\lra K_0(\coh\bbX)\lra K_0(\coh\bbX)/\bbZ\delta.
\endaligned$$
 Then the symmetric Euler form on $K_0(\coh\bbX)$ restricts to bilinear forms on both
$K_0(\co(\vx)\!^{\bot})$ and $K_0({}^{\bot}\!\co(\vx))$. We have the following fact.

\begin{Lem} For each fixed $\vx\in \bbL$, both $\phi_{\vx}^{r}$ and $\phi_{\vx}^{l}$ are isomorphisms.
\end{Lem}

\begin{pf} We only prove that $\phi_{\vx}^{r}$ is an isomorphism. The proof for
$\phi_{\vx}^{l}$ is similar.

By \cite{GL87}, $\{[\co(\vx-\vy)]\mid 0\leq \vy\leq \vc\}$ is a basis of $K_0(\coh\bbX)$. Thus,
$$\{[\co(\vx-\vy)]\mid 0< \vy\leq \vc\}$$
forms a basis of $K_0(\co(\vx)\!^{\bot})$. Then the equality $\delta=[\co(\vx)]-[\co(\vx-\vc)]$ implies
that
$$\{\overline{[\co(\vx-\vy)]}=[\co(\vx-\vy)]+\bbZ\delta\mid 0< \vy\leq \vc\}$$
forms a basis of $K_0(\coh\bbX)/\bbZ\delta$. By the definition, $\phi_{\vx}^{r}$ sends $[\co(\vx-\vy)]$ to
$\overline{[\co(\vx-\vy)]}$ for all $0< \vy\leq \vc$. Therefore, $\phi_{\vx}^{r}$ is an isomorphism.
\end{pf}

For simplicity, we write $\psi_{\vx}^{r}=:(\phi_{\vx}^{r})^{-1}$
and $\psi_{\vx}^{l}:=(\phi_{\vx}^{l})^{-1}$. By \cite{CB10}, there is an isometric isomorphism
$$\phi: K_0(\coh\bbX)/\bbZ\delta{\lra}\bbZ I$$
taking $\overline{[\co]}\mapsto\alpha_{\star}$ and $\overline{[S_{ij}]}\mapsto \alpha_{ij}$ for all
$(i,j)\in \whI$. Then, by \eqref{value-F_c}, $\phi$ is indeed induced by $\fkF_\vc$.

For each $\vx\in\bbL$, by viewing $\phi, \phi_{\vx}^{r}$ and $\phi_{\vx}^{l}$
as identities, we obtain two realizations of the root datum of $\fg=\fg_Q$:
\[\bigl(K_0(\co(\vx)\!^{\bot}), (-,-), \ind\,\cR(\co(\vx)\!^{\bot})\bigr)
\quad\text{and}\quad \bigl(K_0({}^{\bot}\!\co(\vx)), (-,-), \ind\,\cR({}^{\bot}\!\co(\vx))\bigr),\]
where $\cR(\co(\vx)\!^{\bot})$ (resp., $\cR({}^{\bot}\!\co(\vx))$) is the root category
of $\co(\vx)\!^{\bot}$ (resp., ${}^{\bot}\!\co(\vx))$.

Set
$$\scrS_+:=\{\co, S_{ij}\mid (i,j)\in\whI\},\;
\scrS_-=\scrS_+[1],\;\text{ and }\;\scrS=\scrS_+\cup \scrS_-.$$
Their images in the Grothendieck group $K_0(\coh\bbX)$ are denoted
by $[\scrS_+],[\scrS_-]$, and $[\scrS]$, respectively. In other words,
$$[\scrS_{\pm}]=\{\pm[\co], \pm[S_{ij}]\mid (i,j)\in\whI\}\subset K_0(\coh\bbX).$$

\begin{Prop} For each fixed $\vx\in\bbL$, the set $\psi_{\vx}^r(\overline{[\scrS_+]})$
forms a simple root basis of $(K_0(\co(\vx)\!^{\bot}), (-,-), \ind\,\cR(\co(\vx)\!^{\bot})$,
while the set $\psi_{\vx}^l(\overline{[\scrS_+]})$) forms a simple root basis of
$(K_0({}^{\bot}\!\co(\vx)), (-,-), \ind\,\cR({}^{\bot}\!\co(\vx)))$.
\end{Prop}

\begin{pf} In view of $\phi_{\vx}^{r}$ and $\phi_{\vx}^{l}$, it suffices to show
the following assertion: up to a multiple of $\delta$, the class of each indecomposable
coherent sheaf $X$ in $K_0(\coh\bbX)$ can be expressed as a linear combination of elements
in $[\scrS_+]$ with totally non-positive or totally non-negative coefficients.
Indeed, if $X$ is an indecomposable torsion sheaf, then $X=S_x^{(n)}$
for an ordinary simple sheaf $S_x$ or
$X=S_{ij}^{(k+np_i)}$ for an exceptional simple sheaf $S_{ij}$, where
$1\leq k\leq p_i-1, n\in\mathbb{N}$. In the former case, $\overline{[X]}=0$; while in the latter case, $$\overline{[X]}=\overline{[S_{ij}^{(k)}]}=\overline{[S_{ij}]}+\overline{[S_{i,j-1}]}+\cdots+\overline{[S_{i,j-k+1}]}.$$
Further, we have $\overline{[S_{i0}]}=-\sum\limits_{1\leq j\leq p_i-1}\overline{[S_{ij}]}$.
Consequently, the assertion holds for all indecomposable torsion sheaves.

Now let $X$ be an indecomposable vector bundle. Then there is a line bundle
filtration
$$0=X_0\subset X_1\subset X_2\subset \cdots \subset X_m=X,$$
where each factor $X_{k}/X_{k-1}$ is a line bundle for $1\leq k\leq m$.
But for line bundle $\co(\vx)$ with $\vx\in\bbL$ of the normal form $\vx=\sum_{1\leq i\leq t}l_i\vx_i+l\vc$,
we have the equality
$$[L]=[\co]+\sum_{1\leq i\leq t}\sum_{1\leq j\leq l_i}[S_{ij}]+l\delta.$$
 Therefore, the assertion holds for all indecomposable vector bundles.
\end{pf}

\begin{Rem} The proposition above
provides some new simple root bases which do not arise from simple representations of a quiver.
For example, let $\bbX$ be a weighted projective line of weight type ${\bf p}=(1, 3)$. Then $$T=\co(-\vx_1)\oplus\co\oplus\co(\vx_1)\oplus\co(2\vx_1)$$ is a tilting bundle in $D^{b}(\coh\bbX)$.
It is easy to see that the right perpendicular category
$D^b(\co(2\vx_1)^{\bot})$ is triangulated equivalent to the bounded derived category of finite
dimensional modules over the path algebra of type $A_3$. An easy calculation shows
that the indecomposable objects in $D^b(\co(2\vx_1)^{\bot})$ are given by the following coherent
sheaves
$$\co(-\vx_1), \co, \co(\vx_1), S_{10}, S_{11}, S_{11}^{(2)}$$
together with their shifts. Moreover, the Auslander--Reiten quiver
of  $D^b(\co(2\vx_1)^{\bot})$  has the following form:

$$ 
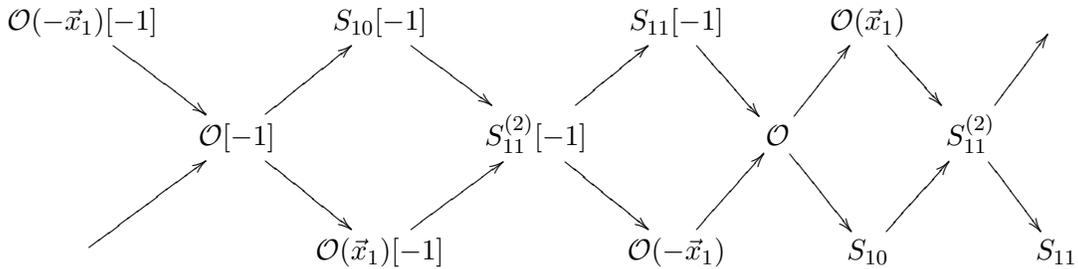
\xymatrix@C=0.8em{
 \co(-\vx_1)[-1]\ar@{->}[rd] && S_{10}[-1] \ar@{->}[rd] && S_{11}[-1] \ar@{->}[rd] &&\co(\vx_1) \ar@{->}[rd] && \\
 &\co[-1]\ar@{->}[rd]\ar@{->}[ru] && S_{11}^{(2)}[-1] \ar@{->}[rd]\ar@{->}[ru] && \co \ar@{->}[rd]\ar@{->}[ru] &&S_{11}^{(2)} \ar@{->}[rd]\ar@{->}[ru]\\
\ar@{->}[ru]&&\co(\vx_1)[-1]\ar@{->}[ru] && \co(-\vx_1)\ar@{->}[ru] && S_{10}\ar@{->}[ru] && S_{11}  }$$

\medskip

\noindent In this case, the isomorphism $\psi_{2\vx_1}^r: K_0(\coh\bbX)/\bbZ\delta\to K_0(\co(\vx)\!^{\bot})$
takes
$$\overline{[\co]}\lmto [\co],\; \overline{[S_{11}]}\lmto [S_{11}],\; \overline{[S_{12}]}\lmto [S_{11}^{(2)}[-1]].$$
This gives a simple root basis
$$\psi_{2\vx_1}^r(\overline{[\scrS_+]})=\{[\co],\; [S_{11}],\; [S_{11}^{(2)}[-1]]\}.$$
In fact, if we simply write $a=[\co], b=[S_{11}], c=[S_{11}^{(2)}[-1]]$, then each root
can be presented by the sum of $a,b$ and $c$ with totally non-positive or totally non-negative
coefficients, which is depicted as follows:
$$ 
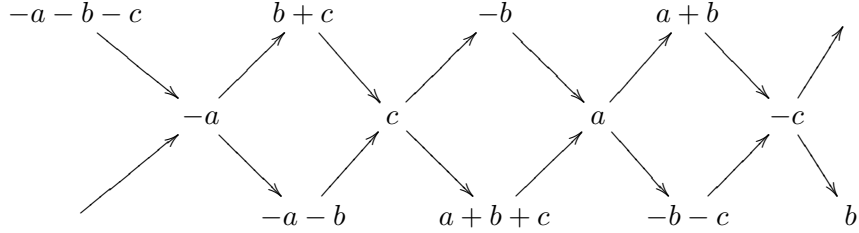
\xymatrix@C=0.8em{
 -a-b-c\ar@{->}[rd] && b+c \ar@{->}[rd] && -b \ar@{->}[rd] &&a+b \ar@{->}[rd] && \\
 &-a\ar@{->}[rd]\ar@{->}[ru] && c \ar@{->}[rd]\ar@{->}[ru] && a \ar@{->}[rd]\ar@{->}[ru] &&-c \ar@{->}[rd]\ar@{->}[ru]\\
\ar@{->}[ru]&&-a-b\ar@{->}[ru] && a+b+c\ar@{->}[ru] && -b-c\ar@{->}[ru] &&b  }$$
\end{Rem}

\subsection{}\label{section for mutation on Grothendieck group}
For each $\vx\in\bbL$, the equivalence
$${\bf{R}}_{\vx}: D^b(\co(\vx)\!^{\bot})\lra D^b({}^{\bot}\!\co(\vx))$$
induces an isomorphism
$$\Ref_{\vx}: K_0(\co(\vx)\!^{\bot})\lra K_0({}^{\bot}\!\co(\vx)).$$
Put $\wtilde{\Ref}_{\vx}:=\phi_{\vx}^l \Ref_{\vx} \psi_{\vx}^r$, i.e., $\wtilde{\Ref}_{\vx}$
is determined by the following commutative diagram
$$\xymatrix{
  K_0(\co(\vx)\!^{\bot})\ar[r]^{\Ref_{\vx}}\ar[d]_{\phi_{\vx}^r}& K_0({}^{\bot}\!\co(\vx))\ar[d]^{\phi_{\vx}^l}\\
              K_0(\coh\bbX)/\bbZ\delta \ar[r]^{\wtilde{\Ref}_{\vx}}&K_0(\coh\bbX)/\bbZ\delta}$$

\begin{Lem}\label{identity on root system}
For each $\vx\in\bbL$ and $X\in\co(\vx)\!^{\bot}$, $$\wtilde{\Ref}_{\vx}(\overline{[X]})=\overline{[X]}-(X,\co(\vx))\overline{[\co(\vx)]}.$$
\end{Lem}

\begin{pf} It is known that $\langle\co(\vx), X\rangle=0$ for $X\in\co(\vx)\!^{\bot}$.
Then by \eqref{approaximation-1},
$$\Ref_{\vx}([X]) = [X]-\langle X, \co(\vx)\rangle[\co(\vx)]
 = [X]-( X, \co(\vx))[\co(\vx)].$$
 This implies that
$$\wtilde{\Ref}_{\vx}(\overline{[X]})=\phi_{\vx}^l \Ref_{\vx} \psi_{\vx}^r(\overline{[X]})
=\phi_{\vx}^l \Ref_{\vx} ([X]) =\overline{[X]}-(X,\co(\vx))\overline{[\co(\vx)]}.$$
\end{pf}

For each vertex $\oz=\oz_{ij}$ in the star shaped quiver $Q=(I,Q_1)$, we have
the corresponding simple reflection $\refl_{ij}:=\refl_\oz$ in $W=W(Q)$.
The main result of this section is as follows.

\begin{Thm}
The following equalities hold in $\text{Aut}(\bbZ I)$:
$$\refl_\star=\phi \wtilde\Ref_0\phi^{-1} \; \text{and\;}\; \refl_{ij}=\phi\wtilde{\Ref}_{(j-1)\vx_i}
 \wtilde{\Ref}_{j\vx_i}\wtilde{\Ref}_{(j-1)\vx_i}\phi^{-1}  \; \text{for\;}\; (i,j)\in\whI,$$
 that is, the following diagram commutes:
 $$\xymatrix{
   K_0(\coh\bbX)/\bbZ\delta\ar[rrr]^{\wtilde{\Ref}_{(j-1)\vx_i}\wtilde{R}_{j\vx_i}\wtilde{\Ref}_{(j-1)\vx_i}}
   \ar[d]_{\phi}&&& K_0(\coh\bbX)/\bbZ\delta\ar[d]^{\phi}\\
              \bbZ I\ar[rrr]^{\refl_{ij}}&&&\bbZ I}$$
\end{Thm}

\begin{pf} We first consider the case $\oz=\oz_{i0}=\star$.
By Lemma \ref{identity on root system},
$$\wtilde{\Ref}_{0}(\overline{[X]})=\overline{[X]}-(X,\co)\overline{[\co]},\;
\forall\,\overline{[X]}\in  K_0(\coh\bbX)/\bbZ\delta.$$
 On the other hand,
$$\refl_\star(\az)=\az-(\az,\az_\star)\az_\star,\;\forall\,\az\in\bbZ I.$$
Then the equality $\refl_\star=\phi \wtilde\Ref_0\phi^{-1}$ follows from
the fact $\phi(\overline{[\co]})=\az_\star$.

For each $(i,j)\in\whI$, by $\bz_{ij}$ we denote the dimension
vector of the indecomposable $\bbF Q$-module $P_{\oz_{ij}}$. Then we have the associated
reflection $\refl_{\bz_{ij}}\in W(Q)$. By \eqref{value-F_c}, $\phi(\overline{[\co(j\vx_i)]})=\bz_{ij}$.
Thus, applying Lemma \ref{identity on root system} gives the commutative diagram
 $$\xymatrix{
    K_0(\coh\bbX)/\bbZ\delta\ar[r]^{\wtilde{\Ref}_{(j-1)\vx_i}}\ar[d]_{\phi}
   &K_0(\coh\bbX)/\bbZ\delta\ar[r]^{\wtilde{\Ref}_{j\vx_i}}\ar[d]^{\phi}
   &K_0(\coh\bbX)/\bbZ\delta\ar[r]^{\wtilde{\Ref}_{(j-1)\vx_i}}\ar[d]_{\phi}
   &K_0(\coh\bbX)/\bbZ\delta\ar[d]_{\phi}\\
    \bbZ I\ar[r]^{\refl_{\bz_{i,j-1}}}
   &\bbZ I\ar[r]^{\refl_{\bz_{ij}}}
   &\bbZ I\ar[r]^{\refl_{\bz_{i,j-1}}} & \bbZ I    }$$
Hence, it remains to show that $\refl_{ij}=\refl_{\bz_{i,j-1}} \refl_{\bz_{ij}}\refl_{\bz_{i,j-1}}$.

Now fix $1\leq i\leq t$. If $j=1$, then it can be directly checked that
$\refl_{i1}=\refl_\star \refl_{\bz_{i1}}\refl_\star$.
For $2\leq j\leq p_i-1$, we are reduced to consider the full subquiver of $Q$ consisting of vertices
$\oz_{ij}$ for $0\leq j\leq p_i-1$, which is a linear quiver $A_{p_i}$. Thus,
the equality $\refl_{ij}=\refl_{\bz_{i,j-1}} \refl_{\bz_{ij}}\refl_{\bz_{i,j-1}}$
 follows from Lemma \ref{A-type}(2).
\end{pf}

\bigskip

\section{Ringel--Hall Lie algebra of $\bbX$ and Tits' automorphisms}
\def\ind{{\rm ind}\,}
\def\ra{\rightarrow}
\def\fkQ{{\frak Q}}

In this section we follow \cite{PX00,LinP,CB10} to define the Ringel--Hall Lie algebra $\fL(\cR_{\bbX})$ of the
root category $\cR_{\bbX}$ of $\coh\bbX$ for a weighted projective line $\bbX$, as well as
the generic composition Lie algebra $\fL_c(\cR_{\bbX})$. We then show that the mutation functors
defined in the previous section induce isomorphisms between certain subalgebras
of $\fL_c(\cR_{\bbX})$. By identifying these Lie subalgebras with the Kac--Moody Lie algebra $\fg_Q$,
these isomorphisms give rise to Tits' automorphisms of $\fg_Q$, where $Q$ is the star shaped quiver
associated with $\bbX$.

\subsection{} We first associate a complex Lie algebra with the weighted
projective line $\bbX$. Let $\bbF$ be a finite field with $q$
elements and $\bbX=\bbX_\bbF$ denote the weighted projective line over $\bbF$ associated with
$\bf p$ and ${\boldsymbol\lambda}$. Consider the root category
$$\cR_{\bbX_\bbF}=D^b(\coh\bbX_\bbF)/[2]$$
 of $\coh\bbX_\bbF$. It is known that the Grothendieck group $K_0(\cR_{\bbX_\bbF})$
is isomorphic to $K_0(\coh\bbX_\bbF)$ which can be identified with $\what\Gamma=\bbZ I\oplus\bbZ \dz$.
Thus, in the following we write $K_0(\cR_{\bbX})$ instead of $K_0(\cR_{\bbX_\bbF})$.

Given three objects
$X,Y,Z\in \ind \cR_{\bbX_\bbF}$, define
$$F^Z_{X,Y}=|\{\text{triangles $Y\ra Z\ra X\ra X[1]$}\}/\Aut(X)\times\Aut(Y)|.$$
 Let $\llz_\bbF$ denote the quotient ring $\bbZ/(|\bbF|-1)\bbZ$. By \cite{PX00,LinP}, the free
$\llz_\bbF$-module
$$\fL(\cR_{\bbX_\bbF})=(\llz_\bbF\otimes_\bbZ K_0(\cR_{\bbX})\oplus \bigoplus_{X\in \ind \cR_{\bbX_\bbF}}\llz_\bbF u_X$$
 admits a Lie algebra structure whose Lie bracket is defined by
$$[u_X,u_Y]=\begin{cases} \displaystyle\sum_{Z\in \ind \cR_{\bbX_\bbF}}(F^Z_{X,Y}-F^Z_{Y,X})u_Z, \quad&\text{if $X\not\cong Y[1]$};\\
                  1\otimes \frac{[X]}{d(X)}, & \text{if $X\cong Y[1]$}\end{cases}$$
and
$$[h_\az,u_X]=-(\az,[X])u_X,\;[h_\az,h_\bz]=0,\;\forall\, \az,\bz\in K_0(\cR_{\bbX}),$$
 where $X,Y\in \ind \cR_{\bbX_\bbF}$, $h_\az:=1\otimes \az$ and $d(X)=\dim_\bbF(\End(X)/ {\rm{rad\,}}\End(X))$.

We call an object $X\in \cR_{\bbX_\bbF}$ {\it exceptional} if
$$\End(X)\cong \bbF\;\text{ and }\; \Hom(X,X[1])=0.$$
A pair $(X,Y)$ of exceptional objects in $\cR_{\bbX_\bbF}$ is called an {\it exceptional
pair} if
$$\Hom(Y,X)=0\;\text{ and }\;\Hom(Y,X[1])=0.$$
 In general, a sequence $(X_1,\ldots,X_r)$ of objects in $\cR_{\bbX_\bbF}$ is called
an {\it exceptional sequence} if each pair $(X_i,X_j)$ with $1\leq i< j\leq r$
is an exceptional pair. Such a sequence is called {\it complete} if the minimal full
triangulated subcategory containing the $X_i$ coincides with $\cR_{\bbX_\bbF}$.

Now fix an algebraic closure $\overline\bbF$ of $\bbF$ and set
$$\Omega=\{E\subset \overline\bbF \mid \text{$\bbF\subseteq E$ is a finite extension}\}.$$
 Then for each $E\in \Omega$, we obtain in a similar way the Lie algebra
$\fL(\cR_{\bbX_E})$ defined over $\llz_E=\bbZ/(|E|-1)\bbZ$. Note that
 if $X$ is an exceptional object in $\cR_{\bbX_\bbF}$, then so is $X^E$ in $\cR_{\bbX_E}$.
 Consider the infinite direct product
$$\prod_{E\in\Omega} \fL(\cR_{\bbX_E})$$
of Lie algebras and define $\fL_c(\cR_{\bbX})_\bbZ$ to be its Lie subalgebra
generated by the elements
$$h_\az:=(h_\az)_E,\; {\bf 1}_{X}:=(u_{X^E})_E,$$
where $\az\in K_0(\cR_{\bbX})$ and $X$ is an exceptional object in $\ind \cR_{\bbX_\bbF}$.
By the construction, $\fL_c(\cR_{\bbX})_\bbZ$
is a Lie algebra over $\bbZ$. Finally, set
$$\fL_c(\cR_{\bbX})=\fL_c(\cR_{\bbX})_\bbC=\bbC\otimes_\bbZ\fL_c(\cR_{\bbX})_\bbZ,$$
called the \emph{{\rm(}generic{\rm )} composition Lie algebra} of $\bbX$.
It is clear that $\fL_c(\cR_{\bbX})$ is $K_0(\cR_{\bbX})$-graded with $\deg \b1_X=[X]$ and
$\deg h_\az=0$. In particular, $\fL_c(\cR_{\bbX})_0=K_0(\cR_{\bbX})\otimes_\bbZ\bbC$, which is
spanned by $h_\az$ for $\az\in K_0(\cR_{\bbX})$. According to
\cite[Prop.~7.3]{RuanZh}, for each complete exceptional sequence
$(X_1,\ldots,X_n)$ in $\cR_{\bbX_\bbF}$, $\fL_c(\cR_{\bbX})$ is generated by
$\b1_{X_i}$ and $\b1_{X_i[1]}$, $1\leq i\leq n$.

For each $r\in\bbN$ and $E\in\Omega$, we denote by ${\cal X}_r(E)$
the set of $X\in\ind \coh \bbX_E$ of type $r\dz$ satisfying
$\Hom(X, S_{ij})=0$ for all $1\leq i\leq t$ and $1\leq j\leq p_i-1$.
By \cite[Lem.~8]{CB10}, $\sum_{X\in {\cal X}_r(E)}d(X)=2$. Finally, set
$$\mu_{r,E}=\sum_{X\in {\cal X}_r(E)} d(X)u_X,\;\mu_{-r,E}=-\sum_{X\in {\cal X}_r(E)} d(X)u_{TX}\in \fL(\cR_{\bbX_E})$$
and define
$$\b1_{r\dz}=(\mu_{r,E})_{E\in\Omega},\,\b1_{-r\dz}=(\mu_{-r,E})_{E\in\Omega}\in \prod_{E\in\Omega} \fL(\cR_{\bbX_E}).$$
 By \cite[Sect.~3]{CB10}, we have the following result.

\begin{Lem}\label{property-av-elt} {\rm(1)} For each $r\in\bbN$ and $X\in{\cal X}_r(E)$,
$$F^{\co(r\vc)}_{X,\co}=1\;\text{ and }\; F^X_{\co[1],\co(r\vc)}=d(X).$$

{\rm(2)} For $r,s\in\bbZ$ with $r\not=0$, $[\b1_{r\dz},\b1_{\co(s\vc)}]=2\b1_{\co((r+s)\vc)}$.

{\rm(3)} For $r,s\in\bbZ$ with $r\not=s$,
$[\b1_{\co(r\vc)},\b1_{\co(s\vc)[1]}]=-\b1_{(r-s)\dz}$.

{\rm(4)} For each $r\in\bbN$, $\b1_{\pm r\dz}\in \fL_c(\cR_{\bbX}).$

\end{Lem}

Moreover, by \cite[Thm.~2]{CB10} and
\cite[Prop.~8.4]{RuanZh}, there is a surjective Lie algebra homomorphism
\begin{equation} \label{epi-loop-comp}
\pi: \cL\fg_{Q}\lra\fL_c(\cR_{\bbX})
\end{equation}
taking
$$\aligned
{}& e_{\oz r}\lmto\left \{ \begin{array}{ll}
               \b1_{S_{ij}^{(rp_i+1)}}, & \text{if $\oz=\oz_{ij}\not=\star$};\\
                \b1_{\co(r\vc)}, & \text{if $\oz=\star$},
             \end{array} \right.
f_{\oz r}\lmto\left \{ \begin{array}{ll}
               \b1_{S_{i,j-1}^{(rp_i-1)}}, & \text{if $\oz=\oz_{ij}\not=\star$};\\
                -\b1_{\co(-r\vc)[1]}, & \text{if $\oz=\star$},
             \end{array} \right.  \\
& c\,\lmto -1\otimes\dz=-h_\dz,\;\; h_{\oz r}\lmto \left \{ \begin{array}{ll}
               -h_{\az_\oz}, & \text{if $r=0$};\\
               \b1_{S_{ij}^{(rp_i)}}-\b1_{S_{i,j-1}^{(rp_i)}}, & \text{if $r\not=0,\,\oz=\oz_{ij}\not=\star$};\\
                \b1_{r\dz}, & \text{if $r\not=0,\,\oz=\star$};\\

             \end{array} \right.
\endaligned$$
 where $\b1_{S_{ij}^{(m)}}=-\b1_{Y[1]}$ with $Y$ the indecomposable torsion sheaf of length $-m$ with
 socle $\tau^{-1}S_{ij}$ if $m<0$. By identifying $K_0(\cR_{\bbX})$ with
$\what\Gamma=\bbZ I\oplus\bbZ\dz$, $\pi$ is a $\what\Gamma$-graded homomorphism.

\begin{Rem}\label{Kac-Moody-Lie-alg}  (1) If we take the category $\cR_{\bbX_\bbF}$
to be the root category $\cR(\co(\vx)\!^{\bot})$ of the subcategory $\co(\vx)\!^{\bot}$
for $\vx\in\bbL$, we obtain the associated (generic) composition
Lie algebra $\fL_c(\cR(\co(\vx)\!^{\bot}))$. For simplicity, write
$\fL_{\vx}^r=\fL_c(\cR(\co(\vx)\!^{\bot}))$.

(2) If we take the category $\cR_{\bbX_\bbF}$ to be the
root category $\cR(\bbF \Gamma)$ of the module category $\mod\,\bbF \Gamma$ over the path
algebra $\bbF \Gamma$ of an acyclic quiver $\Gamma=(\Gamma_0,\Gamma_1)$, then by \cite{PX00},
the associated (generic) composition Lie algebra $\fL_c(\cR(\bbF \Gamma))$
is generated by $\b1_{S_i}$ ($i\in \Gamma_0$), where the $S_i$ are simple $\bbF \Gamma$-modules.
Moreover, it is isomorphic to the Kac--Moody algebra $\fg_{\Gamma}$ of $\Gamma$.
\end{Rem}

\subsection{}\label{section of sub algebra}
In the following we will see that there are several different ways to embed
the Kac-Moody algebra $\fg_{Q}$ into $\fL_c(\cR_{\bbX})$. For each element $\vx\in\bbL$, the embedding
$D^b(\co(\vx)\!^{\bot})\hookrightarrow D^b(\coh\bbX_\bbF)$
induces an embedding of their root categories
$\cR(\co(\vx)\!^{\bot})\hookrightarrow \cR(\coh\bbX_\bbF)$.
By the construction in 4.1, this gives rise to an embedding of Lie algebras
$\fL_{\vx}^r \hookrightarrow\fL_c(\cR_{\bbX})$.

\begin{Prop} \label{grading non-preseved iso}
 For any $\vx\in\bbL$, there is a Lie algebra isomorphism
 $\Omega_{\vx}:\fg_{Q}\cong\fL_{\vx}^{r}$
taking
$$\aligned
{}& e_{\oz}\lmto\left \{ \begin{array}{ll}
               {\bf 1}_{\co(\vx-\vc)}, & \text{if $\oz=\star$};\\
                {\bf 1}_{S_{ij}(\vx)}, & \text{if $\oz=\oz_{ij}$},
             \end{array} \right.
f_{\oz}\lmto\left \{ \begin{array}{ll}
               -{\bf 1}_{\co(\vx-\vc)[1]}, & \text{if $\oz=\star$};\\
                -{\bf 1}_{S_{ij}(\vx)[1]}, & \text{if $\oz=\oz_{ij}$},
             \end{array} \right.\\
& h_{\oz}\lmto\left \{ \begin{array}{ll}
               -h_{[\co(\vx-\vc)]}, & \text{if $\oz=\star$};\\
                -h_{[S_{ij}(\vx)]}, & \text{if $\oz=\oz_{ij}$},
             \end{array} \right.
\endaligned$$
 where $(i,j)\in \whI$.
\end{Prop}

\begin{pf}  By Lemma \ref{derived equivalent to star shaped quiver}, the tilting object
$T_\vx^r=\bigoplus_{0<\vy\leq \vc} \co(\vx-\vy)$ in $\coh\bbX$
induces a triangle equivalence $D^b(\co(\vx)\!^{\bot})\cong D^b({\rm mod}\,\bbF Q)$
which takes the simple torsion sheaf $S_{i,l_i+j}$ with $(i,j)\in\whI$
to the simple $\bbF Q$-module $S_{\oz_{ij}}$  corresponding to the vertex $\oz_{ij}\in I\backslash\{\star\}$, and takes
$\co(\vx-\vc)$ to the simple projective $\bbF Q$-module $S_\star$ corresponding to the vertex $\star$.
The assertion follows from \cite[Thm.~4.7]{PX00}.
\end{pf}

\begin{Cor}\label{basis of right perp sub algebra} For each given $\vx$ of the normal form
$\vx=\sum\limits_{i=1}^tl_i\vx_i+l\vc$, the elements
$${\bf 1}_{X}, {\bf 1}_{X[1]}, h_{[X]}\;\;\text{for}\;
X\in\{\co(\vx-\vc), S_{ij}\mid 1\leq i\leq t,0\leq j\leq p_i-1, j\neq l_i\}$$
form a complete set of Chevalley generators of $\fL_{\vx}^{r}$.
\end{Cor}

\begin{pf} This is a direct consequence of the above proposition.
\end{pf}

Dually, by dealing with the subcategory ${}^{\bot}\!\co(\vx)$ for each $\vx\in\bbL$,
we can define the Lie algebra $\fL_{\vx}^{l}$ and obtain analogously
an embedding of Lie algebras
$\fL_{\vx}^l\hookrightarrow\fL_c(\cR_{\bbX})$, as well as an isomorphism $\fL_{\vx}^{l}\cong \fg_{Q}$.
Moreover, if $\vx$ has the normal form $\vx=\sum\limits_{i=1}^tl_i\vx_i+l\vc$,
then $\fL_{\vx}^{l}$ admits a set of generators
$$\b1_{X}, \b1_{X[1]}\;\text{ with }\;
X\in\{\co(\vx+\vc), S_{ij}\mid 1\leq i\leq t,1\leq j\leq p_i,j\neq l_i+1\}$$
 since the $\co(\vx+\vc)$ and $S_{ij}$ can be ordered in a way that they become a complete
exceptional sequence in ${}^{\bot}\!\co(\vx)$.

\subsection{}
In the following we show that the Kac--Moody Lie algebra $\fg_Q$ can be viewed
as a quotient algebra of $\fL_c(\cR_{\bbX})$. For this purpose, we define a Lie ideal
$\frak I$ of $\fL_c(\cR_{\bbX})$ by setting
$$\fI=\langle \b1_{\co(\vc)}-\b1_{\co}, \b1_{\co(\vc)[1]}-\b1_{\co[1]}\rangle.$$

\begin{Prop}\label{special elts in cI} The following elements lie in $\fI$:
  \begin{itemize}
    \item [(1)] $\b1_{\co(\vx+\vc)}-\b1_{\co(\vx)}$ and $\b1_{\co(\vx+\vc)[1]}-\b1_{\co(\vx)[1]}$,
    where $\vx\in\bbL$;
    \item [(2)] $\b1_{\pm\dz}+h_{[\co]}$;
    \item [(3)] $\b1_{r\dz}-\b1_{\dz}$ for all $r\in\bbZ\backslash\{0\}$;
    \item [(4)] $h_{\delta}$;% for any $r\in\bbN$;
    \item [(5)] $\b1_{S^{(k)}_{ij}}+\b1_{S_{i,j-k}^{(p_i-k)}[1]}$ and
    $\b1_{S^{(k)}_{ij}[1]}+\b1_{S_{i,j-k}^{(p_i-k)}}$, where $1\leq i\leq t, 0\leq j< p_i, 1\leq k\leq p_i-1$;
    \item [(6)] $\b1_{S_{ij}^{(p_i)}}-\b1_{S_{i,j-1}^{(p_i)}}+h_{[S_{ij}]}$ and $\b1_{S_{ij}^{(p_i)}[1]}-\b1_{S_{i,j-1}^{(p_i)}[1]}-h_{[S_{ij}]}$, where $1\leq i\leq t, 0\leq j< p_i$.
  \end{itemize}
\end{Prop}

\begin{pf} For the cases (1), (5) and (6), we only show that one of the given two elements belongs
to $\fI$. The proof for the other one is similar.

(1) Assume $\vx$ has the normal form $\vx=\sum\limits_{i=1}^t l_i\vx_i+l\vc$. We proceed induction
on $\ell(\vx)=|\{i\mid l_i\neq 0\}|$. We first consider the case $\ell(\vx)=0$. If $l=0$, this is trivial.
For $l\not=0$, we have by Lemma \ref{property-av-elt} that
$$2(\b1_{\co(l\vc+\vc)}-\b1_{\co(l\vc)})=[\b1_{l\dz},\b1_{\co(\vc)}-\b1_{\co}]\in \fI.$$
Now suppose the statement holds for all $\vx$ with $\ell(\vx)=k$ and take $\vy\in\bbL$
with $\ell(\vy)=k+1$. Then there exists $1\leq i\leq t$ such that $\vy=\vx+\vx_i$ with $\ell(\vx)=k$.
By the induction hypothesis,
$$\b1_{\co(\vy+\vc)}-\b1_{\co(\vy)}
        =[\b1_{S_{i,l_i+1}}, \b1_{\co(\vx+\vc)}-\b1_{\co(\vx)}]\in \fI.$$

(2) By Lemma \ref{property-av-elt}(3),
$$\aligned
{}& \b1_{\dz}+h_{[\co]}=[\b1_{\co[1]}, \b1_{\co(\vc)}-\b1_{\co}]\in \fI\;\text{ and }\\
&\b1_{-\dz}+h_{[\co]}=[\b1_{\co(\vc)[1]}-\b1_{\co[1]}, \b1_{\co}]\in \fI.\endaligned$$

(3) By (2), $\b1_{\dz}-\b1_{-\dz}\in \fI$. Moreover, for each $r\not=0, 1$, by Lemma \ref{property-av-elt}(3) again, $$\b1_{r\dz}-\b1_{(r-1)\dz}
    =-[\b1_{\co(r\vc)}, \b1_{\co[1]}-\b1_{\co(\vc)[1]}]\in\fI.$$

(4) The equality
$$[\b1_{\co(\vc)}-\b1_{\co}, \b1_{\co(\vc)[1]}+\b1_{\co[1]}]
        =h_{[\co(\vc)]}-h_{[\co]}-(\b1_{\dz}-\b1_{-\dz})\in\fI$$ implies
        that $h_\delta=h_{[\co(\vc)]}-h_{[\co]}\in\fI$.

(5) By the definition, $\b1_{S^{(k)}_{ij}}+\b1_{S_{i,j-k}^{(p_i-k)}[1]}
        =[\b1_{\co(\vc+(j-k)\vx_i)[1]},\b1_{\co(j\vx_i+\vc)}-\b1_{\co(j\vx_i)}]\in \fI$.

(6) We have
$$\begin{array}{rl}
[\b1_{S_{ij}}+\b1_{S_{i,j-1}^{(p_i-1)}[1]}, \b1_{S_{i,j-1}^{(p_i-1)}}]
&=\b1_{S_{ij}^{(p_i)}}-\b1_{S_{i,j-1}^{(p_i)}}-h_{[S_{i,j-1}^{(p_i-1)}]}\\
&=\b1_{S_{ij}^{(p_i)}}-\b1_{S_{i,j-1}^{(p_i)}}+h_{[S_{ij}]}-h_{\delta}\in \fI.
\end{array}$$
This implies that $\b1_{S_{ij}^{(p_i)}}-\b1_{S_{i,j-1}^{(p_i)}}+h_{[S_{ij}]}\in\fI$.

\end{pf}

Given an element $\b1_{X}\in\fL_c(\cR_{\bbX})$, we denote by $\bar{\b1}_{X}$
its image in $\fL_c(\cR_{\bbX})/\fI$. We have the following results.

\begin{Lem} \label{Lie bracket} Assume $\vx$ has the normal form $\vx=\sum\limits_{i=1}^tl_i\vx_i+l\vc$. Then for each
  $(i,j)\in \whI, 1\leq k\leq p_i-1$, the following equalities hold in $\fL_c(\cR_{\bbX})/\fI$:
  \begin{itemize}
    \item[(1)] $[\bar{\b1}_{\co(\vx)}, \bar{\b1}_{S_{ij}^{(k)}}]=-\delta_{j, l_i+k}\bar{\b1}_{\co(\vx+k\vx_i)}$;
    \item[(2)] $[\bar{\b1}_{\co(\vx)}, \bar{\b1}_{S_{ij}^{(k)}[1]}]=\delta_{j, l_i}\bar{\b1}_{\co(\vx-k\vx_i)}$;
    \item[(3)] $[\bar{\b1}_{\co(\vx)}, \bar{\b1}_{\co(\vx-k\vx_i)[1]}]
        =-\bar{\b1}_{S_{i,l_i}^{(k)}}$;
    \item[(4)] $[\bar{\b1}_{\co(\vx)}, \bar{\b1}_{\co(\vx+k\vx_i)[1]}]
        =\bar{\b1}_{S_{i,l_i+k}^{(k)}[1]}$.
  \end{itemize}
\end{Lem}

\begin{pf} Statements (1) and (4) can be deduced from the following triangle in the root category $\cR_{\bbX}$:
  $$\xymatrix{
  S_{i,l_i+k}^{(k)}[1] \ar[r] & \co(\vx)\ar[r]^{x_i^k\quad} & \co(\vx+k\vx_i)\ar[r] &S_{i,l_i+k}^{(k)}
  };$$
  while statements (2) and (3) can be deduced from the following triangle and its rotation:
  $$\xymatrix{
   \co(\vx-k\vx_i)\ar[r]^{\quad\; x_i^k} & \co(\vx)\ar[r] &S_{i,l_i}^{(k)} \ar[r]&\co(\vx-k\vx_i)[1]
  }.$$
\end{pf}

\begin{Prop}\label{property of quotient algebra} Let $\vx\in\bbL$ have the normal form
$\vx=\sum\limits_{i=1}^tl_i\vx_i+l\vc$. Then the Lie algebra $\fL_c(\cR_{\bbX})/\fI$ is generated
by
$$\{\bar{\b1}_{\co(\vx)}, \bar{\b1}_{\co(\vx)[1]},\bar{\b1}_{S_{ij}},\bar{\b1}_{S_{ij}[1]}\mid
1\leq i\leq t, 0\leq j\leq p_i-1, j\neq l_i\}.$$
\end{Prop}

\begin{pf} First of all, the coherent sheaves
$$\co(\vx), \co(\vx+\vc), \;S_{ij}\; (1\leq i\leq t, 0\leq j\leq p_i-1, j\neq l_i)$$
can be ordered so that they form a complete exceptional sequence in $\coh\bbX$,
as well as in $\cR_\bbX$. Then the assertion follows from
\cite[Prop.~8.4]{RuanZh} and Proposition \ref{special elts in cI}(1).
  \end{pf}

Recall from \S2.3 and \S4.1 that both $\cL\fg_Q$ and $\fL_c(\cR_{\bbX})$ are $\what\Gamma$-graded with
$$(\cL\fg_Q)_0=(\fg_Q)_0\oplus\bbC c\;\text{ and }\;\fL_c(\cR_{\bbX})_0 \cong \bbC\otimes_\bbZ\what\Gamma, $$
 where $\what\Gamma=\bbZ I\oplus\bbZ\dz$. By the definition, $\pi: \cL\fg_{Q}\ra\fL_c(\cR_{\bbX})$
restricts to an isomorphism
$$\pi_0:(\cL\fg_Q)_0\lra \fL_c(\cR_{\bbX})_0;\; h_\oz\lmto h_\oz\;(\oz\in I),\, c\,\lmto -h_\dz.$$
  Further, consider the quotient group $\what\Gamma/\bbZ\delta$ which is
 naturally isomorphic to $\bbZ I$. Then $\cL\fg_Q$ and $\fL_c(\cR_{\bbX})$
become $\what\Gamma/\bbZ\delta$-graded with homogeneous spaces defined by $$(\cL\fg_Q)_{\bar{\alpha}}:=\bigoplus_{n\in\bbZ}(\cL\fg_Q)_{\alpha+n\delta}\;\text{ and }\;
\fL_c(\cR_{\bbX})_{\bar{\alpha}}:=\bigoplus_{n\in\bbZ}\fL_c(\cR_{\bbX})_{\alpha+n\delta},$$
  where $\bar{\alpha}=\az+\bbZ\dz\in \what\Gamma/\bbZ\delta$, and $\pi$ becomes a
$\what\Gamma/\bbZ\dz$-graded homomorphism. Moreover, $\fI$ becomes a homogenous ideal of
$\fL_c(\cR_{\bbX})$ with respect to the $\what\Gamma/\bbZ\dz$-grading.

 Let $\iota:\fg_Q \hookrightarrow \cL\fg_Q$ be the natural embedding which takes
$$e_{\oz}\lmto e_{\oz 0},\;f_{\oz}\lmto f_{\oz 0},\;h_{\oz}\lmto h_{\oz 0},\;\forall\,\oz\in I.$$
 By $\Phi$ we denote the composition of $\pi: \cL\fg_Q\to \fL_c(\cR_{\bbX})$ with
the natural embedding and projection
$$\Phi: \fg_Q \stackrel{\iota}{\hookrightarrow} \cL\fg_Q \stackrel{\pi}{\lra} \fL_c(\cR_{\bbX})
\stackrel{\eta}{\twoheadrightarrow} \fL_c(\cR_{\bbX})/\fI.$$
Then by \eqref{epi-loop-comp}, $\Phi$ takes
$$\aligned
{}& e_{\star}\longmapsto \bar{\b1}_{\co},\; f_{\star}\longmapsto -\bar{\b1}_{\co[1]},\;
h_{\star}\longmapsto -h_{\overline{[\co]}},\\
& e_{\oz}\longmapsto \bar{\b1}_{S_{ij}},\; f_{\oz}\longmapsto -\bar{\b1}_{S_{ij}[1]},\;
h_{\oz}\longmapsto -h_{\overline{[S_{ij}]}},\; \forall\,\oz=\oz_{ij}\in I\setminus\{\star\}.
\endaligned$$
 In the following we are going to prove that $\Phi$ is an isomorphism. We need the following
lemma which is well known.

\begin{Lem} \label{ideal-in-KM} Let $I$ be a homogeneous ideal of $\fg_Q$. If
$I\cap(\fg_Q)_0=0$, then $I=0$.
\end{Lem}

\begin{pf} Let $\fg^{\text{\rm KM}}$ denote the Kac--Moody Lie algebra attached to a minimal
realization of $C_Q$ given in \cite[Ch.~1]{Kac90}. Then $\fg_Q$ is the derived
algebra of $\fg^{\text{\rm KM}}$. Since $I$ is a homogeneous
ideal of $\fg_Q$, $I$ is also an ideal of $\fg^{\text{\rm KM}}$. By \cite[Prop.~1.7]{Kac90},
we have
$$\fg_Q=[\fg^{\text{\rm KM}},\fg^{\text{\rm KM}}]\subseteq I\;\text{ or }\; I\subseteq Z(\fg^{\text{\rm KM}}),$$
 where $Z(\fg^{\text{\rm KM}})$ is the center of $\fg^{\text{\rm KM}}$. This implies that
 $I=0$ since $Z(\fg^{\text{\rm KM}})\subseteq (\fg_Q)_0$.
\end{pf}

Following \cite[\S7.2]{Kac90}, let $\what L(\fg_Q)$ be the one-dimensional
central extension of the loop algebra of $\fg_Q$, that is,
$${\what L}(\fg_Q)=\bbC[t,t^{-1}]\otimes\fg_Q\oplus \bbC c$$
 with the Lie bracket
$$[t^l\otimes x, t^m\otimes y]=t^{l+m}\otimes [x,y]+l\dz_{l+m,0}\langle x,y\rangle c,$$
 where $l,m\in\bbZ$, $x,y\in\fg_Q$, and $\langle-,-\rangle$ is a standard invariant bilinear form on $\fg_Q$.
 Then $\what L(\fg_Q)$ is also $\what\Gamma$-graded with
$$\deg (t^k\otimes e_\oz)=\az_\oz+k\dz,\;\deg (t^k\otimes f_\oz)=-\az_\oz+k\dz,\;
\deg (t^k\otimes h_\oz)=k\dz,\;\deg c=0.$$
Furthermore, there is a surjective Lie algebra homomorphism $\Psi:\cL\fg_{Q}\ra \what L(\fg_Q)$
given by
$$e_{\oz,k}\lmto t^k\otimes e_\oz,\;f_{\oz,k}\lmto t^k\otimes f_\oz,\;
h_{\oz,k}\lmto t^k\otimes h_\oz,\;c\,\lmto c,$$
 where $\oz\in I$ and $k\in\bbZ$. Clearly, $\Psi$ is $\what\Gamma$-graded and induces an isomorphism
$$\Psi_0: (\cL\fg_{Q})_0=(\fg_{Q})_0\oplus\bbC c\lra 1\otimes (\fg_{Q})_0\oplus\bbC c=(\what L(\fg_Q))_0,$$
and it is proved in \cite{MRY} that
\begin{equation} \label{kernel-Psi}
\Ker(\Psi)\subset \bigoplus_{\az\in\Delta^{\rm im}(Q)\atop k\in\bbZ}(\cL\fg_{Q})_{\az+k\dz}.
\end{equation}
 Thus, it follows from \eqref{kernel-Psi} that for each $k\in\bbZ\backslash\{0\}$,
$$\dim (\cL\fg_{Q})_{k\dz}=\dim \what L(\fg_Q)_{k\dz}=\dim (\fg_{Q})_0=|I|.$$
 Hence, by considering the $\what\Gamma/\bbZ\delta$-grading of $\what L(\fg_Q)$, $\Psi$
induces an isomorphism
$$\Psi_{\bar0}:(\cL\fg_{Q})_{\bar0}\lra\what L(\fg_Q)_{\bar0}.$$

 On the other hand, by an argument in the proof of \cite[Prop.~8.8]{LinP},
$$\dim \fL_c(\cR_{\bbX})_{k\dz}\geq 1+\sum_{i=1}^t (p_i-1)=|I|.$$
 Since $\pi: \cL\fg_{Q}\ra\fL_c(\cR_{\bbX})$ is $\what\Gamma$-graded and surjective,
we conclude that
$$\dim (\cL\fg_{Q})_{k\dz}=\dim \fL_c(\cR_{\bbX})_{k\dz},\;\forall\,k\in\bbZ\backslash\{0\}.$$
 Therefore, $\pi$ induces an isomorphism
$$\pi_{\bar 0}: (\cL\fg_{Q})_{\bar0}\lra\fL_c(\cR_{\bbX})_{\bar0}.$$

\begin{Prop} \label{iso-g-L(X)/I}
  The homomorphism $\Phi: \fg_Q\to \fL_c(\cR_{\bbX})/\fI$ is an isomorphism.
\end{Prop}

\begin{pf} We simply identify $\what\Gamma/\bbZ\delta$ with $\bbZ I$.
Then $\Phi$ becomes a $\bbZ I$-graded Lie algebra homomorphism.
By the definition of $\iota$ and Proposition \ref{property of quotient algebra},
$\Phi$ is surjective. It remains to show the injectivity of $\Phi$.

Let $\fI'$ be the ideal of $\cL\fg_Q$ generated by $e_{\star,1}-e_{\star,0}$ and $f_{\star,1}-f_{\star,0}$
and $\fI''$ be the ideal of $\what L(\fg_Q)$
generated by $(t-1)\otimes e_{\star}$ and $(t-1)\otimes f_{\star}$.
Then $\Psi(\fI')=\fI''$, and $\fI'$ (resp. $\fI''$) becomes a homogenous ideal of
$\cL\fg_Q$ (resp., $\what L(\fg_Q)$) with respect to the $\what\Gamma/\bbZ\dz$-grading.
By the definition,
$$(t-1)\otimes h_\star=[(t-1)\otimes e_{\star},1\otimes f_\star]\in\fI''\;\text{ and }\;
c=\frac{1}{2}[(t-1)\otimes h_\star,t^{-1}\otimes h_\star]\in\fI''.$$
 Since the ideal of $\fg_Q$ generated by $e_\star$ and $f_\star$
is $\fg_Q$ itself, it follows that
$${\frak K}:=\bigoplus_{m\in\bbZ}(t^m-t^{m-1})\otimes\fg_Q\oplus \bbC c\subseteq \fI''.$$
On the other hand, $\frak K$ is an ideal of $\what L(\fg_Q)$ consisting of
$(t-1)\otimes e_{\star}$ and $(t-1)\otimes f_{\star}$. Therefore,
$$\fI''={\frak K}=\bigoplus_{m\in\bbZ}(t^m-t^{m-1})\otimes\fg_Q\oplus \bbC c.$$
 This implies that
$$\what L(\fg_Q)_{\bar0}/\fI''_{\bar0}\cong (\fg_Q)_0.$$

Clearly, $\Psi_{\bar0}$ and $\pi_{\bar0}$ induce isomorphisms
$$\what L(\fg_Q)_{\bar0}/\fI''_{\bar0}\cong (\cL\fg_{Q})_{\bar0}/\fI'_{\bar0}
\cong \fL_c(\cR_{\bbX})_{\bar{0}}/\fI_{\bar{0}}.$$
 Thus, $\iota$ induces an isomorphism
$$\iota_0:(\fg_Q)_0\lra (\cL\fg_{Q})_{\bar0}/\fI'_{\bar0}.$$
 Consequently, the restriction of $\Phi$ to $(\fg_Q)_0$ gives an isomorphism
$$\Phi_0:(\fg_Q)_0\lra (\fL_c(\cR_{\bbX})/\fI)_{\bar0}=\fL_c(\cR_{\bbX})_{\bar{0}}/\fI_{\bar{0}}.$$

Because $\Phi$ is $\bbZ I$-graded, $\Ker(\Phi)$ is a homogeneous
ideal of $\fg_Q$. The bijectivity of $\Phi_0$ implies that $\Ker(\Phi)\cap (\fg_Q)_0=0$.
 Applying Lemma \ref{ideal-in-KM} shows that $\Ker(\Phi)=0$. Therefore,
$\Phi$ is an isomorphism, as desired.
\end{pf}

\subsection{}
Fix an element $\vx\in\bbL$ and denote by $\theta_{\vx}^r$ the composition of the natural
embedding and projection:
$$\theta_{\vx}^r: \fL_{\vx}^r\hookrightarrow \fL_c(\cR_{\bbX}) \twoheadrightarrow \fL_c(\cR_{\bbX})/\fI.$$

\begin{Prop}\label{iso of sub-quotient}
  The Lie homomorphism $\theta_{\vx}^r: \fL_{\vx}^r\to \fL_c(\cR_{\bbX})/\fI$ is an isomorphism.
\end{Prop}

\begin{pf} Suppose $\vx$ has the normal form $\vx=\sum\limits_{i=1}^tl_i\vx_i+l\vc$.
By Corollary \ref{basis of right perp sub algebra} and Proposition \ref{property of quotient algebra},
$\fL_{\vx}^r$ (resp., $\fL_c(\cR_{\bbX})/\fI$) is generated by the elements
$\b1_{X}$ and $\b1_{X[1]}$ (resp., $\bar{\b1}_{X}$ and $\bar{\b1}_{X[1]}$),
where
$$X\in\{\co(\vx-\vc), S_{ij}\mid 1\leq i\leq t, 0\leq j\leq p_i-1,j\neq l_i\}.$$
Moreover, $\theta_{\vx}^r$ takes $\b1_{X}\mapsto \bar{\b1}_{X}$ and $\b1_{X[1]}\mapsto \bar{\b1}_{X[1]}$.
Consequently, $\theta_{\vx}^r$ is surjective.

On the other hand, the restriction
of $\theta_{\vx}^r$ to the Cartan subalgebra gives an isomorphism
$$(\fL_{\vx}^r)_0=\bbC\otimes_{\bbZ} K_0(\co(\vx)^{\perp}) \cong\bbC\otimes_{\bbZ} (K_0(\coh\bbX)/\bbZ\delta).$$
Thus, $\Ker(\theta_{\vx}^r)\cap (\fL_{\vx}^r)_0=0$. By Lemma \ref{ideal-in-KM}, $\Ker(\theta_{\vx}^r)=0$,
i.e., $\theta_{\vx}^r$ is injective. Hence, $\theta_{\vx}^r$ is an isomorphism.
\end{pf}

\begin{Cor}\label{basis of quotient algebra} For each given $\vx$ of the normal form
$\vx=\sum\limits_{i=1}^tl_i\vx_i+l\vc$, the elements
$$\bar{\b1}_{X}, \bar{\b1}_{X[1]}, h_{\overline{[X]}}\;\;\text{for}\;
X\in\{\co(\vx), S_{ij}\mid 1\leq i\leq t,0\leq j\leq p_i-1, j\neq l_i\}$$
form a complete set of Chevalley generators of $\fL_c(\cR_{\bbX})/\fI$.
\end{Cor}

\begin{pf} Since $\bar{\b1}_{\co(\vx)[1]}=\bar{\b1}_{\co(\vx-\vc)[1]}$, the assertion
follows from Corollary \ref{basis of right perp sub algebra}.
\end{pf}

Dually, by $\theta_{\vx}^l$ we denote the composition of the
natural embedding and projection
$$\theta_{\vx}^l: \fL_{\vx}^l\hookrightarrow \fL_c(\cR(\bbX)) \twoheadrightarrow \fL_c(\cR_{\bbX})/\fI.$$
Then $\theta_{\vx}^l$ is an isomorphism of Lie algebras, too.
As a consequence, $\fL_c(\cR_{\bbX})/\fI$ admits a set of generators
$$\bar{\b1}_{X}, \bar{\b1}_{X[1]}\;\text{ with }\;
X\in\{\co(\vx), S_{ij}\mid 1\leq i\leq t,1\leq j\leq p_i,j\neq l_i+1\}.$$

\begin{Rems} (1) Under the isomorphism $\theta_{\vx}^r$ (resp., $\theta_{\vx}^l$),
we can view the subalgebra $\fL_{\vx}^r$ (resp. $\fL_{\vx}^l$) as the quotient
algebra $\fL_c(\cR_{\bbX})/\fI$. We emphasize that both $\fL_{\vx}^r$ and $\fL_{\vx}^l$
depend on the choice of $\vx$, but $\fL_c(\cR_{\bbX})/\fI$ does not.

(2) Combining with Propositions \ref{iso-g-L(X)/I} and \ref{iso of sub-quotient} gives an
isomorphism between Lie algebras $\fg_Q$ and $\fL_{\vx}^r$. This isomorphism preserves the
$\bbZ I$-grading, which is not the case for the isomorphism given in
Proposition \ref{grading non-preseved iso}.
\end{Rems}

\begin{Prop}\label{iso of Kac-Moody}
For each $(i,j)\in\widehat{I}$, there is a Lie algebra automorphism
$\widetilde{\Omega}_{ij}: \fg_{Q}\ra \fg_{Q}$ taking
$$\aligned
{}& e_{\oz}\lmto\left \{ \begin{array}{ll}
               [e_{\oz_{i,p_i-j}},\cdots,e_{\oz_{i,1}},e_{\star}], & \text{if $\oz=\star$};\\
               (-1)^{p_i}[f_{\oz_{i,p_i-1}},\cdots,f_{\oz_{i,1}}], & \text{if $\oz=\oz_{ij}$};\\
               e_{\oz_{k,l-j\delta_{ik}}}, & \text{if $\oz=\oz_{kl}\neq \oz_{ij}$},
             \end{array} \right. \\
{}& f_{\oz}\lmto\left \{ \begin{array}{ll}
              (-1)^{p_i-j} [f_{\oz_{i,p_i-j}},\cdots,f_{\oz_{i,1}},f_{\star}], & \text{if $\oz=\star$};\\
              {}[e_{\oz_{i,p_i-1}},\cdots,e_{\oz_{i,1}}], & \text{if $\oz=\oz_{ij}$};\\
               f_{\oz_{k,l-j\delta_{ik}}}, & \text{if $\oz=\oz_{kl}\neq \oz_{ij}$},
             \end{array} \right.\\
& h_{\oz}\lmto\left \{ \begin{array}{ll}
              h_{\oz_{i,p_i-j}}+\cdots+h_{\oz_{i,1}}+h_{\star}, & \text{if $\oz=\star$};\\
              {}-(h_{\oz_{i,p_i-1}}+\cdots+h_{\oz_{i,1}}), & \text{if $\oz=\oz_{ij}$};\\
               h_{\oz_{k,l-j\delta_{ik}}}, & \text{if $\oz=\oz_{kl}\neq \oz_{ij}$},
             \end{array} \right.
\endaligned$$
where $(k,l)\in\whI$. Moreover, $\widetilde{\Omega}_{ij}((\fg_{Q})_{\alpha})=(\fg_{Q})_{\varpi(\alpha)},
\;\forall\,\alpha\in\Delta(Q)$, where
$$\varpi=(\refl_{i,p_i-j}\cdots \refl_{i,p_i-1})(\refl_{i,p_i-j-1}\cdots \refl_{i,p_i-2})
\cdots(\refl_{i,1}\cdots \refl_{i,j})\in W(Q).$$
\end{Prop}

\begin{pf} By combining Propositions \ref{grading non-preseved iso}, \ref{iso-g-L(X)/I} and \ref{iso of sub-quotient}
together, we obtain an automorphism
$$\widetilde{\Omega}_{ij}=\Omega_{j\vx_i}^{-1}\circ (\theta_{j\vx_i}^r)^{-1}\circ\Phi: \fg_Q\lra \fL_c(\cR_{\bbX})/\fI\lra\fL_{j\vx_i}^r\lra \fg_Q$$
which takes
$$\aligned
\widetilde{\Omega}_{ij}(e_\star)&=\Omega_{j\vx_i}^{-1}(\b1_{\co})
=\Omega_{j\vx_i}^{-1}([\b1_{S_{i,p_i}},\cdots,\b1_{S_{i,j+1}},\b1_{\co(j\vx_i-\vc)}])
=[e_{\oz_{i,p_i-j}},\cdots,e_{\oz_{i,1}},e_{\star}],\\
\widetilde{\Omega}_{ij}(e_{\oz_{ij}})&=\Omega_{j\vx_i}^{-1}(-\b1_{S_{i,j-1}^{(p_i-1)}[1]})
=-\Omega_{j\vx_i}^{-1}([\b1_{S_{i,j-1}[1]},\b1_{S_{i,j-2}[1]},\cdots,\b1_{S_{i,j+1}[1]}])\\
&=(-1)^{p_i}[f_{\oz_{i,p_i-1}},\cdots,f_{\oz_{i,1}}],\\
\widetilde{\Omega}_{ij}(e_{\oz_{kl}})&=\Omega_{j\vx_i}^{-1}(\b1_{S_{k,l}})=e_{\oz_{k,l-j\delta_{ik}}}.
\endaligned$$
 Dually, we have
$$\aligned
\widetilde{\Omega}_{ij}(f_\star)&=(-1)^{p_i-j} [f_{\oz_{i,p_i-j}},\cdots,f_{\oz_{i,1}},f_{\star}],\\
\widetilde{\Omega}_{ij}(f_{\oz_{ij}})&=[e_{\oz_{i,p_i-1}},\cdots,e_{\oz_{i,1}}],\\
\widetilde{\Omega}_{ij}(f_{\oz_{kl}})&=f_{\oz_{k,l-j\delta_{ik}}}.
\endaligned$$
 Moreover, it is easy to verify that
$$\aligned
{}&\varpi(\alpha_{\star})=\alpha_{\star}+\alpha_{i1}+\cdots+\alpha_{i,p_i-j},
\;\varpi(\alpha_{ij})=-(\alpha_{i1}+\alpha_{i2}+\cdots+\alpha_{i,p_i-1}),\\
& \varpi(\alpha_{kl})=\alpha_{k,l-j\delta_{ik}},\; \forall\,(k,l)\neq (i,j).
\endaligned$$
It follows that $\widetilde{\Omega}_{ij}((\fg_{Q})_{\alpha})=(\fg_{Q})_{\varpi(\alpha)}$
for each $\alpha\in\Delta(Q)$.
\end{pf}

\subsection{} Recall from \S3.1 that for each $\vx\in\bbL$, there is
an equivalence of triangulated categories ${\bf{R}}_{\vx}: D^b(\co(\vx)\!^{\bot})\to D^b({}^{\bot}\!\co(\vx))$
which induces an equivalence of their root categories $\cR(\co(\vx)\!^{\bot})\cong \cR({}^{\bot}\!\co(\vx))$.
Since ${\bf{R}}_{\vx}$ takes exceptional objects to exceptional ones, we obtain
an isomorphism of Lie algebras
$$\Upsilon_{\vx}: \fL_{\vx}^r{\lra} \fL_{\vx}^l$$
taking
$$\b1_{\co(\vx-\vc)[\vez]}\lmto -\b1_{\co(\vx+\vc)[\vez+1]},\;
\b1_{S_{i,l_i+1}[\vez]}\lmto \b1_{\co(\vx+\vx_i)[\vez]},\;
 \b1_{S_{i,j}[\vez]}\lmto \b1_{S_{i,j}[\vez]}$$
 where $\vez\in\{0,1\}$, $1\leq i\leq t$, and $1\leq j\leq p_i$ with $j\not=l_i,l_i+1$. Finally, set
$${\widetilde \Upsilon}_{\vx}:=\theta_{\vx}^l \Upsilon_{\vx} (\theta_{\vx}^r)^{-1}:
\fL_c(\cR_{\bbX})/\fI\lra \fL_c(\cR_{\bbX})/\fI,$$
that is, we have the commutative diagram
$$\xymatrix{
  \fL_{\vx}^r \ar[r]^{\Upsilon_{\vx}}\ar[d]^{\theta_{\vx}^r} &\fL_{\vx}^l \ar[d]^{\theta_{\vx}^l}\\
  \fL_c(\cR_{\bbX})/\fI \ar[r]^{\widetilde\Upsilon_{\vx}} &\fL_c(\cR_{\bbX})/\fI}$$
% Then $\widetilde\Upsilon_{\vx}$ satisfies
%$$\aligned
%{}&\widetilde\Upsilon_{\vx}(\b1_{\co(\vx-k\vx_i)[\vez]})=-\b1_{S_{i,l_i}^{(k)}[\vez+1]},\;
%\widetilde\Upsilon_{\vx}(\b1_{S_{i,l_i+k}^{(k)}[\vez]})=\b1_{\co(\vx+k\vx_i)[\vez]},\\
%& \widetilde\Upsilon_{\vx}(\b1_{S_{i,j}^{(k)}[\vez]})=\b1_{S_{i,j}^{(k)}[\vez]}\;(j\not=l_i,l_i+k),
%\endaligned$$
% where $1\leq i\leq t$, $0\leq j\leq p_i-1$, $1\leq k\leq p_i-1$ and $\vez\in\{0,1\}$.

%\begin{Lem}\label{formula for automor Rx on spetial objects} For any $\vx$ of normal form $\vx=\sum\limits_{i=1}^tl_i\vx_i+l\vc$, we have
%  \begin{itemize}
%    \item[(1)] ${\bf{\tilde{R}}}_{\vx}(\bar{1}_{\co(\vx-k\vx_i)})
%    =\left \{ \begin{aligned} \
%               &\bar{1}_{\co(\vx)[1]}, && k=0 \\
%               & \bar{1}_{S_{i,l_i}^{(k)}[1]}, && 1\leq k\leq p_i-1;
%             \end{aligned}
%    \right.$
%    \item[(2)] ${\bf{\tilde{R}}}_{\vx}(\bar{1}_{S_{ij}})
%    =\left \{ \begin{aligned} \
%               &-\bar{1}_{\co(\vx-\vx_i)[1]}, && j=l_i \\
%               &\bar{1}_{\co(\vx+\vx_i)}, && j=l_i+1 \\
%               & \bar{1}_{S_{ij}}, && \text{else}.
%             \end{aligned}
%    \right.$
%  \end{itemize}
%\end{Lem}

\begin{Lem}\label{formula for automor Rx on spetial objects} For each $\vx\in\bbL$
with the normal form $\vx=\sum\limits_{i=1}^tl_i\vx_i+l\vc$, we have
  \begin{itemize}
    \item[(1)] ${{\wtilde \Upsilon}}_{\vx}(\bar{\b1}_{\co(\vx-k\vx_i)})
    =\left\{ \begin{array}{ll}
               -\bar{\b1}_{\co(\vx)[1]}, & k=0; \\
               -\bar{\b1}_{S_{i,l_i}^{(k)}[1]}, & 1\leq k\leq p_i-1;\\
               \bar{\b1}_{S_{i,l_i-k}^{(-k)}}, & 1\leq -k\leq p_i-1.
             \end{array}    \right.$

    \item[(2)] for $1\leq k\leq p_i-1$, ${\wtilde \Upsilon}_{\vx}(\bar{\b1}_{S_{ij}^{(k)}})
    =\left \{ \begin{array}{ll}
               -\bar{\b1}_{\co(\vx-k\vx_i)[1]}, & j=l_i; \\
               \bar{\b1}_{\co(\vx+k\vx_i)}, & j=l_i+k; \\
                \bar{\b1}_{S_{ij}^{(k)}}, & \text{otherwise.}
             \end{array}     \right.$
    \end{itemize}
\end{Lem}

\begin{pf} Both statements follow from Lemma \ref{mutation formula for special objects}. More precisely,

(1) For the case $k=0$,
$${\wtilde \Upsilon}_{\vx}(\bar{\b1}_{\co(\vx)}) =\theta_{\vx}^l{\Upsilon}_{\vx}(\b1_{\co(\vx-\vc)})
  =\theta_{\vx}^l(-\b1_{\co(\vx+\vc)[1]}) =-\bar{\b1}_{\co(\vx)[1]},$$
while for the case $1\leq k\leq p_i-1$,
$$\aligned
{}& {\wtilde \Upsilon}_{\vx}(\bar{\b1}_{\co(\vx-k\vx_i)})
  =\theta_{\vx}^l{\Upsilon}_{\vx}(\b1_{\co(\vx-k\vx_i)})
  =\theta_{\vx}^l(-\b1_{S_{i,l_i}^{(k)}[1]})
  =-\bar{\b1}_{S_{i,l_i}^{(k)}[1]}\;\text{ and }\\
&{\wtilde \Upsilon}_{\vx}(\bar{\b1}_{\co(\vx+k\vx_i)})
  =\theta_{\vx}^l{\Upsilon}_{\vx}(\b1_{\co(\vx+k\vx_i-\vc)})
  =\theta_{\vx}^l(-\b1_{S_{i,l_i}^{(p_i-k)}[1]})
  =\bar{\b1}_{S_{i,l_i+k}^{(k)}}.
\endaligned$$

(2) In case $j=l_i$, we have by Proposition \ref{special elts in cI}(5),
$$\aligned
{\wtilde \Upsilon}_{\vx}(\bar{\b1}_{S_{i,l_i}^{(k)}})
  &= -{\wtilde \Upsilon}_{\vx}(\bar{\b1}_{S_{i,l_i-k}^{(p_i-k)}[1]})
  =-\theta_{\vx}^l{\Upsilon}_{\vx}(\b1_{S_{i,l_i-k}^{(p_i-k)}[1]})\\
  &=-\theta_{\vx}^l(\b1_{\co(\vx-k\vx_i+\vc)[1]})
  =-\bar{\b1}_{\co(\vx-k\vx_i)[1]}.
  \endaligned$$
In case $j=l_i+k$,
$${\wtilde \Upsilon}_{\vx}(\bar{\b1}_{S_{i,l_i+k}^{(k)}})
  =\theta_{\vx}^l{\Upsilon}_{\vx}(\b1_{S_{i,l_i+k}^{(k)}})
  =\theta_{\vx}^l(\b1_{\co(\vx+k\vx_i)})
  =\bar{\b1}_{\co(\vx+k\vx_i)}.$$
The remaining cases are obvious.
\end{pf}

\begin{Thm}\label{braid relation of Rx} Let $\vx\in\bbL$ have the normal form
$\vx=\sum\limits_{i=1}^tl_i\vx_i+l\vc$. Then for any $1\leq k\leq t$,
$${\wtilde \Upsilon}_{\vx-\vx_k}{\wtilde \Upsilon}_{\vx}{\wtilde \Upsilon}_{\vx-\vx_k}
  ={\wtilde \Upsilon}_{\vx}{\wtilde \Upsilon}_{\vx-\vx_k}{\wtilde \Upsilon}_{\vx}.$$
\end{Thm}

\begin{pf} We show that the automorphisms on both sides take the same image on
the set of generators
$$\{\bar{\b1}_{X}, \bar{\b1}_{X[1]}\mid X=\co(\vx-\vx_k), S_{ij}\;\text{with
$1\leq i\leq t,0\leq j\leq p_i-1$, and $j\neq l_i-\dz_{i,k}$}\}.$$
Since the shift functor $[1]$ induces an isomorphism of Lie algebras, we only check
this fact for the generators $\bar{\b1}_{X}$.

By Lemma \ref{formula for automor Rx on spetial objects}, the automorphisms on both
sides preserve the above generators $\bar{\b1}_{X}$ for $X\not=\co(\vx-\vx_k), S_{k,l_k}$,
or $S_{i,l_i+1} \;(1\leq i\leq t)$. We now apply
Lemma \ref{formula for automor Rx on spetial objects} again to treat the remaining cases.

\medskip

\noindent{\bf Case 1: $X=\co(\vx-\vx_k)$.} In this case, we have
$${ \wtilde \Upsilon}_{\vx-\vx_k}{\wtilde \Upsilon}_{\vx}
  {\wtilde \Upsilon}_{\vx-\vx_k}(\bar{\b1}_{\co(\vx-\vx_k)})
  ={\wtilde \Upsilon}_{\vx-\vx_k}{\wtilde \Upsilon}_{\vx}(-\bar{\b1}_{\co(\vx-\vx_k)[1]})
  ={\wtilde \Upsilon}_{\vx-\vx_k}(\bar{\b1}_{S_{k,l_k}})
  =\bar{\b1}_{\co(\vx)}$$
  and
$${\wtilde \Upsilon}_{\vx}{\wtilde \Upsilon}_{\vx-\vx_k}
  {\wtilde \Upsilon}_{\vx}(\bar{\b1}_{\co(\vx-\vx_k)})
  ={\wtilde \Upsilon}_{\vx}{\wtilde \Upsilon}_{\vx-\vx_k}(-\bar{\b1}_{S_{k,l_k}[1]})
  ={\wtilde \Upsilon}_{\vx}(-\bar{\b1}_{\co(\vx)[1]})
  =\bar{\b1}_{\co(\vx)}.$$

\medskip

\noindent{\bf Case 2: $X=S_{i,l_i+1}$} ($i\not=k$). In this case,
$${\wtilde \Upsilon}_{\vx-\vx_k}{\wtilde \Upsilon}_{\vx}
  {\wtilde \Upsilon}_{\vx-\vx_k}(\bar{\b1}_{S_{i,l_i+1}})
  ={\wtilde \Upsilon}_{\vx-\vx_k}{\wtilde \Upsilon}_{\vx}(\bar{\b1}_{\co(\vx-\vx_k+\vx_i)})
  ={\wtilde \Upsilon}_{\vx-\vx_k}(\bar{\b1}_{\co(\vx-\vx_k+\vx_i)})
  =\bar{\b1}_{S_{i,l_i+1}}$$ and
$${\wtilde \Upsilon}_{\vx}{\wtilde \Upsilon}_{\vx-\vx_k}
  {\wtilde \Upsilon}_{\vx}(\bar{\b1}_{S_{i,l_i+1}})
  ={\wtilde \Upsilon}_{\vx}{\wtilde \Upsilon}_{\vx-\vx_k}(\bar{\b1}_{\co(\vx+\vx_i)})
  ={\wtilde \Upsilon}_{\vx}(\bar{\b1}_{\co(\vx+\vx_i)})
  =\bar{\b1}_{S_{i,l_i+1}}.$$

\medskip

\noindent{\bf Case 3: $X=S_{k,l_k}$.} Then
$${\wtilde \Upsilon}_{\vx-\vx_k}{\wtilde \Upsilon}_{\vx}
  {\wtilde \Upsilon}_{\vx-\vx_k}(\bar{\b1}_{S_{k,l_k}})
  ={\wtilde \Upsilon}_{\vx-\vx_k}{\wtilde \Upsilon}_{\vx}(\bar{\b1}_{\co(\vx)})
  ={\wtilde \Upsilon}_{\vx-\vx_k}(-\bar{\b1}_{\co(\vx)[1]})
  =-\bar{\b1}_{S_{k,l_k}[1]}$$ and
$${\wtilde \Upsilon}_{\vx}{\wtilde \Upsilon}_{\vx-\vx_k}
  {\wtilde \Upsilon}_{\vx}(\bar{\b1}_{S_{k,l_k}})
  ={\wtilde \Upsilon}_{\vx}{\wtilde \Upsilon}_{\vx-\vx_k}(-\bar{\b1}_{\co(\vx-\vx_k)[1]})
  ={\wtilde \Upsilon}_{\vx}(\bar{\b1}_{\co(\vx-\vx_k)})
  =-\bar{\b1}_{S_{k,l_k}[1]}.$$

\medskip

\noindent{\bf Case 4: $X=S_{k,l_k+1}$.} We consider two subcases $p_k=2$ and $p_k>2$. If $p_k=2$, then
$$\aligned
{}&{\wtilde \Upsilon}_{\vx-\vx_k}{\wtilde \Upsilon}_{\vx}
  {\wtilde \Upsilon}_{\vx-\vx_k}(\bar{\b1}_{S_{k,l_k+1}})
  ={\wtilde \Upsilon}_{\vx-\vx_k}{\wtilde \Upsilon}_{\vx}{\wtilde \Upsilon}_{\vx-\vx_k}(-\bar{\b1}_{S_{k,l_k}[1]})\\
   =&{\wtilde \Upsilon}_{\vx-\vx_k}{\wtilde \Upsilon}_{\vx}(\bar{\b1}_{\co(\vx)[1]})
  ={\wtilde \Upsilon}_{\vx-\vx_k}(-\bar{\b1}_{\co(\vx)})=\bar{\b1}_{S_k,l_k}
\endaligned$$
  and
$$\aligned
{}&{\wtilde \Upsilon}_{\vx}{\wtilde \Upsilon}_{\vx-\vx_k}
  {\wtilde \Upsilon}_{\vx}(\bar{\b1}_{S_{k,l_k+1}})
  ={\wtilde \Upsilon}_{\vx}{\wtilde \Upsilon}_{\vx-\vx_k}(\bar{\b1}_{\co(\vx+\vx_k)})={\wtilde \Upsilon}_{\vx}{\wtilde \Upsilon}_{\vx-\vx_k}(\bar{\b1}_{\co(\vx-\vx_k+\vc)})\\
=&{\wtilde \Upsilon}_{\vx}{\wtilde \Upsilon}_{\vx-\vx_k}(\bar{\b1}_{\co(\vx-\vx_k)})
  ={\wtilde \Upsilon}_{\vx}(-\bar{\b1}_{\co(\vx-\vx_k)[1]})
  =\bar{\b1}_{S_{k,l_k}}.
  \endaligned$$
  If $p_k\geq 3$, then
$${\wtilde \Upsilon}_{\vx-\vx_k}{\wtilde \Upsilon}_{\vx}
  {\wtilde \Upsilon}_{\vx-\vx_k}(\bar{\b1}_{S_{k,l_k+1}})
  ={\wtilde \Upsilon}_{\vx-\vx_k}{\wtilde \Upsilon}_{\vx}(\bar{\b1}_{S_{k,l_k+1}})
  ={\wtilde \Upsilon}_{\vx-\vx_k}(\bar{\b1}_{\co(\vx+\vx_k)})
  =\bar{\b1}_{S_{k,l_k+1}^{(2)}}$$
  and
$${\wtilde \Upsilon}_{\vx}{\wtilde \Upsilon}_{\vx-\vx_k}
  {\wtilde \Upsilon}_{\vx}(\bar{\b1}_{S_{k,l_k+1}})
  ={\wtilde \Upsilon}_{\vx}{\wtilde \Upsilon}_{\vx-\vx_k}(\bar{\b1}_{\co(\vx+\vx_k)})
  ={\wtilde \Upsilon}_{\vx}(\bar{\b1}_{S_{k,l_k+1}^{(2)}})
  =\bar{\b1}_{S_{k,l_k+1}^{(2)}}.$$

In conclusion, we have in all the cases the equality
$${\wtilde \Upsilon}_{\vx-\vx_k}{\wtilde \Upsilon}_{\vx}{\wtilde \Upsilon}_{\vx-\vx_k}(\bar{\b1}_{X})
  ={\wtilde \Upsilon}_{\vx}{\wtilde \Upsilon}_{\vx-\vx_k}{\wtilde \Upsilon}_{\vx}(\bar{\b1}_{X}).$$
 This finishes the proof.
\end{pf}

\subsection{} For each exceptional object $X\in\cR_{\bbX}$,
the adjoint operator $\ad(\bar{\b1}_{X})$ on $\fL_c(\cR_{\bbX})/\fI$ is locally nilpotent.
Thus, we can define its exponential map
\[\exp(\ad(\bar{\b1}_{X}))={\rm id}+\ad(\bar{\b1}_{X})+\frac{\b1}{2!}(\ad(\bar{\b1}_{X}))^2+\cdots\]

\begin{Lem}\label{formula for adjoint operator on spetial objects1}
Let $X\in\cR_{\bbX}$ be an exceptional object. Then
\begin{itemize}
  \item[(1)] $\exp(\ad(\bar{\b1}_{X}))(\bar{\b1}_{X})=\bar{\b1}_{X};$
  \item[(2)] $\exp(\ad(\bar{\b1}_{X}))(h_{\overline{[X]}})=h_{\overline{[X]}}+2\bar{\b1}_{X};$
  \item[(3)] $\exp(\ad(\bar{\b1}_{X}))(\bar{\b1}_{X[1]})=\bar{\b1}_{X[1]}+h_{\overline{[X]}}+\bar{\b1}_{X};$
\end{itemize}
\end{Lem}

\begin{pf} Since $[\bar{\b1}_{X}, \bar{\b1}_{X}]=0$, $[\bar{\b1}_{X}, h_{\overline{[X]}}]=(\bar{[X]},\bar{[X]})\bar{\b1}_{X}=2\bar{\b1}_{X}$ and
$[\bar{\b1}_{X}, \bar{\b1}_{X[1]}]=h_{\overline{[X]}}$, statements (1) and (2) follow immediately.

For statement (3), we have
$$\aligned
 {} &\exp(\ad(\bar{\b1}_{X}))(\bar{\b1}_{X[1]})
  =\bar{\b1}_{X[1]}+[\bar{\b1}_{X},\bar{\b1}_{X[1]}]+\frac{\b1}{2}[\bar{\b1}_{X},[\bar{\b1}_{X},\bar{\b1}_{X[1]}]]\\
  =&\bar{\b1}_{X[1]}+h_{\overline{[X]}}+\frac{1}{2}[\bar{\b1}_{X}, h_{\overline{[X]}}]
  =\bar{\b1}_{X[1]}+h_{\overline{[X]}}+\bar{\b1}_{X}.
  \endaligned$$
\end{pf}

\begin{Lem}\label{expexpexp for exc objs} For each exceptional object $X\in\cR_{\bbX}$, the operator $$\exp(\ad(\bar{\b1}_{X}))\exp(\ad(\bar{\b1}_{X[1]}))\exp(\ad(\bar{\b1}_{X}))$$
exchanges the elements $\bar{\b1}_{X}$ and $\bar{\b1}_{X[1]}$.
\end{Lem}

\begin{pf} By Lemma \ref{formula for adjoint operator on spetial objects1}, we have the equalities
  \begin{equation}
 \begin{split}
&\exp(\ad(\bar{\b1}_{X}))\exp(\ad(\bar{\b1}_{X[1]}))
\exp(\ad(\bar{\b1}_{X}))(\bar{\b1}_{X})\\
=&\exp(\ad(\bar{\b1}_{X}))\exp(\ad(\bar{\b1}_{X[1]}))(\bar{\b1}_{X})\\
=&\exp(\ad(\bar{\b1}_{X}))(\bar{\b1}_{X}-h_{\overline{[X]}}+\bar{\b1}_{X[1]})\\
=&\bar{\b1}_{X}-(h_{\overline{[X]}}+2\bar{\b1}_{X})
+(\bar{\b1}_{X}+h_{\overline{[X]}}+\bar{\b1}_{X[1]})\\
=&\bar{\b1}_{X[1]}.
 \end{split}
 \end{equation}
Similarly, we have
\begin{equation}
\exp(\ad(\bar{\b1}_{X}))\exp(\ad(\bar{\b1}_{X[1]}))
\exp(\ad(\bar{\b1}_{X}))(\bar{\b1}_{X[1]})
=\bar{\b1}_{X}.
\end{equation}
\end{pf}

 By Lemma \ref{Lie bracket} we can easily obtain the following result.

\begin{Lem}\label{formula for adjoint operator on spetial objects2} For each $\vx$ of normal form $\vx=\sum\limits_{i=1}^tl_i\vx_i+l\vc$ and $1\leq k\leq p_i-1$, the following hold in $\fL_c(\cR_{\bbX})/\fI$:
\begin{itemize}
  \item[(1)] $\exp(\ad(\bar{\b1}_{\co(\vx)}))(\bar{\b1}_{S_{ij}^{(k)}})
      =\bar{\b1}_{S_{ij}^{(k)}}-\delta_{j,l_i+k}\bar{\b1}_{\co(\vx+k\vx_i)};$
  \item[(2)] $\exp(\ad(\bar{\b1}_{\co(\vx)}))(\bar{\b1}_{S_{ij}^{(k)}[1]})
      =\bar{\b1}_{S_{ij}^{(k)}[1]}+\delta_{j,l_i}\bar{\b1}_{\co(\vx-k\vx_i)};$
  \item[(3)] $\exp(\ad(\bar{\b1}_{\co(\vx)}))(\bar{\b1}_{\co(\vx-k\vx_i)[1]})
      =\bar{\b1}_{\co(\vx-k\vx_i)[1]}-\bar{\b1}_{S_{i,l_i}^{(k)}};$
  \item[(4)] $\exp(\ad(\bar{\b1}_{\co(\vx)}))(\bar{\b1}_{\co(\vx+k\vx_i)[1]})
      =\bar{\b1}_{\co(\vx+k\vx_i)[1]}+\bar{\b1}_{S_{i,l_i+k}^{(k)}[1]}.$
\end{itemize}
\end{Lem}

Recall from Corollary \ref{basis of quotient algebra} that, for each $\vx\in\bbL$ of the normal form $\vx=\sum\limits_{i=1}^tl_i\vx_i+l\vc$, the quotient algebra $\fL_c(\cR_{\bbX})/\fI$ admits a set of
Chevalley generators
$$\bar{\b1}_{X}, \bar{\b1}_{X[1]}, h_{\overline{[X]}},\; X\in\{\co(\vx), S_{ij}\mid 1\leq i\leq t, 0\leq j\leq p_i-1,j\neq l_i\}.$$
Consequently, there is an automorphism
$$\cE_{\vx}: \fL_c(\cR_{\bbX})/\fI\lra \fL_c(\cR_{\bbX})/\fI$$
 taking $\bar{\b1}_{\co(\vx)[\vez]}\mapsto -\bar{\b1}_{\co(\vx)[\vez]}$,
$\bar{\b1}_{S_{i,l_i+1}[\vez]}\mapsto -\bar{\b1}_{S_{i,l_i+1}[\vez]}$ ($\vez\in\{0,1\}$,
$1\leq i\leq t$) and fixing all the other generators.

In the following we show that for each $\vx\in \bbL$, the automorphism
$${\wtilde \Upsilon}_{\vx}:\fL_c(\cR_{\bbX})/\fI\lra \fL_c(\cR_{\bbX})/\fI,$$
can be expressed as a product of certain Tits' automorphisms (up to sign).

\begin{Prop}\label{theorem for Rx} For each $\vx\in\bbL$ of the normal form
$\vx=\sum\limits_{i=1}^tl_i\vx_i+l\vc$, the equality
$${\wtilde \Upsilon}_{\vx}\cE_{\vx}=
\exp(\ad(\bar{\b1}_{\co(\vx)}))\exp(\ad(\bar{\b1}_{\co(\vx)[1]}))\exp(\ad(\bar{\b1}_{\co(\vx)}))$$
holds in $\Aut(\fL_c(\cR_{\bbX})/\fI)$.
\end{Prop}

\begin{pf} We show that the operators on both sides take the same values
on the set of generators
$$\{\bar{\b1}_{X}, \bar{\b1}_{X[1]}\mid  X=\co(\vx), S_{ij}\;\text{with}\; 1\leq i\leq t, 0\leq j\leq p_i-1,j\neq l_i\}.$$
As above, we only consider the generators $\bar{\b1}_{X}$.

\medskip

 \noindent {\bf Case 1:} $X=\co(\vx)$. By Lemma
\ref{formula for automor Rx on spetial objects}(1) and Lemma \ref{expexpexp for exc objs}, we have
   $${\wtilde \Upsilon}_{\vx}\cE_\vx(\bar{\b1}_{\co(\vx)})=\bar{\b1}_{\co(\vx)[1]}
   =\exp(\ad(\bar{\b1}_{\co(\vx)}))\exp(\ad(\bar{\b1}_{\co(\vx)[1]}))
      \exp(\ad(\bar{\b1}_{\co(\vx)}))(\bar{\b1}_{\co(\vx)}).$$

\medskip

\noindent {\bf Case 2:} $X=S_{ij}\; (j\neq l_i+1)$. It is easy to see that
the automorphisms on both sides act as identities on $\bar{\b1}_{S_{ij}}$.

\medskip

\noindent {\bf Case 3:} $X=S_{i,l_i+1} $. By Lemma \ref{formula for automor Rx on spetial objects}(2),
  $${\wtilde \Upsilon}_{\vx}\cE_\vx(\bar{\b1}_{S_{i,l_i+1}})=-\bar{\b1}_{\co(\vx+\vx_i)}.$$
On the other hand,
   \begin{equation}
 \begin{split}
  &\exp(\ad(\bar{\b1}_{\co(\vx)}))\exp(\ad(\bar{\b1}_{\co(\vx)[1]}))
  \exp(\ad(\bar{\b1}_{\co(\vx)}))(\bar{\b1}_{S_{i,l_i+1}})\\
  =&\exp(\ad(\bar{\b1}_{\co(\vx)}))\exp(\ad(\bar{\b1}_{\co(\vx)[1]}))(\bar{\b1}_{S_{i,l_i+1}}-\bar{\b1}_{\co(\vx+\vx_i)})\\
  =&\exp(\ad(\bar{\b1}_{\co(\vx)}))(\bar{\b1}_{S_{i,l_i+1}}-(\bar{\b1}_{\co(\vx+\vx_i)}+\bar{\b1}_{S_{i,l_i+1}}))\\
  =&\exp(\ad(\bar{\b1}_{\co(\vx)}))(-\bar{\b1}_{\co(\vx+\vx_i)})\\
  =&-\bar{\b1}_{\co(\vx+\vx_i)}.
 \end{split}
 \end{equation}
\end{pf}

Fix $\vx$ of the normal form
$\vx=\sum\limits_{i=1}^tl_i\vx_i+l\vc$ and $1\leq k\leq t$. By Corollary
\ref{basis of quotient algebra}, the elements ${\bf 1}_{X}, {\bf 1}_{X[1]}, h_{[X]}$
form a complete set of Chevalley generators of $\fL_c(\cR_{\bbX})/\fI$,
where $$X\in \{\co(\vx-\vx_k), S_{ij}\mid 1\leq i\leq t, 0\leq j\leq p_i-1,j\neq l_i-\dz_{i,k}\}=:{\scr X}.$$
 Further, define an automorphism
$$\cE_{\vx,k}: \fL_c(\cR_{\bbX})/\fI\lra \fL_c(\cR_{\bbX})/\fI$$
 which takes $\bar{\b1}_{S_{k,j}[\vez]}\mapsto -\bar{\b1}_{S_{k,j}[\vez]}$ ($\vez\in\{0,1\}$,
$j=l_k, l_k+1$) and fixes all the other generators.

\begin{Prop}\label{Tit's form for RxRxRx} Let $\vx\in\bbL$ be of the normal
form $\vx=\sum\limits_{i=1}^tl_i\vx_i+l\vc$. Then for any $1\leq k\leq t$,
 $$ {\wtilde \Upsilon}_{\vx-\vx_k}{\wtilde \Upsilon}_{\vx}{\wtilde \Upsilon}_{\vx-\vx_k}\cE_{\vx,k}
  =\exp(\ad(\bar{\b1}_{S_{k,l_k}}))\exp(\ad(\bar{\b1}_{S_{k,l_k}[1]}))\exp(\ad(\bar{\b1}_{S_{k,l_k}})).$$
\end{Prop}

\begin{pf} We show that the automorphisms on both sides take
the same values for the generators $\bar{\b1}_{X}, \bar{\b1}_{X[1]}, X\in \scr X$.
For simplicity, we denote by $\cE$ the automorphism on the right hand side.
We first calculate the images of $\cE$ on the set of generators $\bar{\b1}_{X}$
given above. The calculation for $\bar{\b1}_{X[1]}$ is analogous.
By the definition, for $X\in{\scr X}\backslash\{\co(\vx-\vx_k), S_{k,l_k}, S_{k,l_k+1}\}$,
$\cE(\bar{\b1}_{X})=\bar{\b1}_{X}$. In the following we calculate $\cE(\bar{\b1}_{X})$
for the remaining $X$.

\medskip

\noindent {\bf Case 1: $X=\co(\vx-\vx_k)$}. In this case,
$$\aligned
\cE(\bar{\b1}_{\co(\vx-\vx_k)})&=\exp(\ad(\bar{\b1}_{S_{k,l_k}}))\exp(\ad(\bar{\b1}_{S_{k,l_k}[1]}))
(\bar{\b1}_{\co(\vx-\vx_k)}-\bar{\b1}_{\co(\vx)})\\
&=\exp(\ad(\bar{\b1}_{S_{k,l_k}}))(\bar{\b1}_{\co(\vx-\vx_k)}+(\bar{\b1}_{\co(\vx)}-\bar{\b1}_{\co(\vx-\vx_k)}))\\
&=\exp(\ad(\bar{\b1}_{S_{k,l_k}}))(\bar{\b1}_{\co(\vx)})=\bar{\b1}_{\co(\vx)}.
\endaligned$$

\medskip

\noindent {\bf Case 2: $X=S_{k,l_k}$}. By Lemma \ref{expexpexp for exc objs},
$$\cE(\bar{\b1}_{S_{k,l_k}})=\bar{\b1}_{S_{k,l_k}[1]}.$$

\medskip

\noindent {\bf Case 3: $X=S_{k,l_k+1}$}. If $p_k=2$, we have by Lemma \ref{expexpexp for exc objs},
$$\cE(\bar{\b1}_{S_{k,l_k+1}})=\cE(-\bar{\b1}_{S_{k,l_k}[1]})=-\bar{\b1}_{S_{k,l_k}}.$$
 If $p_k\geq 3$, then
 $$\aligned
\cE(\bar{\b1}_{S_{k,l_k+1}})&=\exp(\ad(\bar{\b1}_{S_{k,l_k}}))\exp(\ad(\bar{\b1}_{S_{k,l_k}[1]}))
(\bar{\b1}_{S_{k,l_{k+1}}}-\bar{\b1}_{S_{k,l_{k+1}}^{(2)}})\\
&=\exp(\ad(\bar{\b1}_{S_{k,l_k}}))(S_{k,l_{k+1}}-(\bar{\b1}_{S_{k,l_{k+1}}^{(2)}}+
\bar{\b1}_{S_{k,l_{k+1}}}))\\
&=\exp(\ad(\bar{\b1}_{S_{k,l_k}}))(-\bar{\b1}_{S_{k,l_{k+1}}^{(2)}})=-\bar{\b1}_{S_{k,l_{k+1}}^{(2)}}.
\endaligned$$

A comparison with the calculations given in the proof of Theorem \ref{braid relation of Rx}
implies the desired equality.
\end{pf}

Finally, for each $\vx\in\bbL$, set
$$\Xi_{\vx}:=\Phi^{-1} {\wtilde \Upsilon}_{\vx} \Phi:\fg_Q\lra\fg_Q.$$
In other words, we have the following commutative diagram:
$$\xymatrix{
  \fg_Q \ar[r]^{\Xi_{\vx}}\ar[d]^{\Phi} &\fg_Q \ar[d]^{\Phi}\\
  \fL_c(\cR_{\bbX})/\fI \ar[r]^{{\wtilde \Upsilon}_{\vx}} &\fL_c(\cR_{\bbX})/\fI}
  $$
As a consequence of Propositions \ref{theorem for Rx} and \ref{Tit's form for RxRxRx}, we obtain the main result of this section.

\begin{Thm}\label{corollary for Rx} The following equalities hold in  $\Aut(\fg_Q)$, up to sign:
$$\Xi_{0}=\exp(\ad(e_{\star}))\exp(\ad(-f_{\star}))\exp(\ad(e_{\star}))$$
and for each vertex $\oz=\oz_{ij}$ with $(i,j)\in \whI$,
$$\Xi_{(j-1)\vx_i} \Xi_{j\vx_i} \Xi_{(j-1)\vx_i}
=\exp(\ad(e_\oz))\exp(\ad(-f_\oz))\exp(\ad(e_\oz)).$$
\end{Thm}

\bigskip

\section{Double Ringel--Hall algebra of $\bbX$ and Lusztig's symmetries}

As in the previous sections, we fix a finite field $\bbF$ with $q$ elements
and set $v=\sqrt{q}$. Let $\bbX=\bbX_\bbF$ be the weighted projective line over $
\bbF$ associated with $\bf p$ and ${\boldsymbol\lambda}$ and let $Q=(I,Q_1)$ be the
associated star shaped quiver given in \S2.3. In this section, we deal with the (reduced Drinfeld)
double composition algebras of subcategories $D^b(\co(\vx)\!^{\bot})$ and $D^b({}^{\bot}\!\co(\vx))$
of $D^b(\coh\bbX_\bbF)$ both of which are naturally
isomorphic to the quantum enveloping algebra ${\bf U}_v(\fg_Q)$ by specializing $\bfv$ to
the square root $v=\sqrt{q}$. By applying a theorem of Cramer \cite{cram},
mutation functors between $D^b(\co(\vx)\!^{\bot})$ and $D^b({}^{\bot}\!\co(\vx))$ induce
isomorphisms of their associated double composition algebras. We finally show that
 these isomorphisms indeed provide a realization of Lusztig's symmetries of $\bfU_v(\fg_Q)$.

\subsection{}

By definition, the Ringel--Hall algebra $H(\bbX)$ of $\bbX$ is the complex
vector space with basis $u_{[X]}$ with $[X]$ running through all the isoclasses of coherent sheaves
in $\coh\bbX$. However, we will write $u_X=u_{[X]}$ for notational simplicity.
The multiplication is given by
$$u_Mu_N=v^{\langle [M], [N]\rangle}\sum_{[L]}F_{M ,N}^L u_L,$$
 where the sum is taken over all the isoclasses $[L]$ in $\coh\bbX$, and $F^L_{M,N}$
denotes the number of subsheaves $X$ of $L$ such that $L\cong N$ and
$L/X\cong M$. Following \cite{X97} (see also \cite{DJX}), we can define the (reduced Drinfeld)
double Ringel--Hall algebra ${\cal D}H(\bbX)$ of $\bbX$ which admits a triangular decomposition
$$\cD H(\bbX)=H^+(\bbX)\otimes_\bbC H^0(\bbX)\otimes_\bbC H^-(\bbX),$$
 where $H^0(\bbX)=\bbC[K_\az\mid \az\in K_0(\coh\bbX)]$ and $H^+(\bbX)$ (resp., $H^-(\bbX)$)
has a basis $\{u_X^+\}$ (resp., $\{u_X^-\}$) such that $H(\bbX)\cong H^+(\bbX)$ taking $u_X\mapsto u_X^+$
(resp., $H(\bbX)\cong H^-(\bbX)$ taking $u_X\mapsto u_X^-$). In particular,
$$K_\az u_X^{\pm}=v^{\pm (\az,[X])}u_X^\pm K_\az,\;
[u_{S_{ij}}^+, u_{S_{ij}}^-]=\frac{K_{[S_{ij}]}-K_{[S_{ij}]}^{-1}}{1-v^2}$$
 where $\az\in K_0(\coh\bbX)$, $X\in\coh\bbX$, $1\leq i\leq t$ and $0\leq j\leq p_i-1$.

In general, for each hereditary abelian
$\bbF$-category $\cal A$ with finite dimensional Hom-spaces and Ext-spaces, we can define similarly
the Ringel--Hall algebra $H(\cal A)$ of $\cal A$ and its (reduced Drinfeld) double $\cD H(\cal A)$.

\subsection{}
By a result of Cramer \cite{cram}, the mutation functor
$${\bf{R}}_{\vx}: D^b(\co(\vx)\!^{\bot})\lra D^b({}^{\bot}\!\co(\vx))$$
defined in Section 3 induces an isomorphism between the double Ringle--Hall algebras
$$\wtilde{R}_{\vx}: \cD H(\co(\vx)\!^{\bot})\lra \cD H({}^{\bot}\!\co(\vx)).$$
 Now let $\cD C(\co(\vx)\!^{\bot})$ (resp., $\cD C({}^{\bot}\!\co(\vx))$) be the subalgebra
of $\cD H(\co(\vx)\!^{\bot})$ (resp., $\cD H({}^{\bot}\!\co(\vx))$) generated by $u_X^\pm$
and $K_\az$ with $X$ exceptional in $\co(\vx)\!^{\bot}$ (resp., ${}^{\bot}\!\co(\vx)$)
and $\az\in K_0(\co(\vx)\!^{\bot})$ (resp., $\az\in K_0({}^{\bot}\!\co(\vx))$), called
the composition subalgebra. By \cite[Cor.~3.6]{RuanZh}, $\cD C(\co(\vx)\!^{\bot})$
(resp., $\cD C({}^{\bot}\!\co(\vx))$) is generated by $K_{[X_i]}^{\pm1}$ and $u_{X_i}^{\pm} (1\leq i\leq n)$
for each complete exceptional sequence $(X_1,\ldots, X_n)$ in $\co(\vx)\!^{\bot}$
(resp., ${}^{\bot}\!\co(\vx)$).

More precisely, in case $\vx=0$, $\cD C(\co\!^{\bot})$ is generated by $u_{\co(-\vc)}^{\pm}$,
$u_{S_{ij}}^{\pm}$ and $K_\az$ for $(i,j)\in \whI$ and $\az\in K_0(\co\!^{\bot})$.
Then $\wtilde{R}_{0}$ restricts to an isomorphism
$$\wtilde{R}_{0}: \cD C(\co\!^{\bot})\lra \cD C({}^{\bot}\!\co)$$
which, by \cite[Prop.~5]{cram}, takes $K_\az\mapsto K_{R_0(\az)}$ for $\az\in K_0(\co\!^{\bot})$,
$$u_{\co(-\vc)}^{\pm}\lmto -v^{-1}u_{\co(\vc)}^{\mp}K_{[\co(\vc)]}^{\pm1},\;\;
u_{S_{i1}}^{\pm}\lmto u_{\co(\vx_i)}^{\pm},\;\forall\,1\leq i\leq t,$$
 and preserves all other generators. For $\vx=j\vx_i$ with $(i,j)\in\whI$,
 $\cD C(\co(j\vx_i)\!^{\bot})$ is generated by $u_{\co}^{\pm}, u_{S_{i0}}^{\pm}$,
$u_{S_{kl}}^{\pm}$, and $K_\az$ for $\oz_{kl}\in {\whI}\backslash\{(i,j)\}$
and $\az\in K_0(\co(j\vx_i)\!^{\bot})$. Then
$\wtilde{R}_{j\vx_i}$ restricts to an isomorphism
$$\wtilde{R}_{j\vx_i}: \cD C(\co(j\vx_i)\!^{\bot})\lra \cD C({}^{\bot}\!\co(j\vx_i))$$
which takes $K_\az\mapsto K_{R_{j\vx_i}(\az)}$ for $\az\in K_0(\co(j\vx_i)\!^{\bot})$,
$$ u_{\co}^{\pm}\lmto -v^{-1}u^{\mp}_{S_{ij}^{(j)}}K^{\pm1}_{[S_{ij}^{(j)}]},\;
u_{S_{k,\dz_{i,k}j+1}}^{\pm}\lmto u_{\co(j\vx_i+\vx_k)}^\pm,\;\forall\,1\leq k\leq t,$$
 and preserves all other generators.

\subsection{} Let $\bfU_v(\cL\fg_Q)$ be the quantum loop algebra of $\fg_Q$ which, by definition,
is an algebra over $\bbC$ generated by $x^{\pm}_{\oz,k}$, $h_{\oz,l}$, and $K^{\pm1}_\oz$ for
$\oz\in I$, $k\in\bbZ$ and $l\in\bbZ^\times$ subject to certain relations; see for example,
\cite[\S1.8]{Sch04}. Also, recall from \S2.2 that the quantum enveloping algebra ${\bf U}_v(\fg_Q)$ is
the $\bbC$-algebra generated by $E_\oz,F_\oz,K_\oz^{\pm 1}$ ($\oz\in I$).

\begin{Lem}\label{a map from Gamma to Z induces an iso of quantum group}
Each map $\sigma: I\to \bbZ$ induces an automorphism
$$\xi_{\sigma}: \bfU_v(\cL\fg_Q)\lra \bfU_v(\cL\fg_Q),\; K^{\pm1}_{\oz}\lmto K^{\pm1}_{\oz},\; h_{\oz,r}
\lmto h_{\oz,r},\; x_{\oz,k}^{\pm}\lmto x_{\oz,k\pm\sigma(\oz)}^{\pm}$$
for all $\oz\in I, r\in \bbZ^\times, k\in\bbZ$. In particular, the subalgebra
generated by the $x_{\oz,\pm\sigma(\oz)}^{\pm}$ and $K_\oz^{\pm1}$ is
isomorphic to $\bfU_v(\fg_Q)$ via taking
$$K^{\pm1}_{\oz}\lmto K^{\pm1}_{\oz},\; x_{\oz,\sigma(\oz)}^+
\lmto E_\oz,\; x_{\oz,-\sigma(\oz)}^-\lmto F_\oz.$$
\end{Lem}

\begin{pf} The fact that $\xi_{\sigma}$ is an automorphism follows from
a direct checking (but tedious) on the generating relations of $\bfU_v(\cL\fg_Q)$.
This induces an isomorphism between the two subalgebras of $\bfU_v(\cL\fg_Q)$
generated respectively by the $x_{\oz,0}^{\pm}$ and by the $x_{\oz,\pm\sigma(\oz)}^{\pm}$ together
with the $K_\oz^{\pm1}$, where the former
one is obviously isomorphic to $\bfU_v(\fg_Q)$.

\end{pf}

We now introduce some notation for later use. Let $A$ be an algebra over $\bbC$.
For $x,y\in A$ and $0\not=c\in\bbC$, set
$$[x,y]_c=xy-cyx\in A.$$
Then $[x,y]_c=-c[y,x]_{c^{-1}}$. Further, for $x_1,\ldots,x_n\in A$ with $n\geq 3$,
we define recursively
$$\aligned
{}& [x_1,x_2,\cdots,x_n]_{v^{-1}}=[[x_1,\ldots,x_{n-1}]_{v^{-1}},x_n]_{v^{-1}}\;\text{ and }\;\\
& [x_1,x_2\cdots,x_n]_v=[[x_1,[x_2\cdots,x_n]_v]_{v}.
\endaligned$$
Then by an inductive argument, we have
$$[x_1,x_2,\cdots,x_n]_{v^{\pm 1}}=(-v^{\pm 1})^{n-1}[x_n,\cdots,x_2,x_1]_{v^{\mp 1}}.$$
The following result can be checked by a direct calculation.

\begin{Lem}\label{skew commutator} {\rm(1)} $[x, [y,z]_{v^{\pm 1}}]_{v^{\pm 1}}
=[[x,y]_{v^{\pm 1}},z]_{v^{\pm 1}}+v^{\pm 1}[y, [x,z]]$;

{\rm(2)} $[x, [y,z]_{v^{\pm 1}}]=[[x,y],z]_{v^{\pm 1}}+[y, [x,z]]_{v^{\pm 1}}$.
\end{Lem}

\begin{Rem} As a consequence of (1), if $[x,z]=0$, then
$$[x, [y,z]_{v^{\pm 1}}]_{v^{\pm 1}}=[[x,y]_{v^{\pm 1}},z]_{v^{\pm 1}}.$$
More generally, if $x_1,\ldots,x_n$ satisfy $[x_{i}, x_{j}]=0$ for $|i-j|\geq 2$, then
$$[[x_1,\ldots,x_{n-1}]_{v^{\pm1}},x_n]_{v^{\pm1}}=[x_1,x_2,\cdots,x_n]_{v^{\pm 1}}
=[x_1,[x_2\cdots,x_n]_{v^{\pm1}}]_{v^{\pm1}}.$$
\end{Rem}

\medskip

By \cite[Thm.~5.2]{Sch04} and \cite[Thm.~5.5]{DJX}, there is a surjective algebra homomorphism
$$\Xi: \bfU_v(\cL\fg_Q)\lra \cD C(\coh\bbX),$$
 where $\cD C(\coh\bbX)$ denotes the composition subalgebra of ${\cal D}H(\bbX)$.
 We refer to \cite[Thm.~5.5]{DJX} for the precise definition of $\Xi$.
 We remark that for any $\alpha\in\Gamma$ and $\lambda_{\alpha}\in\bbZ$,
if we replace the generators $K_{\alpha}$ of $\cD C(\coh\bbX)$ by $K_{\alpha+\lambda_{\alpha}\delta}$,
we still obtain a surjection. In the following, we will identify the following groups for any $\vx\in\bbL$: $$K_0(\Gamma)=K_0(\coh\bbX)/\bbZ\delta=K_0(\co(\vx)^{\perp})=K_0({}^{\perp}\co(\vx)).$$

For each $(i,j)\in \whI$, we simply write $u_{ij}^\pm=u_{S_{ij}}^\pm$ and $K_{ij}=K_{[S_{ij}]}$.

\begin{Lem}\label{cancel ui} {\rm (1)} Let $X,Y\in \coh\bbX$ and $\alpha, \beta\in K_0(\coh\bbX)$. Then
for each $a\in\bbZ$,
$$\aligned
{}&[u_{X}^\pm K_{\alpha}, u_{Y}^\pm K_{\beta}]_{v^{a}}
=v^{(\alpha, [Y])}[u_{X}^\pm, u_{Y}^\pm]_{v^{a'}}K_{\alpha}K_{\beta}\;\text{ and}\\
&[u_{X}^\pm K_{\alpha}, u_{Y}^\mp K_{\beta}]_{v^{a}}
=v^{-(\alpha, [Y])}[u_{X}^\pm, u_{Y}^\mp]_{v^{a''}}K_{\alpha}K_{\beta},
\endaligned$$
where $a'=a\pm([X],\beta)\mp(\alpha,[Y])$ and $a''=a\pm([X],\beta)\pm(\alpha,[Y])$.

{\rm (2)} If $([S_{ij}],[X])=-1$ and $[u_{ij}^{\pm},u_{X}^{\mp}]=0$ for some $(i,j)\in\whI$
and $X\in \coh\bbX$, then
$$[u_{ij}^{-}, [u_X^{+}, u_{ij}^{+}]_{v^{-1}}]=v^{-2}u_X^{+}K_{ij} \text{\quad and\quad }
[[u_{ij}^{-}, u_{X}^{-}]_v, u_{ij}^{+}]=v^{-1}u_{X}^{-}K^{-1}_{ij}.$$
\end{Lem}

\begin{pf} The statement (1) can be checked directly. We now prove (2).

By Lemma \ref{skew commutator}, we have $$[u_{ij}^{-}, [u_X^{+}, u_{ij}^{+}]_{v^{-1}}]
=[[u_{ij}^{-}, u_{X}^{+}], u_{ij}^{+}]_{v^{-1}}+[[u_X^{+}, [u_{ij}^{-}, u_{ij}^{+}]]_{v^{-1}}
=[[u_X^{+}, [u_{ij}^{-}, u_{ij}^{+}]]_{v^{-1}}.$$
 Recall that $ [u_{ij}^{-}, u_{ij}^{+}]=\frac{K_{ij}^{-1}-K_{ij}}{1-v^2}$.
By $([S_{ij}],[X])=-1$, we have
$$[u_X^{+},K_{ij}^{-1}]_{v^{-1}}=0\;\text{ and }\;[u_X^{+},K_{ij}]_{v^{-1}}=(v-v^{-1})K_{ij}u_X^{+}.$$
It follows that
$[u_{ij}^{-}, [u_X^{+}, u_{ij}^{+}]_{v^{-1}}]=v^{-1}K_{ij}u_X^{+}$.

Similarly, we get
$$\aligned
{[{[u_{ij}^-, u_X^-]}_v, u_{ij}^+]}&=-[[u_{ij}^{+}, u_{ij}^{-}],u_X^{-}]_{v}
=-[\frac{K_{ij}-K_{ij}^{-1}}{1-v^2},u_X^{-}]_{v}\\
&=\frac{\b1}{1-v^2}[K_{ij}^{-1}, u_X^{-}]_{v}
=v^{-1}u_X^{-}K^{-}_{ij}.
\endaligned$$
\end{pf}

We recall from \cite[\S5.7]{DJX} the following elements in $\cD C(\coh\bbX)$:
$$\aligned
\eta^+_{ij}
&=(-1)^{j+1}v^{p_i}K_{ij}[u^+_{i,p_i},\cdots,u^+_{i,j+1}, u^+_{i,1},\cdots,u^+_{i,j-1}]_{v^{-1}}\;\text{ and}\\
\eta^-_{ij}
&=(-1)^{p_i-j}v^{-1}[u^-_{i,j-1},\cdots, u^-_{i,1},u^-_{i,j+1},\cdots, u^-_{i,p_i}]_{v}K_{ij}^{-1},
\endaligned$$
 where $(i,j)\in \whI$. In particular,
$$\eta^+_{i1}=u_{S_{i0}^{(p_i-1)}}^{+} K_{[S_{i0}^{(p_i-1)}]}^{-1}\;\text{ and }\;
\eta^-_{i1}=-v^{-1}u_{S_{i0}^{(p_i-1)}}^{-} K_{[S_{i0}^{(p_i-1)}]}.$$
Then by \cite[Thm.~5.5]{DJX}, $\Xi(x^\pm_{\oz,\mp1})=\eta_{ij}^\mp$ for $\oz=\oz_{ij}\in I\backslash\{\star\}$.
Moreover, $[\eta^{\pm}_{ij}, u^{\mp}_{kl}]=0$ whenever
$\oz_{ij}\neq \oz_{kl}$ and $[\eta^{\pm}_{ij}, u^{\mp}_{\co(k\vc)}]=0$ for $j\geq 2$ and
$k\in\bbZ$.

\begin{Lem}\label{from eta to obtain a mod} Let $(i,j)\in\whI$ with $j\geq 2$. Then
 \begin{itemize}
    \item[(1)]
         $[u_{i1}^{-}, \cdots, u_{i,j-2}^{-},u_{i,j-1}^{-}, \eta_{ij}^{+}]_{v^{-1}}
         =(-v)^{-(j-1)}u_{S_{i0}^{(p_i-j)}}^{+} K_{[S_{i0}^{(p_i-j)}]}^{-1};$\\
    \item[(2)] $[\eta_{ij}^{-},u^{+}_{i,j-1}, u^{+}_{i,j-2},\cdots,u^{+}_{i,1}]_v
    =-v^{-1}u_{S_{i0}^{(p_i-j)}}^{-} K_{[S_{i0}^{(p_i-j)}]}$.
  \end{itemize}
Consequently,
  \begin{itemize}
    \item[(3)]
         $[u_{i,p_i-1}^{-}, \cdots, u_{i,j+1}^{-},u_{i1}^{-}, \cdots, u_{i,j-1}^{-}, \eta_{ij}^{+}]_{v^{-1}}
         =(-1)^{j+1}v^{-p_i+2}u_{i0}^{+} K_{i0}^{-1}K_{\delta};$\\
    \item[(4)] $[\eta_{ij}^{-},u^{+}_{i,j-1},\cdots,u^{+}_{i,1},u^{+}_{i,j+1},\cdots,u^{+}_{i,p_i-1}]_v
    =(-1)^{p_i-j}v^{-1}u_{i0}^{-} K_{i0}K_{\delta}^{-1}$.
  \end{itemize}

\end{Lem}

\begin{pf}
(1) By the definition,
$$\aligned
\eta^+_{ij}&
=(-1)^{j+1}v^{p_i}K_{ij}[u^+_{i, p_i},\cdots,u^+_{i,j+1}, u^+_{i,1},\cdots,u^+_{i,j-1}]_{v^{-1}}\\
{}&=(-1)^{j+1}v^{p_i-2}[u^+_{i, p_i},\cdots,u^+_{i,j+1}, u^+_{i,1},\cdots,u^+_{i,j-1}]_{v^{-1}}K_{ij}.
\endaligned$$
A calculation together with Lemma \ref{cancel ui} shows that
$$\aligned
{}&[u^-_{i,j-1},[u^+_{i, p_i},\cdots,u^+_{i,j+1}, u^+_{i,1},\cdots,u^+_{i,j-1}]_{v^{-1}}K_{ij}]_{v^{-1}}\\
=&[u^-_{i,j-1},[u^+_{i, p_i},\cdots,u^+_{i,j+1}, u^+_{i,1},\cdots,u^+_{i,j-1}]_{v^{-1}}] K_{ij}\\
=&v^{-2}[u^+_{i, p_i},\cdots,u^+_{i,j+1}, u^+_{i,1},\cdots,u^+_{i,j-2}]_{v^{-1}} K_{i,j-1}K_{ij}.
\endaligned$$
Thus, applying the operator $[u^-_{i,j-1},-]_{v^{-1}}$ to the expression
\begin{equation}\label{expression for trunct eta}
[u^+_{i, p_i},\cdots,u^+_{i,j+1}, u^+_{i,1},\cdots,u^+_{i,j-1}]_{v^{-1}}K_{ij}
\end{equation}
will cancel the term $u^+_{i,j-1}$ and produce a new factor $v^{-2}K_{i,j-1}$.
Keep on applying the operators $[u^-_{i,k},-]_{v^{-1}}$ for $k$ from $j-2$ to 1 step by step, we finally
obtain from Lemma \ref{drinfeld relations for A2}(3) that
$$\aligned
{}&[u_{i1}^{-}, \cdots, [u_{i,j-2}^{-},[u_{i,j-1}^{-}, \eta_{ij}^{+}]_{v^{-1}}]_{v^{-1}}\cdots]_{v^{-1}}\\
=&(-1)^{j+1}v^{p_i-2}(v^{-2})^{j-1}[u^+_{i, p_i},\cdots,u^+_{i,j+1}]_{v^{-1}} K_{i,1}\cdots K_{i,j-1}K_{ij}\\
=&(-1)^{j+1}v^{p_i-2j}v^{-(p_i-j-1)}u_{S_{i0}^{(p_i-j)}}^{+} K_{i,1}\cdots K_{i,j-1}K_{ij}\\
=&(-v)^{-(j-1)}u_{S_{i0}^{(p_i-j)}}^{+} K_{[S_{i0}^{(p_i-j)}]}^{-1}K_{\delta}.
\endaligned$$
 Moreover, if we further keep on applying the operators $[u^-_{i,k},-]_{v^{-1}}$ for
 $k$ from $j+1$ to $p_i-1$ step by step, we will cancel all the other terms in the
 expression (\ref{expression for trunct eta}) except $u_{i,p_i}^+$. Hence we get the statement (3).

Dually, we can prove the statements (2) and (4).
\end{pf}

\begin{Prop} \label{attach to each x a map}
For each $\vx\in\bbL$, there are uniquely determined maps $\sigma^r_{\vx}, \sigma^l_{\vx}: I\to\bbZ$
such that $\Xi(x^{+}_{\oz,\sigma^r_{\vx}(\oz)})\in \cD C(\co(\vx)\!^{\bot})$ and
$\Xi(x^{+}_{\oz,\sigma^l_{\vx}(\oz)})\in \cD C({}^{\bot}\!\co(\vx))$ for all $\oz\in I$.
\end{Prop}

\begin{pf}
We only show the existence of $\sigma^r_{\vx}$. The proof for the existence of $\sigma^l_{\vx}$ is similar.
It suffices to show that for each $\oz\in I$, there exists a unique integer $k$ such that
$\Xi(x^{+}_{\oz,k})\in \cD C(\co(\vx)\!^{\bot})$. (We then define $\sigma^r_{\vx}(\oz)=k$.)

Write $\vx=\sum\limits_{i=1}^{t}l_i\vx_i+l\vc$ in the normal form.
In case $\oz=\star$, we have $\Xi(x^{+}_{\star,k})=u_{\co(k\vc)}^+$. Then
$k=\sum_{i=1}^{t}l_i+l-1$ is the unique integer
satisfying $\co(k\vc)\in \co(\vx)\!^{\bot}$. Now we consider
the case $\oz=\oz_{ij}\in I\backslash\{\star\}$. Note that for each
torsion sheaf $S$ with $[S]\geq \delta$, we have $\Hom(\co(\vx), S)\neq 0$, i.e.,
$S\not\in\co(\vx)\!^{\bot}$. Thus, for every $k\neq 0,-1$,
$\Xi(x^{+}_{\oz,k})\not\in \cD C(\co(\vx)\!^{\bot})$. We need to show that exactly
one of $\Xi(x^{+}_{\oz,0})=u_{S_{ij}}^+$ and $\Xi(x^{+}_{\oz,-1})=\eta^-_{ij}$ lies in
$\cD C(\co(\vx)\!^{\bot})$. This is true because
$$S_{ij}\not\in\co(\vx)\!^{\bot} \Longleftrightarrow
j=l_i \Longleftrightarrow \eta^-_{ij}\in\cD C(\co(\vx)\!^{\bot}).$$
\end{pf}

\subsection{}
Combining Lemma \ref{a map from Gamma to Z induces an iso of quantum group} and
Proposition \ref{attach to each x a map} gives rise to a homomorphism
$\rho_{\vx}^r: \bfU_v(\fg_Q)\to \cD C(\co(\vx)\!^{\bot})$ taking
$$e_\oz\lmto \Xi(x^{+}_{\oz,\sigma^r_{\vx}(\oz)}),\;
f_{\oz}\lmto \Xi(x^{-}_{\oz,-\sigma^r_{\vx}(\oz)}),\;
K_\oz^{\pm1}\lmto K^{\pm1}_\oz K^{\pm1}_{\sigma^r_{\vx}(\oz)\delta},\;
\forall\,\oz\in I$$
 and a homomorphism $\rho_{\vx}^l: \bfU_v(\fg_Q)\to \cD C({}^{\bot}\!\co(\vx))$ taking
$$e_\oz\lmto \Xi(x^{+}_{\oz,\sigma^l_{\vx}(\oz)}),\;
f_{\oz}\lmto \Xi(x^{-}_{\oz,-\sigma^l_{\vx}(\oz)}),\;
K_\oz^{\pm1}\lmto K^{\pm1}_\oz K^{\pm1}_{\sigma^l_{\vx}(\oz)\delta},\;
\forall\,\oz\in I.$$

\begin{Prop} \label{iso-rho_{jx_i}} For each $(i,j)\in \widehat I$,
$\rho_{j\vx_i}^r$ and $\rho_{j\vx_i}^l$ are isomorphisms.
\end{Prop}

\begin{pf} We only prove that $\rho_{j\vx_i}^r$ is an isomorphism. The proof
for $\rho_{j\vx_i}^l$ is analogous.

For convenience, we use the notation $E_{kl}=E_{\oz}$ and $F_{kl}=F_{\oz}$
for $\oz=\oz_{kl}\in I$. Then $\rho_{j\vx_i}^r$ satisfies
$$\aligned
\rho_{j\vx_i}^r: & E_{\star}\lmto u_{\co(-\delta_{j0}\vc)}^+,\; F_{\star}\lmto -vu_{\co(-\delta_{j0}\vc)}^-,\;  K_{\star}^{\pm1}\lmto K^{\pm1}_{[\co(-\delta_{j0}\vc)]},\\
& E_{ij}\lmto \eta_{ij}^-,\; F_{ij}\lmto \eta_{ij}^+,\;  K_{ij}^{\pm1}\lmto K^{\pm1}_{ij}K_{\mp\delta}\; \text{\;for\;} j\neq 0,\\
\endaligned$$
 and takes $E_{kl}\mapsto u^{+}_{kl},\; F_{kl}\mapsto -vu^{-}_{kl},\; K_{kl}^{\pm1}\mapsto K^{\pm1}_{kl}$ for
the remaining $\oz=\oz_{kl}\in I$.

Recall from \S5.2 that for $j=0$, $\cD C(\co\!^{\bot})$ is generated by $u_{\co(-\vc)}^{\pm}$,
$u_{S_{kl}}^{\pm}$, and $K_\az$ for $(k,l)\in {\whI}$
and $\az\in K_0(\co\!^{\bot})$. This implies that $\rho_{0}^r$ is surjective.
For $j\neq 0$, $\cD C(\co(j\vx_i)\!^{\bot})$ is generated by $u_{\co}^{\pm}, u_{S_{i0}}^{\pm}$,
$u_{S_{kl}}^{\pm}$, and $K_\az$ for $(k,l)\in {\whI}\backslash\{(i,j)\}$
and $\az\in K_0(\co(j\vx_i)\!^{\bot})$. By applying Lemma \ref{cancel ui}(2) to
the element $\eta^{\pm}_{ij}$, we obtain that $u_{S_{i0}}^{\pm}$ can be generated by
$\eta^{\pm}_{ij}$ and $u_{S_{ij'}}^{\pm}$ for $1\leq j'\neq j\leq p_i-1$,
together with $K_\az$'s. In other words, $\cD C(\co(j\vx_i)\!^{\bot})$ is generated by
$u_{\co}^{\pm}, \eta^{\pm}_{ij}$, $u_{S_{kl}}^{\pm}$, and $K_\az$ for $(k,l)\in {\whI}\backslash\{(i,j)\}$
and $\az\in K_0(\co(j\vx_i)\!^{\bot})$. Hence, $\rho_{j\vx_i}^r$ is surjective.

By \cite[Prop.~14]{GL87}, the object
$$T_{ij}:=\bigoplus\limits_{0<\vy\leq \vc}\co(j\vx_i-\vy)$$
is a tilting object in $\co(j\vx_i)\!^{\bot}$
whose endomorphism algebra is isomorphic to $kQ$. By a theorem of Green \cite{Gr95},
 we have an isomorphism
 $$\Lambda_{ij}: \bfU_v(\fg_Q)\lra \cD C(\co(j\vx_i)\!^{\bot})$$
taking $$\aligned
&E_{\star}\lmto u_{\co(j\vx_i-\vc)}^+,\;F_{\star}\lmto -vu_{\co(j\vx_i-\vc)}^-,
\;K_{\star}^{\pm}\lmto K_{[\co(j\vx_i-\vc)]}^{\pm},\;\\
&E_{kl}\lmto u_{S_{kl}(j\vx_i)}^+,\;F_{kl}\lmto -vu_{S_{kl}(j\vx_i)}^-,
\;K_{kl}^{\pm}\lmto K_{[S_{kl}(j\vx_i)]}^{\pm} \text{\;for\;} (k,l)\in{\whI}.\\
\endaligned$$
Consider the composition
$$\Psi_{ij}=\Lambda_{ij}^{-1}\circ \rho_{j\vx_i}^r: \bfU_v(\fg_Q)\lra
\cD C(\co(j\vx_i)\!^{\bot}) \lra \bfU_v(\fg_Q)$$
which is surjective. It is easy to see that
$$\Psi_{ij}(E_{kl})=E_{k,l-j\delta_{ik}},\; \forall\,(k,l)\neq (i,j).$$
Furthermore,
$$\Psi_{ij}(E_{ij})=\Lambda_{ij}^{-1}(\eta_{ij}^-)=(-1)^{j+1}v^{-p_i}
[F_{i,p_i-1},\cdots, F_{i,p_i+1-j},F_{i,1},\cdots, F_{i,p_i-j}]_{v}K_{ij}^{-1}.$$
 Since $u_{\co}^+=v^{-(p_i-j)}[u^+_{i,p_i},\cdots, u^+_{i,j+1},u^+_{\co(\vx-\vc)}]_{v^{-1}}$,
it follows that
$$\Psi_{ij}(E_{\star})=\Lambda_{ij}^{-1}(u_{\co}^+)
=v^{-(p_i-j)}[ E_{i,p_i-j},\cdots,E_{i,1},E_{\star}]_{v^{-1}}.$$
Similarly,
$$\Psi_{ij}(F_{kl})=F_{k,l-j\delta_{ik}},\; \forall\,(k,l)\neq (i,j),$$
$$\Psi_{ij}(F_{ij})=\Lambda_{ij}^{-1}(\eta_{ij}^+)=(-1)^{j+1}v^{p_i}K_{ij}
[E_{i,p_i-j},\cdots, E_{i,1},E_{i,p_i-j+1},\cdots, E_{i,p_i-1}]_{v^{-1}},$$
$$\Psi_{ij}(F_{\star})=\Lambda_{ij}^{-1}(-vu_{\co}^-)
=(-1)^{p_i-j}v^2[ F_{i,p_i-j},\cdots,F_{i,1},F_{\star}]_{v^{-1}}.$$

In order to prove that $\Psi_{ij}$ is an isomorphism, we consider the quantum
enveloping algebra ${\bf U}_\bfv(\fg_Q)$ of $\fg_Q$ over $\bbC(\bfv)$ which is the field of
rational functions in an indeterminate $\bfv$. Let
$\bfPsi_{ij}:{\bf U}_\bfv(\fg_Q)\ra {\bf U}_\bfv(\fg_Q)$ be the map
given by the same formulas as $\Psi_{ij}$ via substituting $\bfv$ for $v$. Since
$\cD C(\co(j\vx_i)\!^{\bot})$ admits a generic form, it follows that $\bfPsi_{ij}$
is a surjective $\bbC(\bfv)$-algebra homomorphism.

Let $\scr A$ be the localization of $\bbC[\bfv,\bfv^{-1}]$ by the multiplicative
set
$$\{[m_1]!\cdots [m_s]!\mid s\geq 1,\\,m_1,\ldots,m_s\geq 0\}.$$
 Then $\scr A$ is again a principal ideal domain. By $\bfU_\scrA$ we denote the $\scrA$-subalgebra
of ${\bf U}_\bfv(\fg_Q)$ generated by
$$E_\oz^{(m)}=\frac{E_\oz^m}{[m]!},\,F_\oz^{(m)}=\frac{F_\oz^m}{[m]!},\,K_\oz^{\pm1}\;\,
\text{for $\oz\in I$}.$$
 Then $\bfU_\scrA$ admits a triangular decomposition
$$\bfU_\scrA=\bfU^-_\scrA\otimes \bfU^0_\scrA\otimes \bfU^+_\scrA,$$
 where $\bfU_\scrA^-$ (resp., $\bfU_\scrA^+$) is the $\scrA$-subalgebra generated by
the $F_\oz^{(m)}$ (resp., $E_\oz^{(m)}$) and $\bfU_\scrA^0=\bfU_\scrA\cap {\bf U}_\bfv(\fg_Q)^0
=\bfU_\scrA\cap\bbC(\bfv)[K_\oz^{\pm1}\mid \oz\in I]$. Moreover, $\bfU_\scrA^-$ and $\bfU_\scrA^+$
are free $\scrA$-modules with decompositions
$$\bfU_\scrA^-=\bigoplus_{\xi\in\bbN I}(\bfU_\scrA^-)_{-\xi}\;\text{ and }\;
\bfU_\scrA^+=\bigoplus_{\xi\in\bbN I}(\bfU_\scrA^+)_\xi,$$
 and $\bfU_\scrA^0$ is also a free $\scrA$-module; see \cite{L90a,L93}.

Let ${\cal U}(\fg_Q)$ be the universal enveloping algebra of $\fg_Q$ which by
definition is a $\bbC$-algebra generated by $e_\oz,f_\oz,h_\oz$ ($\oz\in I$)
with certain relations. For simplicity, we write
$$e_{kl}=e_{\oz_{kl}},\;f_{kl}=f_{\oz_{kl}}\;\text{ for $(k,l)\in\whI$.}$$
By viewing $\bbC$ as an $\scrA$-module with the action of $\bfv$
being $1$, we have by \cite{L93} that the map
$$\phi:\bfU_\scrA\otimes\bbC/\langle K_\oz-1\mid \oz\in I\rangle\lra \calU(\fg_Q)$$
 taking
$$E_\oz\lmto e_\oz,\;F_\oz\lmto f_\oz,\;E_\oz F_\oz-F_\oz E_\oz\lmto h_\oz\;(\oz\in I),$$
is a $\bbC$-algebra isomorphism.

It is clear that $\bfPsi_{ij}$ induces an $\scrA$-algebra homomorphism
$$\bfPsi_{ij}^\scrA: \bfU_\scrA\lra \bfU_\scrA$$
 which further gives rise to a $\bbC$-algebra homomorphism
$$\Psi_{ij}^{(1)}: \bfU_\scrA\otimes\bbC/\langle K_\oz-1\mid \oz\in I\rangle\lra
\bfU_\scrA\otimes\bbC/\langle K_\oz-1\mid \oz\in I\rangle.$$
 Finally, we obtain a $\bbC$-algebra homomorphism
$$\phi\Psi_{ij}^{(1)}\phi^{-1}: \calU(\fg_Q)\lra \calU(\fg_Q)$$
 which induces a Lie algebra homomorphism $\theta_{ij}:\fg_Q\ra\fg_Q$ taking
$$\aligned
{}& \theta_{ij}(e_{kl})=e_{k,l-j\delta_{ik}},\;\theta_{ij}(f_{kl})=f_{k,l-j\delta_{ik}},\;
\forall\,(k,l)\neq (i,j),\\
&\theta_{ij}(e_{\star})=[e_{i,p_i-j},\cdots,e_{i,1},e_{\star}],\;
\theta_{ij}(f_{\star})=(-1)^{p_i-j}[f_{i,p_i-j},\cdots,f_{i,1},f_{\star}],\\
\theta_{ij}(e_{ij})&=(-1)^{j+1} [f_{i,p_i-1},\cdots, f_{i,p_i+1-j},f_{i,1},\cdots, f_{i,p_i-j}]\\
&=(-1)^{j+1} [[f_{i,p_i-1},\cdots, f_{i,p_i+1-j}],(-1)^{p_i-j-1}[f_{i,p_i-j},\cdots,f_{i,1}]]\\
&=(-1)^{p_i} [f_{i,p_i-1},\cdots,f_{i,1}],\\
\theta_{ij}(f_{ij})&=(-1)^{j+1}[e_{i,p_i-j},\cdots, e_{i,1},e_{i,p_i-j+1},\cdots, e_{i,p_i-1}]\\
&=(-1)^{j+1}[[e_{i,p_i-j},\cdots, e_{i,1}],(-1)^{j-2}[e_{i,p_i-1},\cdots, e_{i,p_i-j+1}]]\\
&=[e_{i,p_i-1},\cdots, e_{i,1}].
\endaligned$$
 It is straightforward to check that $\theta_{ij}$ coincides with $\widetilde\Omega_{ij}$
 defined in Proposition \ref{iso of Kac-Moody}.
Hence, $\theta_{ij}$ is an isomorphism. Since
$\widetilde\Omega_{ij}((\fg_Q)_\alpha)=(\fg_Q)_{\varpi(\az)}$ for all
$\az\in\Delta(Q)$, it follows from Poincar$\acute{\rm e}$--Birkhoff--Witt theorem
that $\phi\Psi_{ij}^{(1)}\phi^{-1}$ is an isomorphism, too.

Let $s\geq 1$ and take $\xi_1,\ldots,\xi_s,\eta_1,\ldots,\eta_s\in\bbN I$. Then the
$\scrA$-module
$$L_{\xi_1,\ldots,\xi_s,\eta_1,\ldots,\eta_s}:=\bigoplus_{i=1}^s
(\bfU_\scrA^-)_{\xi_i}(\bfU_\scrA^+)_{\eta_i}$$
 is free of finite rank. Since $\scrA$ is a principal ideal domain, the image
$$\Psi_{ij}^{(1)}(L_{\xi_1,\ldots,\xi_s,\eta_1,\ldots,\eta_s})\subset \bfU_\scrA$$
 is also a free $\scrA$-module. Then $\Psi_{ij}^{(1)}$ restricts to an $\scrA$-module
 homomorphism
$$\zeta: L_{\xi_1,\ldots,\xi_s,\eta_1,\ldots,\eta_s}\lra \Psi_{ij}^{(1)}(L_{\xi_1,\ldots,\xi_s,\eta_1,\ldots,\eta_s})$$
 which splits. Thus, $\Ker(\zeta)$ is a free $\scrA$-submodule of $L_{\xi_1,\ldots,\xi_s,\eta_1,\ldots,\eta_s}$.
Furthermore, $\zeta$ induces a surjective $\bbC$-linear map
$$\bigoplus_{i=1}^s (\calU(\fg_Q))_{\xi_i}(\calU(\fg_Q)_{\eta_i}\lra
\phi\Psi_{ij}^{(1)}\phi^{-1}\big(\bigoplus_{i=1}^s (\calU(\fg_Q))_{\xi_i}(\calU(\fg_Q)_{\eta_i}\big)$$
 which coincides with the restriction of $\phi\Psi_{ij}^{(1)}\phi^{-1}$ and whose
 kernel is isomorphic to $\Ker(\zeta)\otimes_\scrA\bbC$. Since $\phi\Psi_{ij}^{(1)}\phi^{-1}$
 is an isomorphism, it follows that
$$\Ker(\zeta)\otimes_\scrA\bbC=0,$$
 and hence, $\Ker(\zeta)=0$, i.~e., $\zeta$ is injective. By choosing
$\xi_1,\ldots,\xi_s,\eta_1,\ldots,\eta_s$ arbitrarily, we obtain that
the restriction of $\bfPsi_{ij}^\scrA$ to $\bfU_\scrA^-\bfU_\scrA^+$ is injective.
 Therefore, the restriction of $\bfPsi_{ij}$ to $\bfU_\bfv(\fg_Q)^-\bfU_\bfv(\fg_Q)^+$
 gives an isomorphism
$$ \bfU_\bfv(\fg_Q)^-\bfU_\bfv(\fg_Q)^+\lra \bfPsi_{ij}\big(\bfU_\bfv(\fg_Q)^-\bfU_\bfv(\fg_Q)^+\big).$$
By the definition, $\bfPsi_{ij}$ induces an isomorphism $\bfU_\bfv(\fg_Q)^0\ra \bfU_\bfv(\fg_Q)^0$.
Since
$$\bfU_\bfv(\fg_Q)=\bfPsi_{ij}\big(\bfU_\bfv(\fg_Q)^-\bfU_\bfv(\fg_Q)^+\big)\otimes_{\bbC(\bfv)} \bfU_\bfv(\fg_Q)^0,$$
we conclude that $\bfPsi_{ij}$ is an isomorphism.

Specializing $\bfv$ to $v$ shows that $\Psi_{ij}$ is an isomorphism.
This implies that $\rho_{j\vx_i}^r$ is an isomorphism, as desired.
\end{pf}

Similarly, it can been shown that both of $\rho_{\vx}^r$ and $\rho_{\vx}^l$ are
isomorphisms for each $\vx\in\bbL$.

\subsection{} Define
$$\Theta_{\vx}:=(\rho_{\vx}^l)^{-1} \wtilde{R}_{\vx} \rho_{\vx}^r:\bfU_v(\fg_Q)\lra \bfU_v(\fg_Q).$$
In other words, we have
the following commutative diagram:
$$\xymatrix{
  \bfU_v(\fg_Q) \ar[r]^{\Theta_{\vx}}\ar[d]^{\rho_{\vx}^r} &\bfU_v(\fg_Q) \ar[d]^{\rho_{\vx}^l}\\
   \cD C(\co(\vx)\!^{\bot}) \ar[r]^{\wtilde{R}_{\vx}} &\cD C({}^{\bot}\!\co(\vx))  }
  $$

The following results describe the automorphism $\Theta_{\vx}$ for $0\leq \vx<\vc$.

\begin{Prop} \label{thm for operator Rij for star} For each $\oz\in I$, we have
  \begin{itemize}
    \item[(1)]
         $\Theta_{0}(E_\oz)=
         \left\{\begin{array}{lll} -K_{\star}F_{\star}, && \oz=\star;\\
             { [E_{i1}, E_{\star}]_{v}}, && \oz=\oz_{i1};\\
              E_\oz, && \text{otherwise};
             \end{array} \right.$
    \item[(2)] $\Theta_{0}(F_\oz)=
         \left\{\begin{array}{lll} -F_{\star}K_{\star}^{-1}, && \oz=\star;\\
             { [F_{\star}, F_{i1}]_{v^{-1}}}, && \oz=\oz_{i1};\\
              F_\oz, && \text{otherwise};
             \end{array} \right.$
  \end{itemize}
\end{Prop}

\begin{pf} We only prove the statement (1). The second one can be proved similarly.

It is clear that $\Theta_{0}$ preserves all the Chevalley generators $e_\oz$ for
$\oz\not=\star, \oz_{i1}$. Note that
$$\wtilde{R}_0(u^+_{\co(-\vc)})=v^{-1}u^-_{\co(\vc)}K_{[\co(\vc)]}=v K_{[\co(\vc)]}u^-_{\co(\vc)}.$$
Thus
$$\Theta_{0}(E_{\star})= (\rho_{\vx}^l)^{-1} \wtilde{R}_0(u^+_{\co(-\vc)})=(\rho_{\vx}^l)^{-1}(v K_{[\co(\vc)]}u^-_{\co(\vc)})=-K_{\star}F_{\star}.$$
By Lemma \ref{drinfeld relations for A2},
$$\aligned \wtilde{R}_0(u^+_{S_{i1}})&=u^+_{\co(\vx_i)}
=[u^+_{\co(\vc)}, u^-_{S_{i0}^{(p_i-1)}}K_{[S_{i0}^{(p_i-1)}]}]_{v^{-1}}\\
=&-v^{-1}[u^-_{S_{i0}^{(p_i-1)}}K_{[S_{i0}^{(p_i-1)}]}, u^+_{\co(\vc)}]_{v}=[\eta_{i1}^-, u^+_{\co(\vc)}]_{v}.
\endaligned$$
 Therefore,
$$\Theta_{0}(E_{i1})=(\rho_{\vx}^l)^{-1}\wtilde{R}_0(u^+_{S_{i1}})
=(\rho_{\vx}^l)^{-1}([\eta_{i1}^-, u^+_{\co(\vc)}]_{v})=[E_{i1}, E_{\star}]_{v}. $$
\end{pf}

\begin{Prop} \label{thm for operator Rij for ij}
Fix $(i,j)\in\whI$. Then for each $\oz\in I$,
\begin{itemize}
    \item[(1)]
         $\Theta_{j\vx_i}(E_\oz)=
         \left\{\begin{array}{lll} (-1)^j K_{i1}\cdots K_{ij} [F_{ij}, \cdots, F_{i1}]_{v^{-1}}, && \oz=\star;\\
             K_{\star}K_{i1}\cdots K_{i,j-1} [F_{\star}, F_{i1}, \cdots, F_{i,j-1}]_{v^{-1}}, && \oz=\oz_{ij};\\
             { [E_{i,j+1}, \cdots, E_{i1}, E_{\star}]_{v}}, && \oz=\oz_{i,j+1};\\
            (-1)^j  [E_{k1}, E_{\star}, E_{i1}, \cdots, E_{ij}]_{v}, && \oz=\oz_{k1}\neq \oz_{i1};\\
              E_\oz, && \text{otherwise};
             \end{array} \right.$
    \item[(2)] $\Theta_{j\vx_i}(F_\oz)=
         \left\{\begin{array}{lll} (-1)^j[E_{i1}, \cdots, E_{ij}]_{v}K^{-1}_{i1}\cdots K^{-1}_{ij}, && \oz=\star;\\
             { [E_{i,j-1}, \cdots, E_{i1}, E_{\star}]_{v}K^{-1}_{\star}K^{-1}_{i1}\cdots K^{-1}_{i,j-1}}, && \oz=\oz_{ij};\\
             { [F_{\star}, F_{i1}, \cdots, F_{i,j+1}]_{v^{-1}}}, && \oz=\oz_{i,j+1};\\
             (-1)^j [F_{ij}, \cdots, F_{i1}, F_{\star}, F_{k1}]_{v}, && \oz=\oz_{k1}\neq \oz_{i1};\\
              F_\oz, && \text{otherwise}.
             \end{array} \right.$
  \end{itemize}
\end{Prop}

\begin{pf}
We only prove  (1). The proof for (2) is entirely analogous.

Since $$\wtilde{R}_{j\vx_i}(u^+_{\co})=v^{-1}u^-_{S_{ij}^{(j)}} K_{[S_{ij}^{(j)}]}
=v K_{[S_{ij}^{(j)}]}u^-_{S_{ij}^{(j)}}
=v^{j} K_{i1}\cdots K_{ij} [u_{ij}^-, \cdots, u_{i1}^-]_{v^{-1}},$$
 it follows that $\Theta_{j\vx_i}(E_{\star})=(-1)^j K_{i1}\cdots K_{ij} [F_{ij}, \cdots, F_{i1}]_{v^{-1}}.$
Further, by the equalities
$$\aligned
\wtilde{R}_{j\vx_i}(\eta^-_{ij})&
=\wtilde{R}_{j\vx_i}((-1)^{p_i-j}v^{-1}[u^-_{i,j-1},\cdots, u^-_{i,1},u^-_{i,j+1},\cdots, u^-_{i, p_i}]_{v}K_{ij}^{-1})\\
&=(-1)^{p_i-j}v^{-1}[u^-_{i,j-1},\cdots, u^-_{i,1},u^-_{\co((j+1)\vx_i)},\cdots, u^-_{i, p_i}]_{v}K_{[\co((j-1)\vx_i)]}\\
&=(-1)^{p_i-j}v^{-1}[u^-_{i,j-1},\cdots, u^-_{i,1},(-1)^{p_i-j-1}u^-_{\co(\vc)}]_{v}K_{[\co((j-1)\vx_i)]},
\endaligned$$
we obtain that
$$\aligned
\Theta_{j\vx_i}(E_{ij})&=(-v)^{-j-1}[F_{i,j-1},\cdots, F_{i,1},F_{\star}]_{v} K_{\star}K_{i1}\cdots K_{i,j-1}\\
&=K_{\star}K_{i1}\cdots K_{i,j-1} [F_{\star}, F_{i1}, \cdots, F_{i,j-1}]_{v^{-1}}.\endaligned$$
Finally, the equalities
$$\aligned
\wtilde{R}_{j\vx_i}(u^+_{i,j+1})
=&u^+_{\co((j+1)\vx_i)}
=[u^+_{\co(\vc)}, u_{S_{i0}^{(p_i-j-1)}}^{-} K_{[S_{i0}^{(p_i-j-1)}]}]_{v^{-1}}\\
=&-v^{-1}[u_{S_{i0}^{(p_i-j-1)}}^{-} K_{[S_{i0}^{(p_i-j-1)}]}, u^+_{\co(\vc)}]_{v}
=[[\eta^-_{i,j+1}, u^+_{ij}, \cdots, u^+_{i1}]_{v}, u^+_{\co(\vc)}]_{v}\\
=&[\eta^-_{i,j+1}, u^+_{ij}, \cdots, u^+_{i1}, u^+_{\co(\vc)}]_{v}
\endaligned$$
imply that $\Theta_{j\vx_i}(E_{i,j+1})=[E_{i,j+1}, \cdots, E_{i1}, E_{\star}]_{v}$,
and the equalities
$$\aligned
\wtilde{R}_{j\vx_i}(u^+_{k,1})&=u^+_{\co(j\vx_i+\vx_k)}
=[u^+_{\co(j\vx_i+\vc)}, u_{S_{k0}^{(p_k-1)}}^{-} K_{[S_{k0}^{(p_k-1)}]}]_{v^{-1}}\\
&=-v^{-1}[u_{S_{k0}^{(p_k-1)}}^{-} K_{[S_{k0}^{(p_k-1)}]}, u^+_{\co(j\vx_i+\vc)}]_{v}\\
&=[\eta^-_{k,1}, (-1)^j [u^+_{\co(\vc)}, u^+_{i1}, \cdots, u^+_{ij}]_{v}]_{v}\\
&=(-1)^j [\eta^-_{k,1}, u^+_{\co(\vc)}, u^+_{i1}, \cdots, u^+_{ij}]_{v}
\endaligned$$
give that $\Theta_{j\vx_i}(E_{k,1})=(-1)^j  [E_{k1}, E_{\star}, E_{i1}, \cdots, E_{ij}]_{v}$.
It is easy to see that the remaining generators are fixed by $\Theta_{j\vx_i}$. This finishes
the proof.
\end{pf}

\subsection{}
Recall from \cite[Ch.~37]{L93} that for each fixed $\oz\in I$, there is an automorphism
$T_\oz=T'_{\oz,1}$ of $\bfU_v(\fg_Q)$, called the Lusztig's symmetry, which is defined by
$$T_{\oz}(E_u)= \left\{\begin{array}{lll} -K_{\oz}F_{\oz}, && u=\oz;\\
             { [E_{u},E_{\oz}]_{v}}, && a_{\oz u}=-1;\\
              E_u, && \text{otherwise,}
             \end{array} \right.
T_\oz(F_u)=\left\{\begin{array}{lll} -E_\oz K_\oz^{-1}, && u=\oz;\\
             { [F_\oz,F_u]_{v^{-1}}}, && a_{\oz u}=-1;\\
              F_u, && \text{otherwise,}
             \end{array} \right.$$
$$\text{and }\;
T_\oz(K_u)=\left\{\begin{array}{lll} K_\oz^{-1}, && u=\oz;\\
             K_\oz K_u, && a_{\oz u}=-1;\\
              K_u, && \text{otherwise,}
             \end{array} \right.$$
where $u\in I$. For $\oz=\oz_{ij}$ with $(i,j)\in\whI$, we simply write $T_{ij}=T_\oz$.

\begin{Prop} \label{thm for operator psi v} Fix $(i,j)\in \whI$ and set
${\cal J}_{ij}=T_{\star}T_{i1}\cdots T_{ij}\cdots T_{i1}T_{\star}$. Then for $\oz\in I$,
  \begin{itemize}
    \item[(1)]
         ${\cal J}_{ij}(E_\oz)=
         \left\{\begin{array}{lll} -K_{i1}\cdots K_{ij} [F_{ij}, \cdots, F_{i1}]_{v^{-1}}, && \oz=\star;\\
             -K_{\star}K_{i1}\cdots K_{i,j-1} [F_{\star}, F_{i1}, \cdots, F_{i,j-1}]_{v^{-1}}, && \oz=\oz_{ij};\\
             { [E_{i,j+1}, \cdots, E_{i1}, E_{\star}]_{v}}, && \oz=\oz_{i,j+1};\\
             { [E_{k1}, E_{\star}, E_{i1}, \cdots, E_{ij}]_{v}}, && \oz=\oz_{k1}\neq \oz_{i1};\\
              E_\oz, && \text{otherwise};
             \end{array} \right.$
    \item[(2)] ${\cal J}_{ij}(F_\oz)=
       \left\{\begin{array}{lll} -[E_{i1}, \cdots, E_{ij}]_{v}K^{-1}_{i1}\cdots K^{-1}_{ij}, && \oz=\star;\\
       -[E_{i,j-1}, \cdots, E_{i1}, E_{\star}]_{v}K^{-1}_{\star}K^{-1}_{i1}\cdots K^{-1}_{i,j-1}, && \oz=\oz_{ij};\\
       {[F_{\star}, F_{i1}, \cdots, F_{i,j+1}]_{v^{-1}}}, && \oz=\oz_{i,j+1};\\ %why should add { before [? can be begin with [?
             { [F_{ij}, \cdots, F_{i1}, F_{\star}, F_{k1}]_{v}}, && \oz=\oz_{k1}\neq \oz_{i1};\\
              F_\oz, && \text{otherwise};
             \end{array} \right.$
  \end{itemize}
\end{Prop}

\begin{pf} We only prove (1). The statement (2) can be proved similarly.

By Lemma \ref{psi for linear case}, we only need to prove the formulas for $\oz=\oz_{k1}\neq \oz_{i1}$.
In fact,
$$\aligned{\cal J}_{ij}(E_{k1})
&=T_{\star}T_{i1}\cdots T_{ij}\cdots T_{i1}([E_{k1},E_{\star}]_v)\\
&=[T_{\star}T_{i1}\cdots T_{ij}\cdots T_{i1}(E_{k1}), T_{\star}T_{i1}\cdots T_{ij}\cdots T_{i1}\cdot T_{\star}(-F_{\star}K^{-1}_{\star})]_v\\
&=[T_{\star}(E_{k1}), -{\cal J}_{ij}(F_{\star}){\cal J}_{ij}(K^{-1}_{\star})]_v\\
&=[[E_{k1},E_{\star}]_{v}, [E_{i1}, \cdots, E_{ij}]_{v}]_v\\
&=[E_{k1},E_{\star}, E_{i1}, \cdots, E_{ij}]_{v}.\endaligned
$$
\end{pf}

 For each $\oz=\oz_{ij}\in I$, we denote by $\epsilon_{ij}=\epsilon_{\oz}: \bfU_v(\fg_Q)\to \bfU_v(\fg_Q)$ the
involution taking
$$E_\oz\lmto -E_\oz,\; F_\oz\lmto -F_\oz;$$
and preserving all the other generators. For $(i,j)\in\whI$, set $$\kappa_{ij}=\epsilon_{ij}\epsilon_{\star}^{j-1}\prod_{k\neq i}\epsilon_{k1}^{j}.$$

It is easy to check that for each $\oz_{ij}\in I$, both $\Theta_{j\vx_i}$ and ${\cal J}_{ij}$
are compatible with $\Ref_{j\vx_i}$ defined in Section \ref{section for mutation on Grothendieck group},
i.e., for each $\alpha\in \bbZ I$,
$$\Theta_{j\vx_i}(K_{\alpha}^{\pm 1})=K_{\Ref_{j\vx_i}(\alpha)}^{\pm 1}={\cal J}_{ij}(K_{\alpha}^{\pm 1}).$$
Combining Propositions \ref{thm for operator Rij for star}, \ref{thm for operator Rij for ij}
and \ref{thm for operator psi v} together gives the main result of this section
which relates Lusztig's symmetries $T_\oz$ of $\bfU_v(\fg)$ with the $\Theta_{j\vx_i}$.

\begin{Thm}
We have $\Theta_{0}=T_{\star}$, $\Theta_{j\vx_i}=\kappa_{ij} {\cal J}_{ij}$
and
$$T_{ij}=\kappa_{ij} \Theta_{0}^{(-1)^j}\cdots\Theta_{(j-1)\vx_i}^{-1}
\Theta_{j\vx_i}\Theta_{(j-1)\vx_i}^{-1}\cdots \Theta_{0}^{(-1)^j},$$
 where $(i,j)\in \whI$.
\end{Thm}

\bigskip

 \appendix
  \renewcommand{\appendixname}{Appendix~\Alph{section}}

\section{The linear quiver}

\medskip

In this appendix we collect some results concerning the linear quiver
\begin{center}
\begin{pspicture}(-1,0.3)(8,1)
\psset{xunit=1cm,yunit=.8cm}
\uput[d](0,1){$A_n:$}
\psdot*[dotsize=3pt](1,0.5) \psdot*[dotsize=3pt](2,0.5)
\psdot*[dotsize=3pt](3,0.5)
\psdot*[dotsize=3pt](4.5,0.5) \psdot*[dotsize=3pt](5.5,0.5)
\psline{<-}(1,0.5)(2,0.5) \psline{<-}(2,0.5)(3,0.5)
\psline[linestyle=dotted,linewidth=1pt](3,0.5)(4.5,0.5)
\psline{<-}(4.5,0.5)(5.5,0.5)
\uput[d](1,1.2){$_1$} \uput[d](2,1.2){$_2$} \uput[d](3,1.2){$_3$}
\uput[d](4.5,1.2){$_{n-1}$} \uput[d](5.5,1.2){$_n$}
\end{pspicture}
\end{center}

\medskip

\noindent These results have been used in Sections 3 and 5. To this quiver
we have the associated Weyl group $W={\frak S}_{n+1}$, root system $\Delta(A_n)$,
Kac--Moody Lie algebra $\fg=\frak{sl}_{n+1}$ and
quantum enveloping algebra ${\bf U}={\bf U}_v(\frak{sl}_{n+1})$ which is generated by
$E_i,F_i, K_i^{\pm1}$ for $1\leq i\leq n$.

For each vertex $1\leq i\leq n$, we denote by $S_i$ the simple $\bbF A_n$-module corresponding to
$i$ and by $P_i$ the projective cover of $S_i$. Let $\refl_{[P_i]}$
be the reflection in $W={\frak S}_{n+1}$ associated with the dimension vector of $P_i$
which is a positive root in $\Delta(A_n)$ by \cite{Gab}.

\begin{Lem}\label{A-type} For each $2\leq i\leq n$, we have the
equalities in $W={\frak S}_{n+1}$:
  \begin{itemize}
    \item[(1)] $\refl_{[P_i]}=\refl_{i}\refl_{i-1}\cdots \refl_{1}\cdots \refl_{i-1}\refl_{i}=\refl_{1}\refl_{2}\cdots \refl_{i}\cdots \refl_{2}\refl_{1}$;
    \item[(2)]  $\refl_i=\refl_{[P_{i-1}]}\refl_{[P_{i}]}\refl_{[P_{i-1}]}=\refl_{[P_{i}]}\refl_{[P_{i-1}]}\refl_{[P_{i}]}$.
  \end{itemize}
\end{Lem}

\begin{pf} First of all, we have the following known basic facts:
  \begin{itemize}
    \item[(a)] $\refl_i \refl_{\alpha}=\refl_{\refl_{i}(\alpha)}\refl_i$ for each root $\alpha\in\Delta(A_n)$;
    \item[(b)]  $\refl_i^2=1$ and $\refl_i\refl_j=\refl_j\refl_i$ if $|i-j|>1$;
    \item[(c)] $\refl_i\refl_{i+1}\refl_i=\refl_{i+1}\refl_i\refl_{i+1}$ for $1\leq i\leq n-1$.
  \end{itemize}
Since $\refl_i([P_i])=[P_i]-([P_i], \az_i)\az_i=[P_{i-1}]$ for $i\geq 2$,
it follows from (a) that
$$\refl_i \refl_{[P_i]}=\refl_{\refl_{i}([P_i])}\refl_i=\refl_{[P_{i-1}]}\refl_i.$$
 Hence, $\refl_{[P_i]}=\refl_i\refl_{[P_{i-1}]}\refl_i$. An inductive argument gives
$$\refl_{[P_i]}=\refl_{i}\refl_{i-1}\cdots \refl_{1}\cdots \refl_{i-1}\refl_{i}.$$
By applying (b) and (c), we obtain that
$$\aligned \refl_{[P_i]}&=\refl_{i}\refl_{i-1}\cdots \refl_3(\refl_{1}\refl_2\refl_1)\refl_3\cdots \refl_{i-1}\refl_{i}=\refl_{1}\refl_{i}\refl_{i-1}\cdots \refl_3\refl_2\refl_3\cdots \refl_{i-1}\refl_{i}\refl_{1}\\
&=\cdots=\refl_{1}\refl_{2}\cdots \refl_{i}\cdots \refl_{2}\refl_{1}.
\endaligned$$
 For the statement (2), we observe from (1) that
$$\begin{array}{rl} \refl_{[P_{i}]}\refl_{[P_{i-1}]}
=&(\refl_{1}\refl_{2}\cdots \refl_{i}\cdots \refl_{2}\refl_{1})(\refl_{1}\refl_{2}\cdots \refl_{i-1}\cdots \refl_{2}\refl_{1})\\
=&\refl_{1}\refl_{2}\cdots \refl_{i}\refl_{i-2}\cdots \refl_{2}\refl_{1}\\
=&\refl_{1}\refl_{2}\cdots \refl_{i-2}\refl_{i-1}\refl_{i-2}\cdots \refl_{2}\refl_{1} \refl_{i}\\
=&\refl_{[P_{i-1}]}\refl_{i}.
\end{array}$$
Consequently, $\refl_i=\refl_{[P_{i-1}]}\refl_{[P_{i}]}\refl_{[P_{i-1}]}$. On the other hand,
$$\begin{array}{ll} &\refl_{[P_{i}]}\refl_{[P_{i-1}]}\refl_{[P_{i}]}\\
=&(\refl_{1}\refl_{2}\cdots \refl_{i}\cdots \refl_{2}\refl_{1})(\refl_{1}\refl_{2}\cdots \refl_{i-1}\cdots \refl_{2}\refl_{1})(\refl_{1}\refl_{2}\cdots \refl_{i}\cdots \refl_{2}\refl_{1})\\
=&\refl_{1}\refl_{2}\cdots \refl_{i-1}(\refl_{i}\refl_{i-1}\refl_{i})\refl_{i-1}\cdots \refl_{2}\refl_{1}\\
=&\refl_{1}\refl_{2}\cdots \refl_{i-2}\refl_{i}\refl_{i-2}\cdots \refl_{2}\refl_{1}=\refl_i.
\end{array}$$
\end{pf}

Let $\cD H(A_n)$ denote the double Ringel--Hall algebra of the quiver $A_n$. It is known
from \cite{R90,R90a} that $\cD H(A_n)$ is generated by $u^\pm_i=u^\pm_{S_i},\,K^{\pm1}_i =K^{\pm1}_{[S_i]}$
for $1\leq i\leq n$. Moreover, $\cD H(A_n)$ can be identified with ${\bf U}={\bf U}_v(\frak{sl}_{n+1})$
via
$$u_i^+\lmto E_i,\;u_i^-\lmto -v^{-1} F_i,\;K_i^{\pm1}\lmto K_i^{\pm1},\;\forall\,1\leq i\leq n.$$
 A direct calculation gives the following result.

\begin{Lem}\label{drinfeld relations for A2}
{\rm (1)} $[u_1^{\pm}, u_2^{\pm}]_{v}=-u_{P_2}^{\pm}$ and
$[u_2^{\pm}, u_1^{\pm}]_{v^{-1}}=v^{-1}u_{P_2}^{\pm};$\\

{\rm(2)} $[u_{P_2}^{\pm}, u_2^{\mp}K_2^{\pm1}]_{v^{-1}}=u_1^{\pm}$ and
               $[K_1^{\pm1}u_1^{\mp}, u_{P_2}^{\pm}]_{v^{-1}}=u_2^{\pm}$;

{\rm (3)} $[u_1^{\pm}, u_2^{\pm},\cdots, u_n^{\pm}]_{v}=(-1)^{n-1}u_{P_n}^{\pm}\;
\text{ and }\;[u_n^{\pm},\cdots, u_2^{\pm}, u_1^{\pm}]_{v^{-1}}=v^{-(n-1)}u_{P_n}^{\pm}.$
\end{Lem}

As in \S5.6, for each $1\leq i\leq n$, we have the Lusztig's symmetry $T_i$ of
$\bfU_v(\frak{sl}_{n+1})$. In the following we list some basic properties
of the $T_i$. For completeness, we provide some proofs.

\begin{Lem}\label{known property for Ti} The following equalities hold.
\begin{itemize}
    \item[(1)]
         $T_i([E_{i},E_{j}]_{v})=E_j$ and  $T_i([F_{j},F_{i}]_{v^{-1}})=F_j$, where $i=j\pm 1$;
     \item[(2)]
         $T_iT_j(E_{i})=E_j$ and  $T_iT_j(F_{i})=F_j$, where $i=j\pm1$;
     \item[(3)]
         $T_i([E_{i\pm 1},E_{i},E_{i\mp 1}]_{v})=[E_{i\pm 1},E_{i},E_{i\mp 1}]_{v}$;
     \item[(4)] $T_i([F_{i\pm 1},F_{i},F_{i\mp 1}]_{v^{-1}})=[F_{i\pm 1},F_{i},F_{i\mp 1}]_{v^{-1}}$.

  \end{itemize}
\end{Lem}

\begin{pf} These equalities can be easily obtained from the following equalities
$$[-K_iF_i, [E_j,E_i]_v]_v=E_j\;\text{ and }\; [[F_i,F_j]_{v^{-1}}, -E_iK_{i}^{-1}]_{v^{-1}}=F_j,$$
where $a_{ij}=-1$.
\end{pf}

\begin{Lem}\label{property for T1...Tn} We have the following equalities.
\begin{itemize}
    \item[(1)]
        $T_1\cdots T_j(E_i)=
         \left\{\begin{array}{lll} [E_{j+1}, \cdots, E_{1}]_{v}, && i=j+1;\\
             -K_1\cdots K_j[F_1, \cdots, F_j]_{v^{-1}}, && i=j;\\
             E_{i+1}, && i<j;\\
              E_i, && \text{otherwise},
             \end{array} \right.$ \\
     \item[(2)]
             $T_1\cdots T_j(F_i)=
         \left\{\begin{array}{lll} [F_{1}, \cdots, F_{j+1}]_{v^{-1}}, && i=j+1;\\
             -[E_j, \cdots, E_1]_{v}K^{-1}_1\cdots K^{-1}_j, && i=j;\\
             F_{i+1}, && i<j;\\
              F_i, && \text{otherwise}.
             \end{array} \right.$\\
  \end{itemize}
\end{Lem}

\begin{pf} We use induction on $j$ for the proof. For $j=1$, there is nothing to prove.
Now assume statements (1) and (2) hold for $j-1$ and we want to show that the equalities in (1)
hold for $j$. It is trivial that $T_1\cdots T_j(E_i)=E_i$ for $i\geq j+2$.  For the case
$i=j+1$, we have
$$T_1\cdots T_j(E_{j+1})=T_1\cdots T_{j-1}([E_{j+1},E_j]_v)=[E_{j+1}, \cdots, E_{1}]_{v};$$
for the case $i=j$,
$$\aligned
{}&T_1\cdots T_j(E_{j})=T_1\cdots T_{j-1}(-K_jF_j)=-T_1\cdots T_{j-1}(K_j)T_1\cdots T_{j-1}(F_j)\\
=&-K_1\cdots K_j[F_1, \cdots, F_j]_{v^{-1}};
\endaligned$$
and for the case $i<j$,
$$T_1\cdots T_j(E_i)=T_1\cdots T_iT_{i+1}(E_i)=T_1\cdots T_{i-1}(E_{i+1})=E_{i+1}.$$
Similarly, the equalities in (2) hold for $j$.
\end{pf}

\begin{Lem}\label{property for Tn...T1}
{\rm (1)} $T_j\cdots T_1(E_i)=
         \left\{\begin{array}{lll} [E_{j+1}, E_{j}]_{v}, && i=j+1;\\
              E_{i-1}, && 2\leq i\leq j;\\
              -K_1\cdots K_j[F_j, \cdots, F_1]_{v^{-1}}, && i=1;\\
              E_i, && \text{otherwise};
             \end{array} \right.$

{\rm (2)}  $T_j\cdots T_1(F_i)=
         \left\{\begin{array}{lll} [F_{j}, F_{j+1}]_{v^{-1}}, && i=j+1;\\
             F_{i-1}, && 2\leq i\leq j;\\
              -[E_1, \cdots, E_j]_{v}K^{-1}_1\cdots K^{-1}_j, && i=1;\\
                F_i, && \text{otherwise}.
             \end{array} \right.$
\end{Lem}

\begin{pf} We only treat the nontrivial cases for statement (1).
In case $i=j+1$, we have
$$T_j\cdots T_1(E_{j+1})=T_j(E_{j+1})=[E_{j+1}, E_{j}]_{v};$$
in case $2\leq i\leq j$,
$$T_j\cdots T_1(E_{i})=T_j\cdots T_iT_{i-1}(E_i)=T_j\cdots T_{i+1}(E_{i-1})=E_{i-1};$$
in case $i=1$,
$$\aligned
T_j\cdots T_1(E_{1})&=T_j\cdots T_2(-K_1F_1)=-T_j\cdots T_{2}(K_1)T_j\cdots T_{2}(F_1)\\
&=-K_1\cdots K_j[F_j, \cdots, F_1]_{v^{-1}}.
\endaligned$$

\end{pf}

\begin{Lem}\label{psi for linear case}
Set $\psi_{j}=T_1\cdots T_{j-1}T_{j}T_{j-1}\cdots T_{1}$. Then
    \begin{itemize}
    \item[(1)]     $\psi_{j}(E_i)=
         \left\{\begin{array}{lll} -K_2\cdots K_j[F_j, \cdots, F_2]_{v^{-1}}, && i=1;\\
                     -K_1\cdots K_{j-1}[F_1, \cdots, F_{j-1}]_{v^{-1}}, && i=j;\\
              \ [E_{j+1}, \cdots, E_{1}]_{v}, && i=j+1;\\
              E_{i}, && \text{otherwise};
             \end{array} \right.$ \\
            \item[(2)]
             $\psi_{j}(F_i)=
         \left\{\begin{array}{lll}
             -[E_2, \cdots, E_j]_{v}K^{-1}_2\cdots K^{-1}_j, && i=1;\\
             -[E_{j-1}, \cdots, E_1]_{v}K^{-1}_1\cdots K^{-1}_{j-1}, && i=j;\\
            \ [F_{1}, \cdots, F_{j+1}]_{v^{-1}}, && i=j+1;\\
              F_i, && \text{otherwise};
             \end{array} \right.$\\
       \end{itemize}
\end{Lem}

\begin{pf}(1) The proof is based on three lemmas above. Again,
we only consider the nontrivial cases. For the case $i=1$, we have
$$\aligned
\psi_{j}(E_1)&=T_{1}\cdots T_{j-1}(-K_1\cdots K_j[F_j, \cdots, F_1]_{v^{-1}})\\
&=-T_1(K_1\cdots K_j)T_1([F_j, \cdots, F_1]_{v^{-1}})\\
&=-K_2\cdots K_j[F_j, \cdots, F_2]_{v^{-1}};
\endaligned
$$
for the case $i=j$, $$\psi_{j}(E_j)=T_{1}\cdots T_{j-1}(E_{j-1})=-K_1\cdots K_{j-1}[F_1, \cdots, F_{j-1}]_{v^{-1}};$$
for the case $i=j+1$, $$\psi_{j}(E_{j+1})=T_{1}\cdots T_{j-1}([E_{j+1},E_j]_v)
=[E_{j+1},T_{1}\cdots T_{j-1}(E_j)]_v  =[E_{j+1}, \cdots, E_{1}]_{v};$$
 and for the case $2\leq i\leq j-1$, $$\psi_{j}(E_{i})=T_{1}\cdots T_{j-1}(E_{i-1})= E_{i}.$$
(2) The proof is similar to that of (1).
\end{pf}

\end{document}